\documentclass[10pt,oneside]{article}
\usepackage{times,amsmath,amsbsy,amssymb,amscd,mathrsfs,graphicx}
\usepackage{amssymb}
\usepackage{extarrows}
\usepackage{amsmath,bm}
\usepackage{amsmath,amsthm}

\usepackage{xcolor,color}
\usepackage[font=small,labelfont=md,textfont=it]{caption}
\usepackage{floatrow,float}
\usepackage[titletoc, title]{appendix}
\usepackage[colorlinks,linkcolor=blue,citecolor=blue]{hyperref}
\usepackage{longtable}
\usepackage{pgfplots}
\usepackage{diagbox}
\usepackage{booktabs,makecell,multirow}
\usepackage[capitalise, nosort]{cleveref}
\usepackage{algorithm}
\usepackage{algorithmic}
\usepackage{amsmath,bm}
\usepackage{cases}
\usepackage{geometry}
\usepackage{extarrows}
\usepackage{tcolorbox}
\usepackage{bbm}
\usepackage{syntonly}
\usepackage[figuresright]{rotating}
\usepackage{epstopdf}

\crefname{equation}{}{}
\crefname{lem}{Lemma}{Lemmas}
\crefname{Def}{Definition}{Definitions}
\crefname{thm}{Theorem}{Theorems}
\crefname{assum}{Assumption}{Assumptions}

\newcommand{\sgf}{{\mu_{\!f}}}
\newcommand{\Lf}{{L_{\!f}}}
\newcommand{\sgg}{{\mu_{g}}}

\newcommand{\dd}{\,{\rm d}}
\newcommand{\R}{\,{\mathbb R}}

\newcommand{\dual}[1]{\left\langle {#1} \right\rangle}
\newcommand{\prox}[0]{ {\bf prox}}

\newcommand{\dom}[0]{ {\bf dom\,}}
\newcommand{\argmin}[0]{ {\mathop{{\rm  argmin}}\,}}

\newcommand{\st}[0]{ {{\rm  s.t.}\,}}
\newcommand{\argmax}[0]{ {\mathop{{\rm  argmax}}\,}}
\newcommand{\inner}[1]{\left({#1} \right)}

\newcommand{\nm}[1]{\left\lVert {#1} \right\rVert}
\newcommand{\opnm}[1]{\left\lVert {#1} \right\rVert_{\rm op}}
\newcommand{\snm}[1]{\left\lvert {#1} \right\rvert}

\newcommand{\ssnm}[1]
{
	\left\vert\kern-0.25ex
	\left\vert\kern-0.25ex
	\left\vert
	{#1}
	\right\vert\kern-0.25ex
	\right\vert\kern-0.25ex
	\right\vert
}

\makeatletter
\def\spher@harm#1{%
	\vbox{\hbox{%
			\offinterlineskip
			\valign{&\hb@xt@2\p@{\hss$##$\hss}\vskip.2ex\cr#1\crcr}%
		}\vskip-.36ex}%
}
\def\gshone{\spher@harm{.}}
\def\gshtwo{\spher@harm{.&.}}
\def\gshthree{\spher@harm{.&.&.}}
\let\gsh\spher@harm
\makeatother

\newtheorem{prop}{Proposition}[section]
\newtheorem{Def}{Definition}[section]
\newtheorem{assum}{Assumption}

\newtheorem{lem}{Lemma}[section]

\newtheorem{rem}{Remark}[section]

\newtheorem{thm}{Theorem}[section]

\newcolumntype{I}{!{\vrule width 1,5pt}}
\newlength\savedwidth

\newlength\savewidth

\newcounter{mnote}
\let\oldmarginpar\marginpar
\renewcommand\marginpar[1]
{\-\oldmarginpar[\raggedleft\footnotesize #1]{\raggedright\footnotesize #1}}

\makeatletter\def\@captype{table}\makeatother

\newfloatcommand{capbtabbox}{table}[][\FBwidth]
\usepackage[T1]{fontenc} 
\usepackage[utf8]{inputenc} 
\usepackage{authblk}
\begin{document}
\title{
  \Large \bf Accelerated Primal-Dual Proximal Gradient Splitting Methods for Convex-Concave Saddle-Point Problems\thanks{This work was supported by the National Natural Science Foundation of China (Grant No. 12401402), the Science and Technology Research Program of Chongqing Municipal Education Commission (Grant Nos. KJZD-K202300505), the Natural Science Foundation of Chongqing (Grant No. CSTB2024NSCQ-MSX0329) and the Foundation of Chongqing Normal University (Grant No. 22xwB020).}
}
\author[,1,2]{Hao Luo\thanks{Email: luohao@cqnu.edu.cn}}
\affil[1]{National Center for Applied Mathematics in Chongqing, Chongqing Normal University, Chongqing, 401331, China}
\affil[2]{Chongqing Research Institute of Big Data, Peking University,  Chongqing, 401121, China}
\date{\today}
\maketitle

\begin{abstract}
	In this paper, based a novel primal-dual dynamical model with adaptive scaling parameters and Bregman divergences, we propose new accelerated primal-dual proximal gradient splitting methods for solving bilinear saddle-point problems with optimal nonergodic convergence rates. For the first, using the spectral analysis, we show that a naive extension of acceleration to a quadratic game is unstable. Motivated by this, we present an accelerated primal-dual gradient flow which combines acceleration with careful velocity correction. To work with non-Euclidean distances, we also equip our continuous model with general Bregman divergences and prove the exponential decay of a Lyapunov function. Then, new primal-dual splitting methods are developed based on proper semi-implicit Euler schemes of the continuous model, and the theoretical convergence rates are nonergodic and  optimal with respect to the matrix norms, Lipschitz constants and convexity parameters. Moreover, we introduce efficient restart variants to further improve the proposed methods and provide some numerical results to validate the practical performance.
\end{abstract}

\tableofcontents

\section{Introduction}
\label{sec:intro-H-APDG}
Let $\mathcal X$ and $\mathcal Y$ be finite dimensional  real vector spaces equipped with corresponding inner products $\inner{\cdot,\cdot}_\mathcal X$ and  $\inner{\cdot,\cdot}_\mathcal Y$ and  norms $\nm{\cdot}_\mathcal X$ and $\nm{\cdot}_\mathcal Y$.  Denote by $\mathcal X^*$ and $\mathcal Y^*$ the dual spaces containing all continuous linear functionals on $\mathcal X$ and $\mathcal Y$ respectively. The duality pairs between $\mathcal X^*$ and $\mathcal X$ and $\mathcal Y^*$ and $\mathcal Y$ are defined by $\dual{\cdot,\cdot}_\mathcal X$ and $\dual{\cdot,\cdot}_\mathcal Y$,  and the norms of the dual spaces are $\nm{\cdot}_{\mathcal X^*}:=\sup_{\nm{x}_{\mathcal X}\leq 1}\dual{\cdot,x}_{\mathcal X}$ and $\nm{\cdot}_{\mathcal Y^*}:=\sup_{\nm{y}_{\mathcal Y}\leq 1}\dual{\cdot,y}_\mathcal Y$. For simplicity, whenever no confusion arises, we omit the subscripts for inner products, norms, and duality pairings, and use $\nm{\cdot}_*$ to denote the dual norms. Consider the following convex-concave saddle-point problem
\begin{equation}\label{eq:minmax-L}
	\min_{x\in \mathcal X}\max_{y\in \mathcal Y} \,\mathcal L(x,y) =f(x)+ \dual{Ax,y}-g(y), 
\end{equation}
where  $f: \mathcal X\to\R\cup\{+\infty\}$ and $g:\mathcal Y\to\R\cup\{+\infty\}$ are properly closed convex functions and $A:\mathcal X\to \mathcal Y^*$ is a linear operator with the operator norm $\opnm{A}: = \sup_{\nm{x}\leq 1}\nm{Ax}_*$. The model problem \cref{eq:minmax-L} comes from many applications such as image processing \cite{chambolle_introduction_2016}, numerical partial differential equations \cite{benzi_numerical_2005}, machine learning \cite{Lan2020}, and optimal transport \cite{benamou_optimal_2021}. We are mainly interested in first-order methods that involve the (proximal) gradient information of the objective $f/g$; see the {\it proximal algorithm class} \cite[Definition 2.1]{zhang_lower_2022} for instance.

The saddle-point problem \cref{eq:minmax-L} is closely related to many structured convex optimization problems. It contains the unconstrained composite convex minimization
\begin{equation}\label{eq:primal}
	\min_{x\in \mathcal X}\,P(x): = f(x) + g^*(Ax),
\end{equation}
for which the dual problem reads equivalently as 
\begin{equation}\label{eq:dual}
	\min_{y\in \mathcal Y}\,D(y): =  g(y)+f^*(-A^* y),
\end{equation}
where $A^*:\mathcal Y\to \mathcal X^*$ denotes the adjoint operator of $A$ such that $\dual{Ax,y} = \dual{A^*y,x}$ for all $x\in \mathcal X$ and $y\in \mathcal Y$. When $g^*$ (or $f^*$) is the indicator function of a simple convex set $\mathcal K\subset \mathcal Y^*$, the primal problem \cref{eq:primal} (or the dual problem \cref{eq:dual}) also corresponds to the standard affinely constrained problem
\begin{equation}\label{eq:min-eq-1b}
	\min_{x\in \mathcal X}\, f(x)\quad \st \,Ax\in \mathcal K.
\end{equation}
Besides, if we replace $Ax$ (or $A^* y$) with an auxiliary variable, then the unconstrained minimization \cref{eq:primal} (or \cref{eq:dual}) becomes a special case of the two-block separable problem
\begin{equation}\label{eq:min-eq-2b}
	\min_{x_1\in \mathcal X_1,\,x_2\in \mathcal X_2}\, f_1(x_1)+f_2(x_2)\quad \st \,A_1x_1 +A_2x_2 = b\in Z,
\end{equation}
where $\mathcal X_1,\,\mathcal X_2$ and $Z$ are finite dimensional vector spaces; $A_1:\mathcal X_1\to Z$ and $A_2:\mathcal X_2\to Z$ are  linear operators; and $f_1: \mathcal X_1\to\R\cup\{+\infty\}$ and $f_2:\mathcal X_2\to\R\cup\{+\infty\}$ are properly closed convex functions. Note that the formulation \cref{eq:min-eq-2b} can be recast into the saddle-point problem \cref{eq:minmax-L} by using the Lagrange function.
\subsection{Literature review}
As mentioned before,  the minimax formulation \cref{eq:minmax-L} is usually equivalent to the unconstrained composite problem \cref{eq:primal}.
In this regard, classical first-order methods include proximal gradient \cite{parikh_proximal_2014},  forward-backward splitting \cite{Lions1979}, Douglas--Rachford splitting (DRS) \cite{douglas_numerical_1956}, 
and accelerated gradient methods \cite{beck_fast_2009,chen_first_2019,luo_accelerated_2021,luo_differential_2022,nesterov_gradient_2013,tseng_on_accelerated_Seattle_2008}. As for linearly constrained problems \cref{eq:min-eq-1b,eq:min-eq-2b}, prevailing algorithms are quadratic penalty methods \cite{li_convergence_2017,tran-dinh_proximal_2019}, Bregman iteration \cite{cai_linearized_2009}, augmented Lagrangian methods (ALM) \cite{Bai2024,Luo2026b,luo_primal-dual_2022,luo_universal_2024} and alternating direction methods of multipliers (ADMM) \cite{gabay_dual_1976,goldstein_fast_2014,Luo2025,ouyang_accelerated_2015,Sun2024,Xu2017,Yin2024,zhang_faster_2022}. 

\subsubsection{Primal-dual splitting}
In the following, we mainly review existing primal-dual splitting methods that applied directly to the saddle-point problem \cref{eq:minmax-L}. Zhu and Chan \cite{zhu_efficient_2008} proposed the primal-dual hybrid gradient (PDHG) method:
\begin{equation}\label{eq:pdhg}
	\left\{
	\begin{aligned}
		x_{k+1} = {}&\mathop{\argmin}\limits_{x\in \mathcal X}  \mathcal L(x,y_{k}) + \frac{1}{2\tau}\nm{x-x_k}^2 ,\quad\tau>0,\\
		y_{k+1} ={}& \mathop{\argmax}\limits_{y\in \mathcal Y}  \mathcal L(x_{k+1},y) -\frac{1}{2\theta}\nm{y-y_k}^2,\quad\theta>0,
	\end{aligned}
	\right.
\end{equation} 
which actually corresponds to the well-known Arrow--Hurwicz method \cite{arrow_studies_1958}. However, even for extremely small step lengths, PDHG is not necessarily convergent unless one of the objectives is strongly convex; see the counterexamples in \cite{he_convergence_2022}.
Later on, Chambolle and Pock \cite{chambolle_first-order_2011} presented a generalized version (CP for short)
\begin{equation}\label{eq:cp}
	\left\{
	\begin{aligned}
		x_{k+1} = {}&\mathop{\argmin}\limits_{x\in \mathcal X}  \mathcal L(x,y_{k}) + \frac{1}{2\tau}\nm{x-x_k}^2 , \\
		\bar x_{k+1} = {}&x_{k+1}+\alpha(x_{k+1}-x_k),\\
		y_{k+1} ={}& \mathop{\argmax}\limits_{y\in \mathcal Y}  \mathcal L(\bar x_{k+1},y) -\frac{1}{2\theta}\nm{y-y_k}^2 , 
	\end{aligned}
	\right.
\end{equation}
where $\alpha\in(0,1]$ denotes the relaxation parameter. 
When $\alpha=0$, CP amounts to PDHG, and it is also closely related to many existing methods such as  the extra-gradient method \cite{korpelevic_extra_1976}, DRS  \cite{oconnor_equivalence_2018}, and the preconditioned ADMM \cite{shefi_rate_2014,tian_alternating_2018}. In \cite{he_convergence_2012}, He and Yuan discovered that CP can be reformulated as a tight proximal point algorithm (PPA), and the relaxation parameter $\alpha\in(0,1]$ can be further extended to $\alpha\in[-1,1)$, with proper correction steps (see \cite[Algorithms 1 and 2]{he_convergence_2012}). 
However, these extra steps involve additional matrix-vector multiplications. 
For more extensions of CP, with preconditioning, extrapolation/correction and inexact variants, we refer to  \cite{chan_inertial_2014,he_algorithmic_2017,jiang_approximate_2021,liu_acceleration_2021,Zhu2023a}. As discussed in \cite[Remark 2.3]{valkonen_testing_2020}, it is possible to extend the PPA framework of CP to the non-Euclidean setting with general Bregman distances. This is in line with the spirit of the mirror-descent method by Nemirovski \cite{nemirovski_prox-method_2004}; see \cite{chambolle_ergodic_2016,chen_optimal_2014,darbon_accelerated_2022,Lan2020,tran-dinh_smooth_2018} and the references therein.

\renewcommand\arraystretch{1.8}
\begin{table}[h]
	\centering
	\small\setlength{\tabcolsep}{3.5pt}
	\begin{tabular}{c|c|c|c}
		\hline
		Assumption&$\sgf = \sgg=0$	&$\sgf + \sgg>0$	&$\sgf\sgg>0$\\
		\hline
		Rate&$\mathcal O\left(\frac{1}{k}\right)$&$\mathcal O\left(\frac{1}{k^2}\right)$&$\mathcal O\left(\exp\left(-\frac{\sqrt{\mu_f\mu_g}}{\opnm{A}}k\right)\right)$\\
		\hline
		\multirow{4}{*}{		Refs.}&\cite{Bai2023,chambolle_first-order_2011,chambolle_ergodic_2016,davis_convergence_2015,jiang_approximate_2021}&  \cite{chambolle_ergodic_2016,jiang_approximate_2021,malitsky_first-order_2018}& \cite{chambolle_ergodic_2016,jiang_approximate_2021} \\
		&\cite{Liu2025,malitsky_first-order_2018,tran-dinh_unified_2021,tran-dinh_non-stationary_2020,xie_accelerated_2019}&  \cite{Liu2025,tran-dinh_unified_2021,tran-dinh_non-stationary_2020,xie_accelerated_2019}& \cite{rasch_inexact_2020,tran-dinh_unified_2021} \\
		&  \cite{tran-dinh_adaptive_2017,tran-dinh_smooth_2018,valkonen_inertial_2020}*&\cite{chambolle_first-order_2011,He2025,valkonen_inertial_2020}*&\cite{chambolle_first-order_2011,valkonen_inertial_2020}*\\
		&	\textbf{This work}*&	\textbf{This work}*&		\textbf{This work}* \\
		\hline
	\end{tabular}
	\caption{Summary of convergence rates of existing primal-dual methods for solving \cref{eq:minmax-L}. Here, $\mu_f,\,\mu_g\geq 0$ are respectively the strong convexity constants of $f$  and $g$. Nonergodic results are marked with ``*''. }
	\label{tab:rate-pdhg}
\end{table}

\begin{figure}[H]
	\centering 
	\includegraphics[width=13cm]{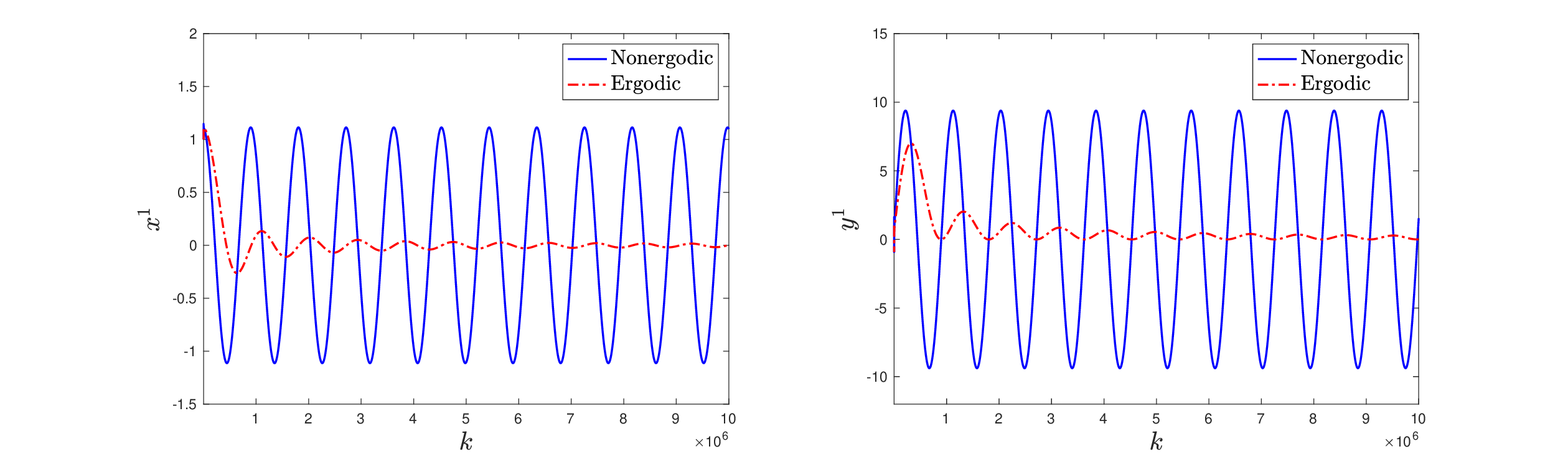}
	\caption{Iterative behavior of CP for solving the matrix game $\min_{x\in \R^3}\max_{y\in  \R^3}\,y^\top Ax$, where $A\in\R^{3\times 3}$ is given by \cref{eq:counter-A} with $\epsilon=10^{-4}$ in \cref{sec:counter}. The initial points are $x_0=y_0=[1,1,1]^\top$ and the parameters are $\alpha=1,\,\theta=1,\,\tau = 0.5/\nm{A}^2$. The left part is the first components of the nonergodic sequence $\{x_k\}$ and the ergodic sequence $\{X_k\}$, and the right part is the first components of the nonergodic sequence $\{y_k\}$ and the ergodic sequence $\{Y_k\}$.}
	\label{fig:cp-limit-cycle}
\end{figure}

According to \cite{chambolle_first-order_2011,chambolle_ergodic_2016}, for generic convex-concave problems, CP converges with an {\it ergodic} sublinear rate $\mathcal O(1/k)$ for the partial primal-dual gap, and fast (ergodic) sublinear/linear rates can also be achieved for (partially) strongly convex case $\sgf+\sgg>0$; see \cref{tab:rate-pdhg} for a summary of related works on convergence rates. However, from the literature, it is rare to see {\it nonergodic} results, which, as mentioned by Chambolle and Pock in \cite[Section 8]{chambolle_ergodic_2016}, converge very often much faster than the sense of ergodic. We mention that linear convergence can still be obtained even for non-strongly convex problems \cite{du_linear_2019,kovalev_accelerated_2022,luke_globally_2018} and nonergodic convergence with respect to  fixed-point type residual can be found in \cite{chan_inertial_2014,he_algorithmic_2017,kim_accelerated_2021}. However, it seems not easy to obtain nonergodic convergence from ergodic results. Additionally, in \cref{fig:cp-limit-cycle}, we give a simple example to show that the nonergodic sequence (last iterate) of CP can be dramatically slower than the ergodic sequence (averaged iterate); see \cref{sec:counter} for more details. This agrees well with the theoretical validation in \cite{Golowich2020}.

Based on the PPA structure of CP and the {\it testing} approach \cite{valkonen_testing_2020}, Valkonen \cite{valkonen_inertial_2020} proposed the  inertial corrected primal-dual proximal splitting (IC-PDPS), which achieves the optimal sublinear/linear nonergodic rates (cf.\cref{tab:rate-pdhg}). In \cite{Luo2025y}, we provided a continuous-time perspective on IC-PDPS, which can be viewed as an implicit-explicit time discretization of a novel primal-dual dynamics. 
When the primal (dual) objective $f$ ($g$) in \cref{eq:minmax-L} possesses the composite structure $f = f_1+f_2$ ($g = g_1+g_2$), where the smooth part $f_1$($g_1$) has Lipschitz continuous gradient, some extensions  can be found in \cite{chambolle_ergodic_2016,chen_primal-dual_2016,condat_primaldual_2013,Condat2026,latafat_asymmetric_2017,Salim2022,vu_splitting_2013,Yan2018}. There also exist primal-dual splitting methods \cite{Chen2023,chen_optimal_2014,Driggs2024,he_accelerated_2016,kovalev_accelerated_2022,xie_accelerated_2019} possessing optimal mixed-type complexity bounds; see \cref{tab:mix-rate-pdhg}. Compared with these works, our methods cover both convex ($\mu_f=\mu_g=0$) case and partially convex case ($\mu_f+\mu_g>0$), and the use of Bregman divergences provides more flexible capability for both theoretical analysis and practical computations.

\renewcommand\arraystretch{1.8}
\begin{table}[H]
	\centering
	\small\setlength{\tabcolsep}{5pt}
		\begin{tabular}{c|cc|c}
			\hline
			&	\multicolumn{2}{c|}{$\sgf=0$}	&$\sgf>0$\\
			\hline
			\multirow{3}{*}{$\sgg=0$}	&\multirow{3}{*}{	$\mathcal O\left(\frac{\opnm{A}}{\epsilon}+\sqrt{\frac{\Lf}{\epsilon}}\right)$} &\cite{chen_optimal_2014,Driggs2024,nesterov_excessive_2005,nesterov_smooth_2005}&
			\multirow{2}{*}{	$\mathcal O\left(\frac{\opnm{A}}{\sqrt{\sgf\epsilon}}+\sqrt{\frac{L_f}{\mu_f}}\snm{\ln\epsilon}\right)$}\\
			& &  \cite{he_accelerated_2016,xie_accelerated_2019}$^\star$&
			\\
			& &  \textbf{This work}&\cite{Driggs2024}\,\textbf{This work}
			\\
			\hline 
			\multirow{2}{*}{	$\sgg>0$}&				\multirow{2}{*}{$\mathcal O\left(\frac{\opnm{A}}{\sqrt{\sgg\epsilon}}+\sqrt{\frac{\Lf}{\epsilon}}\right)$}&\cite{nesterov_excessive_2005} \cite{xie_accelerated_2019}$^\star$&			$\mathcal O\left(\sqrt{\frac{\opnm{A}^2}{\mu_f\mu_g}+\frac{L_f}{\mu_f}}\snm{\ln\epsilon}\right)$
			\\
			& &  \textbf{This work}& \cite{Chen2023,Driggs2024,kovalev_accelerated_2022}\,\,\textbf{This work}
			\\
			\hline
		\end{tabular}
		\caption{Summary of existing mixed-type complexity bounds for solving \cref{eq:minmax-L} under the composite structure $f = f_1+f_2$, where $f_1$ is $\sgf$-convex with  $\Lf$-Lipschitz continuous gradient. Ergodic results are marked with ``$\star$''.}
		\label{tab:mix-rate-pdhg}
\end{table}

\subsubsection{Smoothing}
Except primal-dual splitting methods based on the PPA framework mentioned above, let us also review other works with smoothing technique. For nonsmooth composite convex optimization problems with the minimax structure (like \cref{eq:minmax-L}), Nesterov \cite{nesterov_smooth_2005} combined his accelerated gradient \cite{Nesterov1983} with smoothing and improved the iteration complexity of subgradient-type methods from $\mathcal O(1/\epsilon^2)$ \cite{nemirovsky_problem_1983} to $\mathcal O(1/\epsilon)$. Afterwards, there are primal-dual methods based on the excessive gap technique (double smoothing approach) \cite{devolder_double_2012,nesterov_excessive_2005} and homotopy strategy (adaptively changing smoothing parameters) \cite{tran-dinh_adaptive_2017}.

In \cite{tran-dinh_smooth_2018}, Tran-Dinh {\it et al.} proposed two accelerated smoothed gap reduction (ASGARD) algorithms and established the optimal sublinear rate $\mathcal O(1/k)$ for the primal objective residual.  The unified variant of ASGARD in \cite{tran-dinh_unified_2021} and two primal-dual algorithms in \cite{tran-dinh_non-stationary_2020} achieve optimal sublinear/linear rates for the primal-dual gap, the primal objective residual, and the dual objective residual. However, it should be noticed that (i) both \cite[Algorithm 2]{tran-dinh_smooth_2018} and \cite[Algorithm 2]{tran-dinh_non-stationary_2020} involve {\it three} proximal operations per iteration; (ii) the convergence results in \cite{tran-dinh_unified_2021,tran-dinh_non-stationary_2020} are {\it semi-ergodic}. In other words, the primal sequence is in nonergodic sense while the dual sequence is averaging; (iii) the extra dual steps in \cite[Algorithms 1 and 2]{tran-dinh_non-stationary_2020} require additional matrix-vector multiplications. 


\subsection{Main contribution}
In this work, combining acceleration with velocity correction, we propose a novel accelerated primal-dual gradient (APDG) flow \cref{eq:APDG-sys}. Based on proper time discretizations, we obtain new primal-dual splitting methods for solving the saddle-point problem \cref{eq:minmax-L}. More specifically, we highlight our main contributions as follows.
\begin{itemize}
	\item[(i)] By using the spectral analysis, we prove that a naive extension of acceleration model to a quadratic game does not work. With careful velocity correction, we propose a novel accelerated primal-dual dynamical system \cref{eq:APDG-spec} and prove the improved spectral radials of the Gauss-Seidel splitting ; see \cref{thm:spec-gs}.
	\item[(ii)] We then equip the template model \cref{eq:APDG-spec}  with adaptively scaling parameters and general Bregman divergence (cf.\cref{eq:APDG-sys,eq:ABPDG-sys}). This allows us to handle both convex and (partially) strongly convex cases in a unified way, and make the dynamical model work for non-Euclidean distances. Also, we establish the exponential decay rate of a tailored Lyapunov function, which provides important guidance for the discrete analysis; see \cref{thm:conv-APDG,thm:conv-ABPDG}.
	\item[(iii)] After that, new primal-dual splitting methods are developed based on proper semi-implicit Euler schemes of the continuous model \cref{eq:ABPDG-sys}.  Our methods (\cref{algo:ABPD-PS,algo:ABPD-PGS}) belong to the first-order (proximal) algorithm classes (see \cite[Section 2]{zhang_lower_2022}), which involve the matrix-vector multiplications of $A$ and $A^*$ and the proximal calculations of $f$ and $g$ only {\it once} in each iteration. In addition, the theoretical rates are {\it nonergodic} and {\it optimal} with respect to the matrix norm $\opnm{A}$,\, the Lipschitz constant $\Lf$ and the convexity parameter $\sgf(\sgg)$, as mentioned in \cref{tab:rate-pdhg,tab:mix-rate-pdhg}. For instance, to achieve the accuracy $	\mathcal L(x_k,\widehat{y})-	\mathcal L(\widehat{x},y_k)\leq \epsilon$, the iteration complexity bound of \cref{algo:ABPD-PS} is (cf.\cref{thm:rate-ABPD-PS})
	\[
	\mathcal O	\left(\min\left\{\frac{\opnm{A}}{\epsilon},\,\frac{\opnm{A}}{\sqrt{\sgf \epsilon}},\,\frac{\opnm{A}}{\sqrt{\sgg \epsilon}},\,\frac{\opnm{A}}{\sqrt{\mu_f\mu_g}}\snm{\ln \epsilon}\right\}\right),
	\]
	and for the composite objective $f = f_1+f_2$, we have (cf.\cref{thm:rate-sym-ABPD-PGS})
	\[
	\begin{aligned}
		{}&	\mathcal O	\Bigg(\min\left\{\frac{\opnm{A}}{\epsilon},\,\frac{\opnm{A}}{\sqrt{\sgf  \epsilon}},\,\frac{\opnm{A}}{\sqrt{\sgg  \epsilon}},\,\frac{\opnm{A}}{\sqrt{ \sgf\sgg }}\snm{\ln \epsilon}\right\}\\
		{}&\qquad\qquad +\min\left\{\sqrt{\frac{L_f }{\epsilon}},\,\sqrt{\frac{L_f}{\mu_f}}\snm{\ln \epsilon}\right\}+\min\left\{\sqrt{\frac{L_g }{\epsilon}},\,\sqrt{\frac{L_g}{\mu_g}}\snm{\ln \epsilon}\right\}\Bigg).
	\end{aligned}
	\]
	\item[(iv)] Thanks to the primal and dual scaling parameters, both the algorithm designing and convergence analysis cover automatically the convex and (partially) strongly convex objectives. The use of Bregman divergences not only unifies the analysis of standard Euclidean distances and general cases in an elegant way, but also makes our methods more flexible and adaptive to problem-dependent metrics that might better than the Euclidean one (like the entropy). Moreover, we propose an efficient fixed restart strategy for our methods (\cref{algo:ABPD-PS,algo:ABPD-PGS}) and establish the linear convergence rates (cf.\cref{thm:conv-ABPD-PS-R,thm:conv-ABPD-PGS-R}).
\end{itemize}
\subsection{Organization}
The rest of this paper is organized as follows. In \cref{sec:pre}, we prepare essential preliminaries including basic notations, the definition of Bregman distance, two auxiliary functions and some useful Lyapunov functions. Next, in \cref{sec:spec}, we study a simple quadratic game via spectral analysis and find some dynamical template system with provable stability and acceleration. This leads to two novel primal-dual flow models in \cref{sec:APDG} with adaptive scaling parameters and Bregman divergences. After that, in \cref{sec:ABPD-PS,sec:ABPD-PGS}, we propose two new primal-dual splitting algorithms from numerical discretizations of the continuous model and prove the nonergodic mixed-type convergence rate via a tailored discrete Lyapunov function. Also, we introduce a fixed restart idea in \cref{sec:restart} and prove the linear convergence, and several numerical examples are reported in \cref{sec:numer}. Finally, some concluding remarks are given in \cref{sec:conclu}.

\section{Preliminaries}
\label{sec:pre}
\subsection{Notations}
For a proper, closed and convex function $h:\mathcal X\to\R\cup\{+\infty\}$,
denote by $\dom h $ and $\partial h (x)$ its domain and its subdifferential at $x$, respectively. The standard $\mu$-convexity means that there exists $\mu\geq0$ such that
\begin{equation}\label{eq:mu-ineq}
	h(v) 			\geq{} h( x) 
	+	\dual{p,v- x}  + \frac{\mu}{2}\nm{ x-v}^2
	\quad  \forall\,x,v \in \dom h,\quad p\in \partial h(v).
\end{equation}
We say $h :\dom h \to\R$ is $C^1$-smooth if it is continuous differentiable, and it is called $C_L^1$-smooth when its gradient $\nabla h $ is $L$-Lipschitz continuous:
\[
\nm{\nabla h(x)-\nabla h(v)}_*\leq L\nm{x-v}\quad\forall\,x,v\in \dom h,
\]
which implies that
\begin{equation}\label{eq:Lip-ineq}
	h(v) 			\leq{} h( x) 
	+	\dual{\nabla h(x),v- x}  + \frac{L}{2}\nm{ x-v}^2
	\quad\forall\,x,v\in \dom h.
\end{equation}
Given any $R:\mathcal X\to \mathcal X^*$, if $\dual{Rx,x}\geq 0$ for all $x\in \mathcal X$ and $R=R^*$, then we say $R$ is symmetric positive semi-definite (SPSD), and the induced  semi-norm is defined by $\nm{x}_R := \sqrt{\dual{Rx,x}}$. If $R$ is SPSD and $\dual{Rx,x}> 0$ for all $x\in \mathcal X/\{0\}$, then it is also called symmetric positive definite (SPD) and the induced semi-norm becomes a norm. In particular, if $R=I:\mathcal X\to \mathcal X^*$ is the natural identity mapping, then the induced norm is $\nm{\cdot}_I:=\sqrt{\inner{\cdot,\cdot}}$. For notation simplicity, $\nm{\cdot}$ also stands for the 2-norm of matrices and vectors whenever there is no ambiguity. For any $M\in\R^{n\times n}$, we use $\sigma(M)$ to denote the set of all eigenvalues of $M$ and let $\rho(M):=\max_{\lambda\in\sigma(M)}\snm{\lambda}$ be the spectral radius of $M$. If $M\in\R^{n\times n}$ is SPSD, then the spectral condition number is defined by $\kappa(M) := \lambda_{\max}(M)/\lambda_{\min}(M)$, where  $\lambda_{\max}(M)$ and $\lambda_{\min}(M)$ are respectively the largest and smallest eigenvalues.
\subsection{Bregman distance}
Let us first recall the definitions of prox-function and Bregman distance.
\begin{Def}[Prox-function]\label{def:prox-fun}
	We call $\phi:\mathcal X\supset\dom \phi\to \R$ a prox-function  if it is $C^1$-smooth and 1-convex.
\end{Def}
%
\begin{Def}[Bregman distance]\label{def:bregman}
	Let $\phi:\mathcal X\supset\dom\phi\to\R$ be a prox-function.
	The Bregman distance induced by $\phi$ is $	D_\phi(x,v):=		\phi(x)-\phi(v)-\langle\nabla \phi(v), x-v\rangle$ for all $x,v\in\dom\phi$.
	\end{Def}
	
	Since $\phi$ is 1-convex, we have $D_\phi(x,v)\geq 1/2\nm{x-v}^2$ for all $x,v\in\dom\phi$. For $\phi(x)=1/2\nm{x}_R^2$  with an SPSD operator $R:\mathcal X\to \mathcal X^*$, we obtain $D_\phi(x,v) = 1/2\nm{x-v}_R^2$. In addition, we have the following three-term identity \cite[Lemma 2.1]{luo_universal_2024}.
	\begin{lem}[\cite{luo_universal_2024}]
		\label{lem:3-id}
		Let $\phi:\mathcal X\supset\dom\phi\to\R$ be a prox-function, then
		\begin{equation}\label{eq:3-id}
			\dual{\nabla \phi(x)-\nabla \phi(v), v-\bar x}=D_\phi(\bar x, x)-D_\phi(\bar x, v)-D_\phi(v, x) ,\quad\forall\, x,v,\bar x\in\dom\phi.
		\end{equation}
		In particular, for $\phi(x)=1/2\nm{x}_R^2$ with an SPSD operator $R:\mathcal X\to \mathcal X^*$, we have 
		\begin{equation}\label{eq:3-id-R}
			\dual{R(x-v), v-\bar{x}}=\frac{1}{2}\nm{x-\bar{x}}_R^2-\frac{1}{2}\nm{v-\bar{x}}_R^2-\frac{1}{2}\nm{x-v}_R^2.
		\end{equation}
	\end{lem}
	
	For convergence rate analysis, we introduce the conception of {\it relative convexity} \cite[Definitin 1.2]{lu_relatively_2018}. Clearly, since $\phi$ is 1-convex, the relative $\mu$-convexity implies the standard $\mu$-convexity (cf.\cref{eq:mu-ineq}). 
	\begin{Def}[Relatively convexity]\label{def:rela-strng-conv} 
		Let $\phi:\mathcal X\supset\dom\phi\to\R$ be a prox-function and $h:\mathcal X\to\R\cup\{+\infty\}$ a closed proper convex function such that $\dom h\subset\dom\phi$. 
		We say $h$ is relatively $\mu$-convex with respect to (w.r.t. for short) $\phi$ if there exists $\mu\geq 0$ such that
		\begin{align}
			\label{eq:scv-f}
			h(x)\geq{}  &h(v)+\dual{p,x-v}+\sgf D_\phi(x,v),\quad  \forall\,x,v \in \dom h,\quad p\in \partial h(v).
		\end{align}
	\end{Def}
	\begin{lem}\label{lem:L-sg-descent}
		Let $\phi:\mathcal X\supset\dom\phi\to\R$ be a prox-function and $h:\dom h\to \R$ a closed proper convex function such that $\dom h\subset\dom\phi$. If  $h$ is $C_{L}^1$-smooth and relatively $\mu$-convex w.r.t. $\phi$, then  
		\[
		-	\dual{\nabla h(\bar x),v-x}\leq h(x) - h(v) - \mu   D_\phi(x,\bar x) + \frac{L}{2}\nm{\bar x-v}^2,\quad\forall\, x,v,\bar x\in\dom h.
		\]
	\end{lem}
	\begin{proof}
		It follows from \cref{eq:Lip-ineq} and $			-\dual{\nabla h(\bar x),\bar x-x}\leq{} h(x) - h(\bar x) -  \mu   D_\phi(x,\bar x)$ immediately.
	\end{proof}
	\subsection{Two auxiliary functions}
	For later use, we introduce two auxiliary functions related to the min-part $f$ and the max-part $g$ of the saddle-point problem \cref{eq:minmax-L}. Denote by $\dom \mathcal L := \dom f\times \dom g$ the domain of the minimax function. Given $\widehat{\bm z}=(\widehat{x},\widehat{y})\in\dom\mathcal L$, define
	\begin{equation}\label{eq:hat-f-hat-g}
		\widehat{f}(x;\widehat{\bm z}): = {} f(x)+\dual{A^*\widehat{y}, x-\widehat{x}},\quad
		\widehat{g}(y;\widehat{\bm z}): ={}  g(y)-\dual{A \widehat{x},y-\widehat{y}},
	\end{equation}
	for all $(x,y)\in\dom \mathcal L$. It is clear that both $\widehat{f}(\cdot,\widehat{\bm z}):\mathcal X\to\R\cup\{+\infty\}$ and $\widehat{g}(\cdot,\widehat{\bm z}):\mathcal Y\to\R\cup\{+\infty\}$ are properly closed convex, and we have the relations: 
	\begin{equation}\label{eq:subg}
		\partial \widehat{f}(x;\widehat{\bm z})=\partial_x\mathcal L(x,\widehat{y}) ,\quad
		\partial \widehat{g}(y;\widehat{\bm z})=-\partial_y\mathcal L(\widehat{x},y).
	\end{equation}
	Since  $\widehat{f}(\widehat{x};\widehat{\bm z}) = f(\widehat{x})$ and $\widehat{g}(\widehat{y};\widehat{\bm z}) = g(\widehat{y})$, a simple calculation implies
	\begin{equation}\label{eq:gap-obj} 
		\begin{aligned}
			\mathcal L(x, \widehat{y})-\mathcal L(\widehat{x}, y) 
			=		{}&\widehat{f}(x;\widehat{\bm z})-f(\widehat{x}) + \widehat{g}(y;\widehat{\bm z})-g(\widehat{y}).
		\end{aligned}
	\end{equation} 
	
	Throughout we equip two prox-functions $\phi$ and $\psi$ respectively for the primal and dual variables, and make the following assumption: \begin{assum}\label{assum:ABPD-PS-prox-fun}
		$\phi:\mathcal X\supset \dom\phi\to\R$ and $\psi:\mathcal Y\supset \dom\psi\to\R$ are two prox-functions such that $\dom f\subset\dom \phi$ and $\dom g\subset\dom \psi$.  
	\end{assum}
	
	Invoking \cref{def:prox-fun,def:rela-strng-conv,eq:subg,eq:gap-obj}  gives the lemma below.
\begin{lem}
	\label{lem:L-gap-obj-res}
	Assume $(f,g)$ are relatively $(\mu_f,\mu_g)$-convex w.r.t. $(\phi,\psi)$ with $\sgf,\,\sgg\geq0$.
	\begin{itemize}
		\item[(i)]  For any $\widehat{\bm z}=(\widehat{x},\widehat{y})\in\dom \mathcal L$, $\widehat{f}(\cdot,\widehat{\bm z})$ and $\widehat{g}(\cdot,\widehat{\bm z})$ are also relatively convex w.r.t. $\phi$ and $\psi$ with constants $\sgf,\,\sgg\geq0$.
		\item[(ii)] If  $\widehat{\bm z}=(\widehat{x},\widehat{y})$ is a saddle point to \cref{eq:minmax-L}, then we have
		\[
		\begin{aligned}
			\widehat{f}(x;\widehat{\bm z}) - f(\widehat{x})\geq{}& \mu_fD_\phi(x,\widehat{x})\geq \sgf/2\nm{x-\widehat{x}}^2,\\ \widehat{g}(y;\widehat{\bm z})- g(\widehat{y})\geq{}& \mu_gD_{\psi}(y,\widehat{y})\geq \sgg/2\nm{y-\widehat{y}}^2,
		\end{aligned}
		\] 
		for all $(x,y)\in\dom \mathcal L$. In particular, by \cref{eq:gap-obj}, it holds that
		\[
		\mathcal L(x, \widehat{y})-\mathcal L(\widehat{x}, y) \geq\mu_fD_\phi(x,\widehat{x})+\mu_gD_{\psi}(y,\widehat{y})\geq  \frac{\sgf}{2}\nm{x-\widehat{x}}^2+\frac{\sgg}{2}\nm{y-\widehat{y}}^2.
		\]
	\end{itemize} 
\end{lem}
	\subsection{Lyapunov functions}
	To prove the convergence rates of our continuous models and discrete algorithms, we need some technical Lyapunov functions. 
	
	For simplicity, let $\mathcal Z: = (\dom f)^2\times (\dom g)^2$.
	Given  $\widehat{\bm z}=(\widehat{x},\widehat{y})\in\dom \mathcal L$, define a Lyapunov function $\mathcal E:\R_+^2\times \mathcal Z\to \R_+$ by that 
	\begin{equation}\label{eq:Et}
		\mathcal E(\Theta,Z;\widehat{\bm z}): = \mathcal L(x, \widehat{y})-\mathcal L(\widehat{x}, y)+\frac{\gamma}{2}\nm{v-\widehat{x}}^2
		+\frac{\beta}{2}\nm{ w- \widehat{y}}^2,
	\end{equation}
	where $\Theta=(\gamma,\beta)\in\R_+^2$ and $Z = (x,v,y,w)\in\mathcal Z$. To work with Bregman distances, we also introduce a generalized Lyapunov function $	\mathcal E_{ D}:\R_+^2\times \mathcal Z\to \R_+$:
	\begin{equation}\label{eq:Ft}
		\mathcal E_{ D}(\Theta,Z;\widehat{\bm z}): = \mathcal L(x, \widehat{y})-\mathcal L(\widehat{x}, y)+\gamma D_{\phi}(\widehat{x},v)
		+\beta D_\psi(\widehat{y},w).
	\end{equation}
	
	The previous two functions \cref{eq:Et,eq:Ft} are mainly applied to analyzing continuous models (cf.\cref{thm:conv-APDG,thm:conv-ABPDG}). For discrete algorithms, however, some additional careful correction is required. Let 
	$\{\alpha_k\}\subset\R_+,\,\{\Theta_k=(\gamma_k,\beta_k)\}\subset\R_+^2$  and $\{Z_k=(x_k,v_k,y_k,w_k)\}\subset \mathcal Z$ be given, define a discrete Lyapunov function 
	\begin{equation}\label{eq:Hk}
		\mathcal H_k(\widehat{\bm z}): = \mathcal E_{ D}(\Theta_k,Z_k;\widehat{\bm z})-\alpha_k\dual{A(v_k-\widehat{x}),w_k-\widehat{y}},
	\end{equation}
	with any $\widehat{\bm z} \in\dom \mathcal L$.
	For convergence analysis, the key is to establish the upper bound of
	\begin{equation}\label{eq:diff-Hk}
		\mathcal H_{k+1}(\widehat{\bm z})-\mathcal H_k(\widehat{\bm z})+\alpha_k\mathcal H_{k+1}(\widehat{\bm z}) = \mathbb I_{11}+ \mathbb I_{21}+\mathbb I_{22}+ \mathbb I_{31}+ \mathbb I_{32},
	\end{equation}
	where the decomposition on the right hand side consists of 
	\[
	\begin{aligned} 
		\label{eq:Hk-I1}
		\mathbb I_{11}:={}&\alpha_k\dual{A(v_k-\widehat{x}), w_k-\widehat{y}}-\alpha_{k+1}(1+\alpha_k)\dual{A(v_{k+1}-\widehat{x}), w_{k+1}-\widehat{y}}, \\
		\mathbb I_{21}: ={}&\widehat{f}(x_{k+1};\widehat{\bm z})-\widehat{f}(x_k;\widehat{\bm z})+ \alpha_k\big(\widehat{f}(x_{k+1};\widehat{\bm z})-f(\widehat{x})\big), \\ 
		\mathbb I_{22}:={}& 
		\gamma_{k+1}(1+\alpha_k)D_\phi(\widehat{x},v_{k+1})-	\gamma_{k}D_\phi(\widehat{x},v_{k}) ,\\
		\mathbb I_{31}: ={}&\widehat{g}(y_{k+1};\widehat{\bm z})-\widehat{g}(y_k;\widehat{\bm z})+ \alpha_k\left(\widehat{g}(y_{k+1};\widehat{\bm z})-g(\widehat{y})\right),\\ 
		\mathbb I_{32}:={}& 
		\beta_{k+1}(1+\alpha_k )D_\psi(\widehat{y},w_{k+1})-	\beta_{k}D_\psi(\widehat{y},w_{k}).
	\end{aligned}
	\]
	
	\begin{lem}\label{lem:decomp-diff-Hk}
		Let $\{\alpha_k\}\subset\R_+,\,\{\Theta_k=(\gamma_k,\beta_k)\}\subset\R_+^2$  and $\{Z_k=(x_k,v_k,y_k,w_k)\}\subset \mathcal Z$ be such that
		\begin{align}
			\label{eq:assum-Axv-D}
			&\alpha_k \snm{\dual{A(v_k-\widehat{x}), w_k-\widehat{y}}}\leq  \gamma_k D_\phi(\widehat{x},v_k)	+ \beta_k  D_\psi(\widehat{y}, w_k),\\
			\label{eq:diff-Hk-uni}
			& \mathcal H_{k+1}(\widehat{\bm z})-\mathcal H_k(\widehat{\bm z})+\alpha_k\mathcal H_{k+1}(\widehat{\bm z})  \leq 0, 
		\end{align}
		for all $\widehat{\bm z}=(\widehat{x},\widehat{y})\in\dom \mathcal L$,
		then we have
		\begin{equation}\label{eq:L-gap-x-y-uni}
			\mathcal L(x_k,\widehat{y})-	\mathcal L(\widehat{x},y_k)
			+	\frac{\sgf  }{2}\nm{x_k-\widehat{x}}^2+\frac{\sgg }{2}\nm{y_k-\widehat{y}}^2
			\leq 2\theta_k\mathcal H_0(\widehat{\bm z}),
		\end{equation} 
		where $\{\theta_k\}$ is defined by 
		\begin{equation}\label{eq:thetak}
			\theta_0=1,\quad\theta_k = \prod_{i=0}^{k-1} \frac{1}{1+\alpha_i}\quad\forall\,k\geq 1.
		\end{equation} 
	\end{lem} 
	\begin{proof}
		It is not hard to obtain from \cref{eq:thetak} that 
		\begin{equation}\label{eq:diff-tk}
			\theta_{k+1}-\theta_k = -\alpha_k\theta_{k+1}.
		\end{equation}  
		According to our assumption \cref{eq:diff-Hk-uni}, it follows that 
		\[
		\mathcal H_k(\widehat{\bm z}) + \theta_k\sum_{i=0}^{k-1}\frac{\Omega_i}{\theta_i}
		\leq \theta_k\mathcal H_0(\widehat{\bm z}),\quad k\geq 0.
		\]
		Then by \cref{eq:assum-Axv-D} we get
		$	\mathcal L(x_k,\widehat{y})-	\mathcal L(\widehat{x},y_k)\leq \mathcal H_k(\widehat{\bm z}) $. Consequently, invoking  \cref{lem:L-gap-obj-res} (ii) proves \cref{eq:L-gap-x-y-uni}.
	\end{proof} 
	
	\section{Spectral Analysis of a Quadratic Game}
	\label{sec:spec}
	In this section, focusing on the standard Euclidean setting, we follow the main idea from \cite[Section 2]{luo_differential_2022} to seek proper lifting system of a quadratic game. By using the tool of spectral analysis, we verify the stability and better dependence on the condition number. This provides the key guidance for desinging our continuous models for saddle-point problems.

	\subsection{The unconstrained case}
	\label{sec:spec-uncon} 
	Let us start from the quadratic programming
	\begin{equation}\label{eq:min-quad}
		\min_{x\in \R^n}\,f(x) = \frac{1}{2}\nm{x}_Q^2,
	\end{equation}
	where $Q\in\R^{n\times n}$ is SPD. Note that the optimal solution is $x^* = 0$. The classical gradient flow reads as $x' = -\nabla f(x) = -Qx$, 
	which is converges exponentially.
	The explicit Euler discretization, i.e., the gradient descent method $x_{k+1} = x_k - \alpha\nabla f(x_k) = (I-\alpha Q)x_k$, converges at a linear rate $\mathcal O(1-1/\kappa(Q))$ with $\alpha=O(1/\kappa(Q))$, which is suboptimal.
	
	To seek acceleration, the basic idea in \cite[Section 2]{luo_differential_2022} is to find some lifting vector field $\bm x' = \mathcal G\bm x$, where $\bm x=(x,v)$ and $\mathcal G$ is preferably a block matrix  with $\Re(\lambda)<0$ for all $\lambda\in \sigma(\mathcal G)$. The block matrix $\mathcal G$ transforms $\sigma(-Q)$ from the negative real line to the left half complex plane and hopefully provides better dependence on the condition number.
	
	Especially, a nice candidate has been given in \cite[Section 2]{luo_differential_2022}:
	\begin{equation}\label{eq:Gnag}
		\mathcal G = \begin{bmatrix}
			-I&I\\
			I-Q/\lambda_{\min}(Q)& -I
		\end{bmatrix},
	\end{equation}
	which leads to the continuous NAG flow model
	\begin{equation}\label{eq:nag}
		\begin{aligned}
			x' = {}v-x,\quad v' = {}x-v-\frac{Qx}{\lambda_{\min}(Q)}.
		\end{aligned}
	\end{equation}
	It is proved that \cite[Theorem 2.1]{luo_differential_2022}, a suitable Gauss--Seidel splitting of \cref{eq:nag} possesses a smaller spectral radius $\mathcal O(1-1/\sqrt{\kappa(Q)})$ with larger step size $\alpha=\mathcal O(1/\sqrt{\kappa(Q)})$, which matches the optimal rate \cite{nesterov_introductory_2004}. 
	\subsection{A quadratic game}
	Let us then consider the saddle-point problem \cref{eq:minmax-L} with $\mathcal X=\R^n,\,\mathcal Y=\R^m$ and $A\in\R^{m\times n}$. The objectives are simple quadratic functions $f(x)=\mu/2\nm{x}^2$ and $g(y) = \mu/2\nm{y}^2$, which lead to a quadratic game
	\begin{equation}\label{eq:quad-game}
		\min_{x\in \R^n}\max_{y\in \R^m}\frac{1}{2}(x^\top,-y^\top) \mathcal Q\begin{pmatrix}
			x\\ y
		\end{pmatrix},\quad \text{where }\,\mathcal Q = 	
		\mu I + \begin{bmatrix}
			O&A^\top \\ -A&O
		\end{bmatrix}.
	\end{equation}
	Note that $\mathcal Q$ is strongly monotone and $\sigma(\mathcal Q) = \mu \pm i\sqrt{\sigma(A A^\top)}$. The optimality condition is $
	\mathcal Q z^* = 0$, and the unique saddle point is $z^*=(x^*,y^*)=(0,0)$.  Analogously to the unconstrained case, the explicit Euler scheme 
	\begin{equation}\label{eq:gd}
		\frac{z_{k+1}-z_{k}}{\alpha} = -\mathcal Qz_{k}
	\end{equation}
	of the saddle-point dynamics $z' = -\mathcal Qz$ results in a large spectral radius $\mathcal O(1-\mu^2/\nm{A}^2)$ with the step size $\alpha = \mathcal O(\mu^2/\nm{A}^2)$. However, by the theoretical lower complexity bound \cite{azizian_accelerating_2020,ibrahim_linear_2020,zhang_lower_2022}, this is suboptimal.

	Following  \cite[Section 2]{luo_differential_2022}, we aim to find some lifting  system $\bm \xi' = \mathcal G\bm \xi$ with $\bm \xi = (z,u)$, which realizes the improved spectral radius $\mathcal O(1-\mu/\nm{A})$ with larger step size $\alpha = \mathcal O(\mu/\nm{A})$. Motivated by the NAG model \cref{eq:Gnag}, a natural one shall be
	\[
	\mathcal G_{\rm nag} = 
	\begin{bmatrix}
		-I&I\\
		I-\mathcal Q/\mu& -I
	\end{bmatrix}.
	\]
	Expanding the system $\bm \xi' = \mathcal G_{\rm nag}\bm \xi$ with $z=(x,y)$ and $u = (v,w)$ gives
	\begin{equation}\label{eq:nag-sys}
		\begin{pmatrix}
			x'\\y'
		\end{pmatrix}
		= \begin{pmatrix}
			v-x\\w-y
		\end{pmatrix},\quad 
		\begin{pmatrix}
			v'\\ w'
		\end{pmatrix}
		=\begin{pmatrix}
			x-v\\y-w
		\end{pmatrix}
		-\frac{1}{\mu }\begin{pmatrix}
			\mu x+A^\top y\\
			\mu y-Ax
		\end{pmatrix}.
	\end{equation}
	This is very close to the saddle-point dynamics in \cite[Eq.(3.2) in Chapter 3]{mccreesh_accelerated_2019}, which has one more time rescaling factor. Asymptotic stability of \cref{eq:nag-sys} has been proved by \cite[Theorem 4.1.1]{mccreesh_accelerated_2019}, under the assumption that {\it the coupling matrix $A$ is not dominated}. Otherwise, we will see {\it instability}. In other words, the {\it naive} model \cref{eq:nag-sys} is unstable; see \cref{fig:spec-phase}.
	\begin{prop}\label{prop:nag}
		If $\nm{A}>2\mu$, then there must be some $\lambda\in \sigma(\mathcal G_{\rm nag})$ with $\Re(\lambda)>0$. 
	\end{prop}
	\begin{proof}
		Let $\lambda\in\sigma( \mathcal G_{\rm nag})$, then a primal calculation gives 
			$(\lambda+1)^2 = \pm i\sqrt{\delta}/\mu$, where $\delta\in\sigma(A A^\top )$. Since $\lambda_{\max}(A A^\top  )=\nm{A}^2>4\mu^2$, there exists at least one $\delta^*\in\sigma(AA^\top)$ satisfying $\delta^*>4\mu^2$. This results in a trouble eigenvalue $\lambda^*\in\sigma(\mathcal G_{\rm nag})$ such that
			\[
			\lambda^*  = -1+\sqrt{\frac{\sqrt{\delta^*}}{2\mu}}\pm i\sqrt{\frac{\sqrt{\delta^*}}{2\mu}},\quad \text{with }\,\Re(\lambda^*) = -1+\sqrt{\frac{\sqrt{\delta^*}}{2\mu}}>0.
			\]
			This finishes the proof.
		\end{proof}
		\subsection{NAG with velocity correction}
\begin{figure}[H]
			\centering 
			\includegraphics[width=13cm]{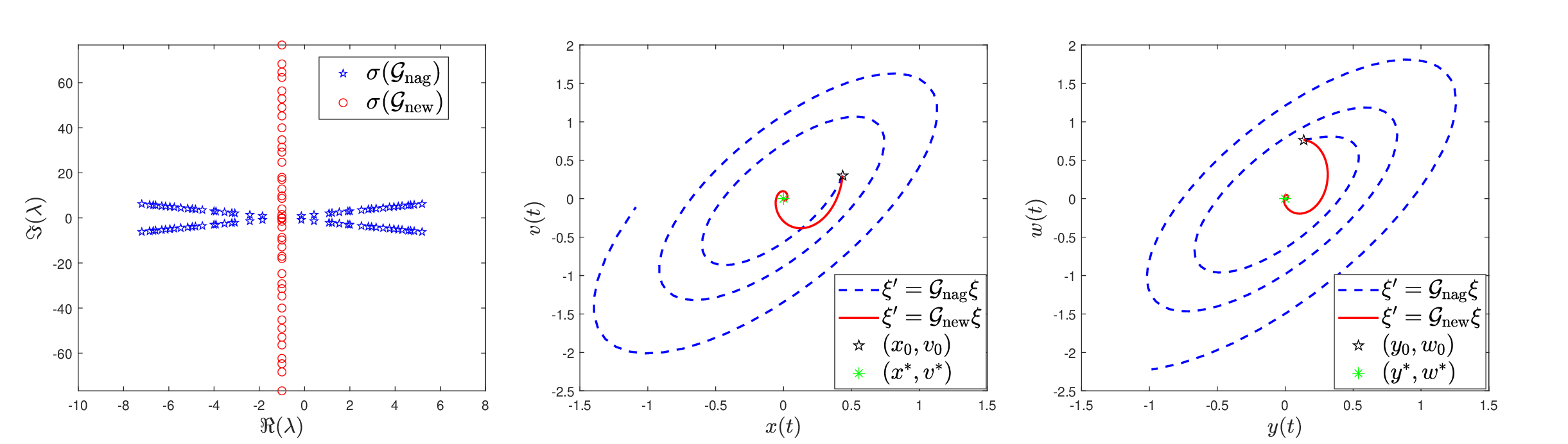}
			\caption{ Numerical illustrations of the spectrum  and the trajectories of \cref{eq:nag-sys,eq:APDG-spec}. }
			\label{fig:spec-phase}
		\end{figure}

		Motivated by the {\it velocity correction} of the continuous models in \cite{Luo2026b,luo_primal-dual_2022}, we modify \cref{eq:nag-sys} as follows
		\begin{equation}\label{eq:APDG-spec}
			\begin{pmatrix}
				x'\\y'
			\end{pmatrix}'
			= \begin{pmatrix}
				v-x\\w-y
			\end{pmatrix},\quad 
			\begin{pmatrix}
				v'\\ w'
			\end{pmatrix}
			=\begin{pmatrix}
				x-v\\y-w
			\end{pmatrix}
			-\frac{1}{\mu }\begin{pmatrix}
				\mu x+A^\top (y+y')\\
				\mu y-A(x+x')
			\end{pmatrix},
		\end{equation}
		which is related to $\bm \xi' = \mathcal G_{\rm new}\bm \xi$ with
		\begin{equation}\label{eq:Gnew}
			\mathcal G_{\rm new}= 
			\begin{bmatrix}
				-I&I\\
				O& -\mathcal Q/\mu
			\end{bmatrix}.
		\end{equation}
		Similarly with \cref{prop:nag}, any $\lambda\in\sigma(\mathcal G_{\rm new})$ satisfies the quadratic equation 
		\begin{equation}\label{eq:eig-Gnew}
			(\lambda+1)(\lambda+1\pm i\sqrt{\delta}/\mu) = 0,
		\end{equation}
		with some $\delta\in\sigma(A A^\top  )$.
		Thus the modified one \cref{eq:APDG-spec} is stable since  $ \sigma(\mathcal G_{\rm new})$ is located at the straight line $-1 + i\R$; see \cref{fig:spec-phase}. 
		
		Write $\mathcal G_{\rm new} = -I +\mathcal L + \mathcal U$, where $\mathcal L$ and $\mathcal U$ are respectively the strictly lower triangular and upper triangular part of $\mathcal G_{\rm new}$. Let $\alpha>0$ be the step size, we now consider the Gauss--Seidel splitting
		\begin{equation}\label{eq:gs}
			\frac{\bm \xi_{k+1}-\bm \xi_k}{\alpha} = (-I+\mathcal L)\bm \xi_{k+1}+\mathcal U\bm \xi_k,
		\end{equation}
		where the lower triangular part adopts the implicit discretization.
		A direct computation gives 
		\[
		\bm \xi_{k+1} = \widehat{\mathcal G}\bm \xi_k,\quad \widehat{\mathcal G} = \left[(1+\alpha)I-\alpha \mathcal L\right]^{-1}\left(I+\alpha \mathcal U\right).
		\]   
		\begin{thm}
			\label{thm:spec-gs}
			If $0<\alpha\leq 2\mu/\nm{A}$, then $\rho(\widehat{\mathcal G}) = \frac{1}{1+\alpha}$. In particular, taking $\alpha = \mu/\nm{A}$ gives
			\[
			\rho(\widehat{\mathcal G}) = \frac{1}{1+\mu/\nm{A}}.
			\]
		\end{thm}
		\begin{proof}
			See \cref{app:spec-gs}.
		\end{proof}
		
		Recall that for the unconstrained case \cref{eq:min-quad}, we see acceleration from $  \mathcal O(1-1/\kappa(Q))$ to $\mathcal O(1-1/\sqrt{\kappa(Q)})$. 
		For the quadratic game \cref{eq:quad-game}, it has $\mu\leq \Re(\lambda)$ and $\snm{\lambda}\leq \sqrt{\mu^2+\nm{A}^2}$ for all $\lambda\in\sigma(\mathcal Q)$. From \cref{thm:spec-gs}, we see a larger step size $\alpha=\mathcal O(\mu/\nm{A})$ leads to the improved spectral radius $\mathcal O(1-\mu/\nm{A})$, which is optimal \cite{azizian_accelerating_2020,ibrahim_linear_2020} and brings acceleration.

\section{Continuous Model and Exponential Decay}
\label{sec:APDG}
Based on the spectral analysis of the quadratic case, we propose our continuous time models for the saddle-point problem \cref{eq:minmax-L}. We consider firstly the Euclidean case (standard convexity) and then non-Euclidean case (relative convexity). The main ingredients are dynamically changing parameters and Lyapunov analysis. 
\subsection{Accelerated primal-dual gradient flow}
Let us firstly focus on the special (Euclidean) case that the norms on $\mathcal X$ and $\mathcal Y$ are induced respectively by its inner products. In other words, we have $\nm{\cdot} = \sqrt{\inner{\cdot,\cdot}}$. By the template \cref{eq:APDG-spec}, we propose an {\it accelerated primal-dual gradient (APDG) flow}:
\begin{subnumcases}{\label{eq:APDG-sys}}
	{}	x' = v-x,
	\label{eq:APDG-sys-x}\\
	{}			y' = w-y,
	\label{eq:APDG-sys-y}\\	
	{}	\gamma v'  \in \sgf (x-v)-\partial_x \mathcal L(x,w),
	\label{eq:APDG-sys-v}\\
	{}	\beta w'  \in\sgg(y-w)+\partial_y \mathcal L(v,y).
	\label{eq:APDG-sys-w}
\end{subnumcases}
Here and in the sequel, the primal parameter $\gamma$ and the dual parameter $\beta$ satisfy
\begin{equation}\label{eq:gama-beta}
	\gamma' = \sgf-\gamma,\quad 		\beta' =\sgg-\beta,
\end{equation} 
with arbitrary positive initial conditions $\gamma(0)=\gamma_0>0$ and $\beta(0)=\beta_0>0$. 
It is not hard to find the exact solutions of \cref{eq:gama-beta}:
\begin{equation}\label{eq:gama-beta-sol}
	\gamma(t) = \gamma_0e^{-t}+\sgf(1-e^{-t}),\quad \beta(t) = \beta_0e^{-t}+\sgg(1-e^{-t}).
\end{equation}
Eliminating the inertial variables $v$ and $w$ from \cref{eq:APDG-sys} leads to a novel second-order ODE model
\begin{subnumcases}{\label{eq:APDG}}
	{}	\gamma x''  +(\sgf+\gamma)x'+\partial_x \mathcal L(x,y+y') = 0,
	\label{eq:APDG-x}\\
	{}	\beta y''  +(\sgg+\beta)y'-\partial_y \mathcal L(x+x',y) = 0.
	\label{eq:APDG-y}
\end{subnumcases}

Under the smooth setting, our APDG flow \cref{eq:APDG-sys} admits a unique global $C^1$-smooth solution $Z = (x,v,y,w):[0,\infty)\to \mathcal Z$. Under the standard convexity (cf.\cref{eq:mu-ineq}), we are able to establish the exponential decay property of the Lyapunov function along the continuous trajectory.
\begin{thm}\label{thm:conv-APDG}
	Suppose the norms on $\mathcal X$ and $\mathcal Y$ are induced respectively by its inner products.
	Assume $\dom f = \mathcal X,\,\dom g = \mathcal Y$ and $f/g$ are $C^1_{L_f}/C^1_{L_g}$-smooth and  $\mu_f/\mu_g$-convex . There is 	a unique global $C^1$-smooth solution $ Z :[0,\infty)\to \mathcal Z$ to the APDG flow \cref{eq:APDG-sys}.
	Moreover, we have
	\begin{equation}\label{eq:conv-APDG}
		\frac{\dd }{\dd t}\mathcal E(\Theta,Z;\widehat{\bm z})\leq  -\mathcal E(\Theta,Z;\widehat{\bm z}) ,
	\end{equation}
	where  $\mathcal E(\Theta,Z;\widehat{\bm z})$ is defined by \cref{eq:Et} with any $\widehat{\bm z}=(\widehat{x},\widehat{y})\in \mathcal X\times \mathcal Y$. This leads to the exponential decay
	\[
	\mathcal E(\Theta(t),Z(t);\widehat{\bm z})\leq  	e^{-t}\mathcal E(\Theta(0),Z(0);\widehat{\bm z}),\quad \forall\,0<t<\infty.
	\]
\end{thm}
\begin{proof} 
	Under the smooth assumption, the differential inclusion  \cref{eq:APDG-sys} becomes
	\begin{equation}\label{eq:APDG-sys-smooth}
		Z'(t) = F(t,Z(t)),\quad\text{with}\quad F(t,Z) = 
		\begin{pmatrix}
			v-x\\
			\frac{1}{	\gamma(t)}\left( \sgf (x-v)-\nabla f(x)-A^* u\right)\\
			w-y\\
			\frac{	1}{\beta(t)}\left(\sgg(y-w)+Av-\nabla g(y)\right)
		\end{pmatrix}.
	\end{equation}
	According to \cref{eq:gama-beta-sol} and the Lipschitz assumptions on $\nabla f$ and $\nabla g$, we claim that $F(t,\cdot)$ is Lipschitz continuous over $\mathcal Z$ for all $0<t<\infty$. In particular, the Lipschitz constant $L(t)$ is bounded by
	\[
	L(t) \leq C\left(1+\frac{\mu_f+L_f+\opnm{A}}{\gamma(t)}+\frac{\mu_g+L_g+\opnm{A}}{\beta(t)}\right)\quad\forall\, t>0,
	\]
	with some universal constant $C>0$. Hence, $L(t)$ is locally integrable over $(0,\infty)$, and by the Picard--Lindel\"{o}f theorem (cf.\cite[Theorem 2.5 and Corollary 2.6]{Teschl2012}), there exists a unique global $C^1$-smooth solution to \cref{eq:APDG-sys}.  
	
	It remains to establish \cref{eq:conv-APDG}.
	By \cref{eq:APDG-sys,eq:gama-beta} and the assumption $\nm{\cdot} = \sqrt{\inner{\cdot,\cdot}}$, a direct computation gives
	\begin{equation}\label{eq:dtE}
		\begin{aligned}
			\frac{\dd }{\dd t}\mathcal E(\Theta,Z;\widehat{\bm z})=	\dual{\nabla_\Theta\mathcal E, \Theta'}+\dual{\nabla_Z\mathcal E, Z'}		:={}&\mathbb A_1+\mathbb A_2+\mathbb A_3+\mathbb A_4,
		\end{aligned}
	\end{equation}
	where 
	\[
	\left\{
	\begin{aligned}
		\mathbb A_1={}&	 \frac{\sgf-\gamma}{2}\nm{v-\widehat{x}}^2
		+\frac{\sgg-\beta}{2}\nm{w-\widehat{y}}^2,\\
		\mathbb A_2={}&   \sgf\inner{x-v,v-\widehat{x}}+ \sgg\inner{y-w,w-\widehat{y}},\\
		\mathbb A_3={}&	\dual{\nabla_y\mathcal L(v ,y),w-\widehat{y} }	-\dual{\nabla_y\mathcal L(\widehat{x},y),w-y},\\
		\mathbb A_4={}&	\dual{ \nabla_x\mathcal L(x,\widehat{y}),v-x}-\dual{\nabla_x\mathcal L(x,w),v-\widehat{x}}.
	\end{aligned}
	\right.
	\]
	Applying the three-term identity \cref{eq:3-id-R} with the natural identity mapping $R=I:\mathcal X\to \mathcal X^*$ gives that
	\[
	\begin{aligned}
		\mathbb A_2 = {}&\frac{\sgf}{2}\left(\nm{x-\widehat{x}}^2-\nm{v-\widehat{x}}^2-\nm{x-v}^2\right)
		+\frac{\sgg}{2}\left(\nm{y-\widehat{y}}^2-\nm{w-\widehat{y}}^2-\nm{y-w}^2\right).
	\end{aligned}
	\]
	For $\mathbb A_3$, by \cref{eq:subg} and the  strong convexity of $\widehat{g}(\cdot;\widehat{\bm z})$, there follows 
	\[
	\begin{aligned}
		\mathbb A_3={}&\dual{A (v-\widehat{x}),w-\widehat{y}}+\dual{\nabla_y\mathcal L(\widehat{x} ,y),w-\widehat{y} }	-\dual{\nabla_y\mathcal L(\widehat{x},y),w-y}\\ 
		={}&\dual{A (v-\widehat{x}),w-\widehat{y}}+\dual{ \nabla_y\mathcal L(\widehat{x},y),y-\widehat{y}}\\		
		={}&\dual{A (v-\widehat{x}),w-\widehat{y}}-\dual{ \nabla\widehat{g}(y),y-\widehat{y}}\quad\left(\text{by \cref{eq:subg}}\right)\\				
		\leq  {}& \dual{A (v-\widehat{x}),w-\widehat{y}}+
		g(\widehat{y})-\widehat{g}(y;\widehat{\bm z})-\frac{\sgg}{2}
		\nm{y-\widehat{y}}^2, 
	\end{aligned}
	\]
	and  similarly
	\[
	\begin{aligned}
		\mathbb A_4 		
		\leq  {}& -\dual{A^* (w-\widehat{y}),v-\widehat{x}}+
		f(\widehat{x})-\widehat{f}(x;\widehat{\bm z})-\frac{\sgf}{2}
		\nm{x-\widehat{x}}^2.
	\end{aligned}
	\]
	This together with the fact \cref{eq:gap-obj} gives
	\begin{equation}\label{eq:i3-i4}
		\begin{aligned}
			\mathbb A_3+		\mathbb A_4\leq{}&
			g(\widehat{y})-\widehat{g}(y;\widehat{\bm z})-\frac{\sgg}{2}
			\nm{y-\widehat{y}}^2+ f(\widehat{x})-\widehat{f}(x;\widehat{\bm z})-\frac{\sgf}{2}
			\nm{x-\widehat{x}}^2\\
			=&\mathcal L(\widehat{x}, y)-\mathcal L(x, \widehat{y})
			-\frac{\sgg}{2}
			\nm{y-\widehat{y}}^2-\frac{\sgf}{2}
			\nm{x-\widehat{x}}^2.
		\end{aligned}
	\end{equation}
	Therefore, combining this with $\mathbb A_1$ and $\mathbb A_2$ leads to
	\[
	\begin{aligned}
		\frac{\dd }{\dd t}\mathcal E(\Theta,Z;\widehat{\bm z})\leq  &-\mathcal E(\Theta,Z;\widehat{\bm z})
		-		\frac{\sgg}{2}\nm{y-w}^2
		-		\frac{\sgf}{2}\nm{x-v}^2.
	\end{aligned}
	\]
	Using the substitutions $x' = v-x$ and $y' = w-y$, we obtain \cref{eq:conv-APDG} and complete the proof.
\end{proof}
\subsection{Accelerated Bregman primal-dual gradient flow}
We then move to the general non-Euclidean case that the norms on $\mathcal X$ and $\mathcal Y$ are not necessarily induced by its inner products. In this case, we consider the Bregman distances and the relative convexity w.r.t. the prox-functions $\phi:\dom\phi\to\R$ and $\psi:\dom\psi\to\R$. More precisely, motivated by the accelerated Bregman primal-dual flow \cite[Eq.(24)]{luo_universal_2024} for solving linearly constrained optimization, we propose the following {\it accelerated Bregman primal-dual gradient (ABPDG) flow}:
\begin{subnumcases}{\label{eq:ABPDG-sys}}
	{}	x' = v-x,
	\label{eq:ABPDG-sys-x}\\
	{}			y' = w -y,
	\label{eq:ABPDG-sys-y}\\	
	{}	\gamma \frac{\dd }{\dd t}\nabla \phi (v) \in\sgf(\nabla \phi(x)-\nabla \phi (v))-\partial_x \mathcal L(x,w),
	\label{eq:ABPDG-sys-v}\\
	{}	\beta \frac{\dd }{\dd t}\nabla \psi(w)  \in \sgg(\nabla \psi(y)-\nabla \psi(w) )+\partial_y\mathcal L(v,y),
	\label{eq:ABPDG-sys-w}
\end{subnumcases}
where $\gamma$ and $\beta$ are still governed by \cref{eq:gama-beta}. Analogously to \cref{eq:APDG}, we can drop $v$ and $w$ to obtain
\[
\begin{aligned}
	{}	\gamma \frac{\dd }{\dd t}\nabla \phi (x+x') \in\sgf(\nabla \phi(x)-\nabla \phi (x+x'))-\partial_x \mathcal L(x,y+y'),\\
	{}	\beta \frac{\dd }{\dd t}\nabla \psi(y+y')  \in \sgg(\nabla \psi(y)-\nabla \psi(y+y') )+\partial_y\mathcal L(x+x',y).
\end{aligned}
\]
\begin{rem}
	If the norms on $\mathcal X$ and $\mathcal Y$ are induced respectively by its inner products, then  the ABPDG flow  \cref{eq:ABPDG-sys}  amounts to the APDG flow  \cref{eq:APDG-sys} by letting $\phi(x)=1/2\nm{x}^2$ and $\psi(y)=1/2\nm{y}^2$.
\end{rem}
\begin{rem}
	Let us discuss the relation and difference among our models and existing works. We mainly focus on primal-dual dynamical systems for the saddle-point problem \cref{eq:minmax-L} and the equality constrained problem \cref{eq:min-eq-1b} with $Ax=b$, which is a special case of \cref{eq:minmax-L}  with $g(y) = \dual{b,y}$.
	\begin{itemize}
		\item Linearly constrained problem \cref{eq:min-eq-1b}: Looking at the equations for the primal variable and the dual variable (Lagrange multiplier), most works adopt either second-order + second-order \cite{he_convergence_rate_2021,zeng_distributed_2018,Zhao2023} or first-order + first-order structures \cite{chen_transformed_2023,luo_primal-dual_2022}. Besides, in \cite{he_second_2022,Luo2026b,Luo2025,luo_universal_2024}, the hybrid  second-order + first-order combination has been considered. When $g(y) = \dual{b,y}$, both of our APDG flow \cref{eq:APDG} and ABPDG flow \cref{eq:ABPDG-sys} provide new second-order + second-order models for the equality constrained problem \cref{eq:min-eq-1b}. Especially, the latter is equipped with general Bregman distances respectively for the primal variable and the Lagrange multiplier, which can be adapted to proper metrics and brings further preconditioning effect. 
		\item Saddle-point problem \cref{eq:minmax-L}: We are aware of some recent works  \cite{Chen2023,He2024,Li2026}. Compared with these, ours combine Bregman distances with subtlety time rescaling parameters for primal and dual variables. This paves the way for designing new (nonlinear) primal-dual splitting methods for both convex and (partially) strongly convex problems in a unified way.  
	\end{itemize}
\end{rem}
\begin{thm}\label{thm:conv-ABPDG}
	Assume $\dom f = \mathcal X,\,\dom g = \mathcal Y$ and $f/g$ are $C^1$-smooth and relatively $\mu_f/\mu_g$-convex with respect to $\phi/\psi$. Let $Z = (x,v,y,w):[0,\infty)\to\mathcal Z$ be a global $C^1$-smooth solution to the ABPDG flow \cref{eq:ABPDG-sys}, then for the Lyapunov function \cref{eq:Ft}, we have
	\begin{equation}\label{eq:conv-ABPDG}
		\frac{\dd }{\dd t}\mathcal E_D (\Theta,Z;\widehat{\bm z})\leq  -\mathcal E_D (\Theta,Z;\widehat{\bm z})
		-\sgf D_\phi(v,x)
		-	\sgg D_\psi(w,y),
	\end{equation}
	which implies 
	\begin{equation}\label{eq:exp-EDt}
		\mathcal E_D (\Theta(t),Z(t);\widehat{\bm z}) 
		\leq e^{-t}	\mathcal E_D (\Theta(0),Z(0);\widehat{\bm z}),\quad\forall\,0<t<\infty,
	\end{equation}
	for any $\widehat{\bm z}=(\widehat{x},\widehat{y})\in \mathcal X\times \mathcal Y$. 
\end{thm}
\begin{proof}
	Observing \cref{def:bregman}, a direct calculation gives 
	\[
	\frac{\dd}{\dd t} D_{\phi}(\widehat{x},v)  
	=\dual{\frac{\dd}{\dd t}\nabla\phi(v),v-\widehat{x}},\quad 
	\frac{\dd}{\dd t} D_\psi(\widehat{y},w)  
	= \dual{\frac{\dd}{\dd t}\nabla\psi(w),w-\widehat{y}}.
	\]
	Therefore, similarly with \cref{eq:dtE}, we take the time derivative on \cref{eq:ABPDG-sys} and use the right hand side of \cref{eq:ABPDG-sys} to obtain the expansion
	\[
	\begin{aligned}
		\frac{\dd }{\dd t}\mathcal E_D (\Theta,Z;\widehat{\bm z}) ={}&\mathbb A_1+\mathbb A_2+\mathbb A_3+\mathbb A_4,
	\end{aligned}
	\]
	where 
	\[
	\left\{
	\begin{aligned}
		\mathbb A_1={}&	 (\sgf-\gamma)D_\phi(\widehat{x},v)+(\sgg-\beta)D_\psi(\widehat{y},w),\\
		\mathbb A_2={}&   \sgf\dual{\nabla\phi(x)-\nabla\phi(v),v-\widehat{x}}+ \sgg\dual{\nabla\psi(y)-\nabla\psi(w),w-\widehat{y}},\\
		\mathbb A_3={}&	\dual{\nabla_y\mathcal L(v ,y),w-\widehat{y} }	-\dual{\nabla_y\mathcal L(\widehat{x},y),w-y},\\
		\mathbb A_4={}&	\dual{ \nabla_x\mathcal L(x,\widehat{y}),v-x}-\dual{\nabla_x\mathcal L(x,w),v-\widehat{x}}.
	\end{aligned}
	\right.
	\]
	Following the proof of \cref{eq:conv-APDG} and using the relatively strong convexity  (cf.\cref{def:rela-strng-conv}), the estimate of $\mathbb A_3$ and $\mathbb A_4$ becomes (cf.\cref{eq:i3-i4})
	\[
	\begin{aligned}
		\mathbb A_3+		\mathbb A_4\leq{}&\mathcal L(\widehat{x}, y)-\mathcal L(x, \widehat{y})
		-\sgg
		D_\psi(\widehat{y},y)-\sgf
		D_\phi(\widehat{x},x).
	\end{aligned}
	\]
	For the second term $\mathbb A_2$, applying \cref{lem:3-id} yields that 
	\[
	\begin{aligned}
		\mathbb A_2 = {}& \sgf \left(D_\phi(\widehat{x}, x)-D_\phi(\widehat{x}, v)-D_\phi(v, x) \right)+ \sgg \left(D_\psi(\widehat{y}, y)-D_\psi(\widehat{y},w)-D_\psi(w, y) \right).
	\end{aligned}
	\] 
	Consequently, collecting all the terms implies \cref{eq:conv-ABPDG} and finishes the proof.
\end{proof}
\begin{rem}\label{lem:partial-gap}
	Since \cref{eq:exp-EDt} holds for any $\widehat{\bm z}=(\widehat{x},\widehat{y})\in \mathcal X\times \mathcal Y$, a further argument leads to the exponential convergence of the Lagrangian gap $\mathcal L(x(t), \widehat{y})-\mathcal L(\widehat{x}, y(t))$ and the partial duality gap $	\sup_{\widehat{\bm z}\in  \widehat{\mathcal X}\times \widehat{\mathcal Y} } \mathcal L(x(t), \widehat{y})-\mathcal L(\widehat{x}, y(t))$ over any compact subset $ \widehat{\mathcal X}\times \widehat{\mathcal Y}\subset \mathcal X \times  \mathcal Y $. For simplicity, we omit the details.
\end{rem} 

\section{Accelerated Bregman Primal-Dual Proximal Splitting}
\label{sec:ABPD-PS}
In this section, we propose an accelerated Bregman primal-dual proximal splitting for solving the saddle-point problem \cref{eq:minmax-L} under the following setting.
\begin{assum}\label{assum:ABPD-PS-f-g} 
	$f:\mathcal X\to\R\cup\{+\infty\}$ is relatively $\mu_f$-convex w.r.t. $\phi:\dom\phi\to\R$ with  $\sgf\geq0$, and $g:\mathcal Y\to\R\cup\{+\infty\}$ is relatively $\mu_g$-convex w.r.t. $\psi:\dom\psi\to\R$ with  $\sgg\geq0$. 
\end{assum}
\subsection{Discretization}
Let us consider a semi-implicit Euler (SIE) scheme for the ABPDG flow \cref{eq:ABPDG-sys}. Given $Z_k = (x_k,v_k,y_k,w_k)\in\mathcal Z$ and $\Theta_k = (\gamma_k,\beta_k)\in\R_+^2$, find $\alpha_k>0$ and update $Z_{k+1} = (x_{k+1},v_{k+1},y_{k+1},w_{k+1})\in\mathcal Z$ by that 
\begin{subnumcases}{\label{eq:ABPD-PS}} 
	{}	\frac{x_{k+1}-x_k}{\alpha_k}=v_{k+1}-x_{k+1},
	\label{eq:ABPD-PS-x}\\	
	{}	\frac{y_{k+1}-y_k}{\alpha_{k}}= \eta_k(w_{k+1}-y_{k+1}),
	\label{eq:ABPD-PS-y}\\
	{}	\gamma_k\frac{\nabla\phi(v_{k+1})-\nabla\phi(v_k)}{\alpha_k}\in
	\sgf (\nabla\phi(x_{k+1})-\nabla\phi(v_{k+1}))-
	\partial_x \mathcal L(x_{k+1}, 	  w_k),
	\label{eq:ABPD-PS-v}\\
	{}	\beta_k\frac{\nabla\psi(w_{k+1})-\nabla\psi(w_k)}{\alpha_k} \in\sgg ( \nabla\psi(y_{k+1})-\nabla\psi(w_{k+1}))+ \eta_k \partial_y \mathcal L({\bar v_{k+1}},y_{k+1}),
	\label{eq:ABPD-PS-w}
\end{subnumcases} 
where  
\begin{equation}\label{eq:etak-b-vk1}
	\eta_k:=\frac{\alpha_{k+1}}{\alpha_k}(1+\alpha_k),\quad 	\bar v_{k+1} := v_{k+1}+\frac{v_{k+1} -v_k}{\eta_k}.
\end{equation}
To update 
$\Theta_{k+1} = (\gamma_{k+1},\beta_{k+1})\subset\R_+^2$, the continuous parameter equation \cref{eq:gama-beta} is discretized by 
\begin{equation}\label{eq:ABPD-PS-gk-bk}
	\frac{\gamma_{k+1}-\gamma_k}{\alpha_k}={}\sgf -\gamma_{k+1},
	\quad
	\frac{	\beta_{k+1}-	\beta_k}{\alpha_k} ={}\sgg -\beta_{k+1},
\end{equation}
with arbitrary positive initial values $\gamma_0,\,\beta_0>0$. 
Later in \cref{algo:ABPD-PS}, we reformulate the SIE scheme \cref{eq:ABPD-PS} as a primal-dual proximal splitting method. 

\subsection{Lyapunov analysis}
\label{sec:ABPD-PS-lyap}
In this part, we aim to establish the Lyapunov analysis of the SIE scheme \cref{eq:ABPD-PS} by using the discrete Lyapunov function $\mathcal H_k(\widehat{\bm z})$ defined by \cref{eq:Hk}. Motivated by the decomposition given in \cref{lem:decomp-diff-Hk}, in the next, we shall focus on $\mathbb I_{21}+\mathbb I_{22}$ and  $\mathbb I_{31}+\mathbb I_{32}$ and prove the corresponding upper bounds separately in \cref{lem:ABPD-PS-est-i2,lem:ABPD-PS-est-i3}.
\begin{lem}
	\label{lem:ABPD-PS-ak}
	Let $\{\gamma_k\}$ and $\{\beta_k\}$ be defined by \cref{eq:ABPD-PS-gk-bk} with any $\alpha_k>0$, then
	\begin{subequations}	\label{eq:gk-bk-lb-ub}
		\begin{align}
			\label{eq:gk-lb-ub}
			{}&	\max\{\gamma_0\theta_k,\,	\min\{\sgf ,\,\gamma_0\}\}\leq \gamma_k\leq 	\max\{\sgf ,\,\gamma_0\},\\
			\label{eq:bk-lb-ub}
			{}&	\max\{\beta_0\theta_k,\,	\min\{\sgg ,\,\beta_0\}\}\leq\beta_k\leq \max\{\sgg ,\,\beta_0\},
		\end{align}
	\end{subequations}
	where $\{\theta_k\}$ is defined by \cref{eq:thetak}. If
	$\opnm{A}^2\alpha_k^2= \gamma_k\beta_k$, then $\eta_k = \alpha_{k+1}(1+\alpha_k)/\alpha_k\geq 1$. 
\end{lem}
\begin{proof}
	By \cref{eq:ABPD-PS-gk-bk}, we have 
	\begin{equation}\label{eq:gk1-bk1}
		\gamma_{k+1} = \frac{\sgf \alpha_k+\gamma_k}{1+\alpha_k},\quad
		\beta_{k+1}=\frac{\sgg \alpha_k+\beta_k}{1+\alpha_k}.
	\end{equation}
	In view of this and \cref{eq:diff-tk}, a standard induction argument verifies \cref{eq:gk-bk-lb-ub}. Furthermore, if $\opnm{A}^2\alpha_k^2= \gamma_k\beta_k$, then it follows that
	\begin{equation}\label{eq:etak}
		\begin{aligned}
			\eta_k ={} \frac{1+\alpha_k}{\alpha_k}\alpha_{k+1}=\frac{1+\alpha_k}{\alpha_k\opnm{A}} \sqrt{\gamma_{k+1}\beta_{k+1}} 
			={}&  \frac{1}{\alpha_k\opnm{A}}\sqrt{(\sgf \alpha_k+\gamma_k)(\sgg \alpha_k+\beta_k)}\\
			={}&\frac{1}{\sqrt{\gamma_k\beta_k}}\sqrt{(\sgf \alpha_k+\gamma_k)(\sgg \alpha_k+\beta_k)}
			\geq 1.
		\end{aligned}
	\end{equation}
	This completes the proof of this lemma.
\end{proof}

\begin{lem}\label{lem:ABPD-PS-est-i2}
	Given $\Theta_k=(\gamma_k,\beta_k)\in\R_+^2$ and $Z_k=(x_k,v_k,y_k,w_k)\in\mathcal Z$, let $\Theta_{k+1}=(\gamma_{k+1},\beta_{k+1})\in\R_+^2$ and $Z_{k+1}=(x_{k+1},v_{k+1},y_{k+1},w_{k+1})\in\mathcal Z$ be the outputs of the SIE scheme \cref{eq:ABPD-PS} and the parameter equation \cref{eq:ABPD-PS-gk-bk} at the $k$-th step with any step size $\alpha_k>0$. Then under  \cref{assum:ABPD-PS-f-g} we have 
	\begin{equation}\label{eq:ABPD-PS-est-i2}
		\begin{aligned}
			%
			\mathbb I_{21}+\mathbb I_{22}\leq & -\alpha_k\dual{A^*(w_k-\widehat{y}),v_{k+1}-\widehat{x}}
			-\gamma_k D_\phi(v_{k+1},v_{k}),
		\end{aligned}
	\end{equation}
	where $\mathbb I_{21}$ and $\mathbb I_{22}$ are defined in \cref{eq:diff-Hk} with any $\widehat{\bm z}=(\widehat{x},\widehat{y})\in\dom \mathcal L$.
\end{lem}
\begin{proof}
	According to \cref{assum:ABPD-PS-f-g}, $f$ is relatively $\mu_f$-convex w.r.t. $\phi$ and so is  $\widehat{f}(\cdot;\widehat{\bm z})$ by \cref{lem:L-gap-obj-res} (i). This implies that
	\begin{align} 
		\mathbb I_{21}=	{}&	\widehat{f}(x_{k+1};\widehat{\bm z})-\widehat{f}(x_k;\widehat{\bm z})+	\alpha_k\left(\widehat{f}(x_{k+1};\widehat{\bm z})-f(\widehat{x})\right)\notag \\
		\overset{\text{by \cref{eq:scv-f}}}{	\leq} &
		\dual{\partial \widehat{f}(x_{k+1};\widehat{\bm z}), x_{k+1} -x_k}+
		\alpha_k	\dual{\partial \widehat{f}(x_{k+1};\widehat{\bm z}), x_{k+1} -\widehat{x}}-\sgf \alpha_k D_\phi(\widehat{x},x_{k+1}) \notag \\
		\overset{\text{by \eqref{eq:ABPD-PS-x}}}{	=} &\alpha_k	\dual{\partial \widehat{f}(x_{k+1};\widehat{\bm z}), v_{k+1} -\widehat{x}}-\sgf\alpha_kD_\phi(\widehat{x},x_{k+1}) .\label{eq:est-1-I1-ABPD-PS}
	\end{align}
	Notice that $\partial \widehat{f}(x_{k+1};\widehat{\bm z}) = \partial_x\mathcal L(x_{k+1},\widehat{y})= \partial_x\mathcal L(x_{k+1},w_k)-A^*(w_k-\widehat{y})$, which together with \eqref{eq:ABPD-PS-v} gives  
	\[
	\begin{aligned}
		{}& \alpha_k	\dual{\partial \widehat{f}(x_{k+1};\widehat{\bm z}), v_{k+1} -\widehat{x}}\\
		=&-\alpha_k\dual{A^*(w_k-\widehat{y}),v_{k+1}-\widehat{x}}+\gamma_k\dual{\nabla\phi(v_{k})-\nabla\phi(v_{k+1}),v_{k+1}-\widehat{x}}\\
		{}&\quad+\sgf\alpha_k\dual{\nabla\phi(x_{k+1})-\nabla\phi(v_{k+1}),v_{k+1}-\widehat{x}}.
	\end{aligned}
	\]
	Plugging this into the previous estimate \cref{eq:est-1-I1-ABPD-PS} and using the three-term identity \cref{eq:3-id}, we obtain
	\[
	\begin{aligned}
		\mathbb I_{21}\leq &-\sgf\alpha_k D_\phi(\widehat{x},x_{k+1}) -\alpha_k\dual{A^*(w_k-\widehat{y}),v_{k+1}-\widehat{x}}\\
		{}&\quad + \gamma_k \left( D_\phi(\widehat{x},v_{k})-D_\phi(\widehat{x},v_{k+1})-D_\phi(v_{k+1},v_{k})\right)\\
		{}&\qquad + \sgf \alpha_k \left( D_\phi(\widehat{x},x_{k+1})-D_\phi(\widehat{x},v_{k+1})-D_\phi(v_{k+1},x_{k+1})\right)\\
		=&	-\sgf \alpha_kD_\phi(v_{k+1},x_{k+1}) -\alpha_k\dual{A^*(w_k-\widehat{y}),v_{k+1}-\widehat{x}}\\
		{}&\quad+ \gamma_k D_\phi(\widehat{x},v_{k}) -(\gamma_{k}+\sgf\alpha_k)D_\phi(\widehat{x},v_{k+1})-\gamma_k D_\phi(v_{k+1},v_{k}).
	\end{aligned}
	\]
	In view of \cref{eq:gk1-bk1}, it follows that  $\gamma_k+\sgf\alpha_k=\gamma_{k+1}(1+\alpha_k)$ and 
	\[
	\begin{aligned}
		\mathbb I_{21}\leq &
		-\gamma_k D_\phi(v_{k+1},v_{k}) -\alpha_k\dual{A^*(w_k-\widehat{y}),v_{k+1}-\widehat{x}}\\
		{}&\quad
		-\sgf \alpha_kD_\phi(v_{k+1},x_{k+1})+ \underbrace{\gamma_k D_\phi(\widehat{x},v_{k}) -\gamma_{k+1}(1+\alpha_k)D_\phi(\widehat{x},v_{k+1})}_{=-\mathbb I_{22}}.
	\end{aligned}
	\]
	This implies \cref{eq:ABPD-PS-est-i2} and completes the proof of this lemma.
\end{proof}
\begin{lem}\label{lem:ABPD-PS-est-i3}
	Given $\Theta_k=(\gamma_k,\beta_k)\in\R_+^2$ and $Z_k=(x_k,v_k,y_k,w_k)\in\mathcal Z$, let $\Theta_{k+1}=(\gamma_{k+1},\beta_{k+1})\in\R_+^2$ and $Z_{k+1}=(x_{k+1},v_{k+1},y_{k+1},w_{k+1})\in\mathcal Z$ be the outputs of the SIE scheme \cref{eq:ABPD-PS} and the parameter equation \cref{eq:ABPD-PS-gk-bk} with  $\opnm{A}^2\alpha_k^2=\gamma_k\beta_k$. Then under  \cref{assum:ABPD-PS-f-g} we have 
	\begin{equation}\label{eq:ABPD-PS-est-i3}
		\begin{aligned}
			%
			\mathbb I_{31}+	\mathbb I_{32}\leq {}&\alpha_k\eta_k\dual{A(\bar v_{k+1}-\widehat{x}),w_{k+1}-\widehat{y}}
			-\beta_k D_\psi(w_{k+1},w_{k}) ,
		\end{aligned}
	\end{equation}  
	where $\mathbb I_{31}$ and $\mathbb I_{32}$ are defined in \cref{eq:diff-Hk} with $\widehat{\bm z}=(\widehat{x},\widehat{y})$ being a saddle point to \cref{eq:minmax-L}.
\end{lem}
\begin{proof}
	By \cref{lem:L-gap-obj-res} (ii) we have $\widehat{g}(y_{k+1};\widehat{\bm z})\geq g(\widehat{y})$. Analogously to \cref{eq:est-1-I1-ABPD-PS}, using \cref{eq:scv-f,lem:L-gap-obj-res} (i) and \eqref{eq:ABPD-PS-y} implies that 
	\[
	\begin{aligned}
		\mathbb I_{31}={}&	\widehat{g}(y_{k+1};\widehat{\bm z})-\widehat{g}(y_k;\widehat{\bm z})+ \alpha_k\left(\widehat{g}(y_{k+1};\widehat{\bm z})-g(\widehat{y})\right)\\
		\leq {}&	\widehat{g}(y_{k+1};\widehat{\bm z})-\widehat{g}(y_k;\widehat{\bm z})+ \alpha_k\eta_k\left(\widehat{g}(y_{k+1};\widehat{\bm z})-g(\widehat{y})\right)\\	
		\leq {}& 	\dual{\partial \widehat{g}(y_{k+1};\widehat{\bm z}), y_{k+1} -y_k}+	\alpha_k\eta_k	\dual{\partial \widehat{g}(y_{k+1};\widehat{\bm z}), y_{k+1} -\widehat{y}}  -\sgg\alpha_k\eta_kD_\psi(\widehat{y},y_{k+1}) \\
		\leq{}&\alpha_k\eta_k	\dual{\partial \widehat{g}(y_{k+1};\widehat{\bm z}), w_{k+1} -\widehat{y}}-\sgg\alpha_kD_\psi(\widehat{y},y_{k+1}),
	\end{aligned}
	\]
	where in the second and the last lines, we used the fact $\eta_k\geq 1$, which holds true indeed by  \cref{lem:ABPD-PS-ak} and the step size condition $\opnm{A}^2\alpha_k^2=\gamma_k\beta_k$. Observing that $\partial \widehat{g}(y_{k+1};\widehat{\bm z}) = \partial_y\mathcal L(\bar v_{k+1},y_{k+1})+A(\bar v_{k+1}-\widehat{x})$, we get 	 
	\[
	\begin{aligned}
		\mathbb I_{31}
		\overset{\text{by \eqref{eq:ABPD-PS-w}}}{\leq } &-\sgg\alpha_kD_\psi(\widehat{y},y_{k+1}) +\alpha_k\eta_k\dual{A(\bar v_{k+1}-\widehat{x}),w_{k+1}-\widehat{y}}\\
		{}&\quad + \sgg\alpha_k \left( D_\psi(\widehat{y},y_{k+1})-D_\psi(\widehat{y},w_{k+1})-D_\psi(w_{k+1},y_{k+1})\right)\\
		{}&\qquad + \beta_k \left( D_\psi(\widehat{y},w_{k})-D_\psi(\widehat{y},w_{k+1})-D_\psi(w_{k+1},w_{k})\right)\\
		\overset{\text{by \cref{eq:3-id}}}{=}&-\sgg \alpha_kD_\psi(w_{k+1},y_{k+1}) +\alpha_k\eta_k\dual{A(\bar v_{k+1}-\widehat{x}),w_{k+1}-\widehat{y}}\\
		{}&\quad +\underbrace{\beta_k D_\psi(\widehat{y},w_{k}) -(\beta_k+\sgg \alpha_k )D_\psi(\widehat{y},w_{k+1})}_{=-\mathbb I_{32}}-\beta_k D_\psi(w_{k+1},w_{k}).
	\end{aligned}
	\]
	We mention that  $\beta_k+\sgg\alpha_k=\beta_{k+1}(1+\alpha_k)$ (cf.\cref{eq:gk1-bk1}).  This finishes the proof of \cref{eq:ABPD-PS-est-i3}.
\end{proof} 
\begin{thm}
	\label{thm:conv-ABPD-PS}  
	Given $\Theta_k=(\gamma_k,\beta_k)\in\R_+^2$ and $Z_k=(x_k,v_k,y_k,w_k)\in\mathcal Z$, let $\Theta_{k+1}=(\gamma_{k+1},\beta_{k+1})\in\R_+^2$ and $Z_{k+1}=(x_{k+1},v_{k+1},y_{k+1},w_{k+1})\in\mathcal Z$ be the outputs of the SIE scheme \cref{eq:ABPD-PS} and the parameter equation \cref{eq:ABPD-PS-gk-bk} with  $\opnm{A}^2\alpha_k^2=\gamma_k\beta_k$.  Then under \cref{assum:ABPD-PS-prox-fun,assum:ABPD-PS-f-g} we have 
	\begin{equation}\label{eq:diff-ABPD-PS}
		\begin{aligned}
			{}&\mathcal H_{k+1}(\widehat{\bm z})-\mathcal H_k(\widehat{\bm z})+\alpha_k\mathcal H_{k+1}(\widehat{\bm z})
			\leq{}0,
		\end{aligned}
	\end{equation}
	%
	where $\mathcal H_k(\widehat{\bm z})$ is define by \cref{eq:Hk} with $\widehat{\bm z}=(\widehat{x},\widehat{y})$ being a saddle point to \cref{eq:minmax-L}.
\end{thm}
\begin{proof} 
	Thanks 	to  \cref{lem:ABPD-PS-est-i2,lem:ABPD-PS-est-i3} and the decomposition \cref{eq:diff-Hk},  we obtain 
	\begin{equation}\label{eq:diff-Hk-ABPD-PS}
		\begin{aligned} 
			\mathcal H_{k+1}(\widehat{\bm z})-\mathcal H_k(\widehat{\bm z})+\alpha_k\mathcal H_{k+1}(\widehat{\bm z}) 
			\leq{}&  \Delta_k
			-	\gamma_{k}
			D_\phi(v_k,v_{k+1})- \beta_{k} 
			D_\psi(w_k, w_{k+1}),
		\end{aligned}
	\end{equation}
	where
	\begin{equation}\label{eq:Deltak}
		\begin{aligned}
			\Delta_k={}&\alpha_k\eta_k\dual{A(\bar v_{k+1}-\widehat{x}), w_{k+1}-\widehat{y}}
			-\alpha_k\dual{A(v_{k+1}-\widehat{x}), w_k-\widehat{y}}\\
			{}&\quad+
			\alpha_k\dual{A(v_k-\widehat{x}), w_k-\widehat{y}}-\alpha_{k+1}(1+\alpha_{k})\dual{A(v_{k+1}-\widehat{x}), w_{k+1}-\widehat{y}}.
		\end{aligned}
	\end{equation} 
	Recalling \cref{eq:etak-b-vk1}, we expand the first term as follows
	\[
	\begin{aligned}
		{}&\alpha_k\eta_k\dual{A(\bar v_{k+1}-\widehat{x}), w_{k+1}-\widehat{y}}=\alpha_{k+1}(1+\alpha_{k})\dual{A(\bar v_{k+1}-\widehat{x}), w_{k+1}-\widehat{y}}\\
		={}&\alpha_k\eta_k\dual{A(v_{k+1}-\widehat{x}), w_{k+1}-\widehat{y}}
		+\alpha_k\dual{A(v_{k+1}-v_k), w_{k+1}-\widehat{y}},
	\end{aligned}
	\]
	which yields that $	\Delta_k=\alpha_k\dual{A(v_{k+1}-v_k), w_{k+1}- w_k}$. Hence, using the fact $\opnm{A}^2\alpha_k^2=\gamma_k\beta_k$ and the mean value inequality gives $\Delta_k\leq \gamma_k/2\nm{v_{k+1}-v_k}^2+\beta_k/2\nm{w_{k+1}-w_k}^2$.
	Combining this with \cref{eq:diff-Hk-ABPD-PS} and \cref{assum:ABPD-PS-prox-fun},
	we obtain  \cref{eq:diff-ABPD-PS} and conclude the proof of this theorem.
\end{proof}
\subsection{The overall algorithm and iteration complexity}
\label{sec:ABPD-PS-algo-conv}
The overall SIE scheme \cref{eq:ABPD-PS} with the step size condition $\opnm{A}^2\alpha_k^2= \gamma_k\beta_k$ is listed in \cref{algo:ABPD-PS}, which is called the Accelerated Bregman Primal-Dual Proximal Splitting (ABPD-PS). By \cref{assum:ABPD-PS-f-g}, the objectives of the proximal steps in Lines \ref{line:ABPD-PS-xk1} and \ref{line:ABPD-PS-yk1} of \cref{algo:ABPD-PS}  are (strongly) convex. With the initial setting: $x_0,v_0\in\dom f$ and  $y_0,w_0\in\dom g$, it is clear that  $x_k,v_k\in\dom f$ and  $y_k,w_k\in\dom g$. 
\begin{algorithm}[H]
	\caption{ABPD-PS for solving \cref{eq:minmax-L} under \cref{assum:ABPD-PS-prox-fun,assum:ABPD-PS-f-g} }
	\label{algo:ABPD-PS} 
	\begin{algorithmic}[1] 
		\REQUIRE  Problem parameters: $\sgf  ,\,\sgg  \geq 0$ and $\opnm{A}$.			\\
		\quad ~~~Initial guesses: $x_0,v_0\in\dom f$ and  $y_0, w_0\in\dom g$.\\
		\quad ~~~Initial parameters: $\gamma_0,\,\beta_0>0$.
		\FOR{$k=0,1,\ldots$} 
		\STATE Take $\alpha_k = \sqrt{\gamma_k\beta_k}/\opnm{A}$. \label{line:ABPD-PS-ak}
		\STATE Update $\gamma_{k+1} = (\sgf  \alpha_k+\gamma_k)/(1+\alpha_k)$ and $
		\beta_{k+1}=(\sgg  \alpha_k+\beta_k)/(1+\alpha_k)$.	
		\STATE Compute $\alpha_{k+1} = \sqrt{\gamma_{k+1}\beta_{k+1}}/\opnm{A}$ and $	\eta_k= \alpha_{k+1}(1+\alpha_k)/\alpha_k$. 
		\STATE
		\STATE Define $f_k(x,v):=	f(x)+\gamma_{k+1}D_\phi(v,v_k)+\dual{A^*  w_k,x}$.		
		\STATE Update $x_{k+1}=\mathop{\argmin}\limits_{x\in \mathcal X}\left\{
		f_k(x,v)-\sgf  D_\phi(x,v_k) :\,v =x+ \frac{x-x_k }{\alpha_k}\right\}$. 		\label{line:ABPD-PS-xk1}
		\STATE Update $ v_{k+1}=x_{k+1}+( x_{k+1}- x_{k})/\alpha_k$ and $\bar v_{k+1}  = v_{k+1}+(v_{k+1} -v_k)/\eta_k$. 
		\STATE
		\STATE Define $g_k(y,w):=	g(y)+\frac{1+\alpha_k}{1+\alpha_k\eta_k}\beta_{k+1}D_\psi(w, w_k)-\dual{ A\bar v_{k+1},y}$.
		\STATE Update  $y_{k+1}=\mathop{\argmin}\limits_{y\in \mathcal Y}\left\{
		\eta_k	g_k(y,w)-\sgg  D_\psi(y, w_k) :\,w =y+ \frac{y-y_k}{\alpha_k\eta_k}\right\}$. 	 			 	\label{line:ABPD-PS-yk1}
		\STATE Update $  w_{k+1}=y_{k+1}+( y_{k+1}- y_{k})/(\alpha_k\eta_k)$.
		\ENDFOR
	\end{algorithmic}
\end{algorithm}
\begin{rem}\label{rem:xk1-yk1-eulidean}
	Let $R:\mathcal X\to \mathcal X^*$ and $S:\mathcal Y\to \mathcal Y^*$ be SPD. If we choose quadratic prox-functions
	\begin{equation}\label{eq:phi-R-psi-S}
		\phi(x) = 1/2\nm{x}_R^2,\quad\psi(y) = 1/2\nm{y}_S^2,
	\end{equation}
	such that $\phi$ and $\psi$ are 1-convex, 
	then the proximal steps in Lines \ref{line:ABPD-PS-xk1} and \ref{line:ABPD-PS-yk1} become
	\begin{align}
		x_{k+1}={}&\mathop{\argmin}\limits_{x\in \mathcal X}\!\!\left\{
		f(x)+\dual{A^*  w_k,x}+\frac{\delta_k}{2\alpha_k^2}\nm{x-\widetilde{x}_k}_{R}^2\right\},		\label{eq:ABPD-PS-xk1}\\
		y_{k+1}={}&\mathop{\argmin}\limits_{y\in \mathcal Y}\!\!\left\{
		g(y)-\dual{A\bar v_{k+1},y}+\frac{\tau_k}{2\eta_k\alpha_k^2}\nm{y-\widetilde{y}_k}_{S}^2\right\},
		\label{eq:ABPD-PS-yk1}
	\end{align} 
	where  
	\[
	\begin{aligned}
		\delta_k={}&\sgf \alpha_k+\gamma_k(1+\alpha_k),{}&&
		\widetilde{x}_k=\frac{1}{\delta_k}\left((\sgf \alpha_k+\gamma_k)x_k+\gamma_k\alpha_kv_k\right),\\
		\tau_k = {}&\sgg \alpha_k+\beta_k(1+\eta_k\alpha_k), {}&&\widetilde{y}_k =\frac{1}{\tau_k}((\sgg \alpha_k+\beta_k)y_k+\eta_k\beta_k\alpha_k  w_k).
	\end{aligned}
	\]
\end{rem}
%
\begin{thm}
	\label{thm:rate-ABPD-PS}
	Assume that $\{(x_k,y_k)\}\subset\dom \mathcal L$ is generated by \cref{algo:ABPD-PS}. Then under \cref{assum:ABPD-PS-prox-fun,assum:ABPD-PS-f-g} we have 
	\begin{equation}\label{eq:L-gap-x-y-ABPD-PS}
		\begin{aligned}
			{}&	\mathcal L(x_k,\widehat{y})-	\mathcal L(\widehat{x},y_k)+	\frac{\sgf  }{2}\nm{x_k-\widehat{x}}^2+\frac{\sgg }{2}\nm{y_k-\widehat{y}}^2\leq 2 \theta_k\mathcal H_0(\widehat{\bm z}),
		\end{aligned}
	\end{equation}
	where $\{\theta_k\}$ is defined by \cref{eq:thetak} and  $\widehat{\bm z}=(\widehat{x},\widehat{y})$ a saddle point to \cref{eq:minmax-L}.
	Moreover, if $\gamma_0\geq \sgf ,\,\beta_0\geq \sgg $ and  $\gamma_0\beta_0\leq \opnm{A}^2$, then 
	\begin{equation}\label{eq:L-gap-x-y}
		\begin{aligned} 
			\theta_k
			\leq{}&   \min\left\{\frac{\opnm{A}}{\opnm{A}+\sqrt{\gamma_0\beta_0}k},\,\frac{64\opnm{A}^2}{\sgg \gamma_0k^2},\,\frac{64\opnm{A}^2}{\sgf \beta_0k^2},\,
			\left(1+\frac{\sqrt{\sgf \sgg }}{\opnm{A}}\right)^{-k}
			\right\}.
		\end{aligned}
	\end{equation} 
\end{thm}
\begin{proof} 
	By \cref{assum:ABPD-PS-f-g} and the fact $\opnm{A}^2\alpha_k^2=\gamma_k\beta_k$ (cf. Line \ref{line:ABPD-PS-ak} of \cref{algo:ABPD-PS}), it follows that 
	\begin{equation}\label{eq:cross}
		\alpha_k \snm{\dual{A(v_k-\widehat{x}), w_k-\widehat{y}}}\leq  \gamma_k D_\phi(\widehat{x},v_k)	+ \beta_k  D_\psi(\widehat{y}, w_k).
	\end{equation}
	Applying \cref{lem:decomp-diff-Hk,thm:conv-ABPD-PS} proves \cref{eq:L-gap-x-y-ABPD-PS}. 
	
	Recall the difference equation \cref{eq:diff-tk} of $\{\theta_k\}$. To establish \cref{eq:L-gap-x-y}, the key is to find the upper bound and lower bound of $\alpha_k$.
	Since $\gamma_0\geq \sgf $ and $\beta_0\geq \sgg $, by \cref{eq:gk-bk-lb-ub} we have  $\gamma_k\leq\gamma_0$ and $\beta_k\leq \beta_0$, which implies $\alpha_k=\opnm{A}^{-1}\sqrt{\gamma_k\beta_k}\leq  \sqrt{\gamma_0\beta_0}/\opnm{A}\leq 1$.
	The lower bound of $\alpha_k$ is more subtle. From \cref{eq:gk-bk-lb-ub}, both $\gamma_k$ and $\beta_k$ have two under estimates. Hence, we consider four situations.
	\begin{itemize}
		\item[(i)] First of all, for $\alpha_k\geq \sqrt{\gamma_0\beta_0}\theta_k/\opnm{A}=\alpha_0\theta_{k}$, it follows from \cref{eq:diff-tk} that $	\theta_{k+1}-\theta_k \leq  -\alpha_0\theta_k\theta_{k+1}$. This also implies $1/\theta_{k+1}-1/\theta_k\geq \alpha_0$ and yields that $\theta_k\leq 1/(1+\alpha_0k)\leq \opnm{A}/(\opnm{A}+\sqrt{\gamma_0\beta_0}k)$.
		\item[(ii)] Secondly, for $\alpha_k\geq \sqrt{\sgf \beta_0\theta_k}/\opnm{A}$, we have $\theta_{k+1}-\theta_k \leq  - \sqrt{\sgf \beta_0}\sqrt{\theta_k}\theta_{k+1}/\opnm{A}$. Notice that $\theta_{k+1}/\theta_k = 1/(1+\alpha_k)\geq 1/(1+\alpha_0)\geq 1/2$. Hence, invoking \cref{lem:app-est} implies $\theta_k\leq 64\opnm{A}^2/(\sgf \beta_0k^2)$. 
		\item[(iii)] Similarly, for $\alpha_k\geq \sqrt{\sgg \gamma_0\theta_k}/\opnm{A}$, we get $\theta_k\leq 64\opnm{A}^2/(\sgg \gamma_0k^2)$. 
		\item[(iv)] Finally, for $\alpha_k\geq \sqrt{\sgf \sgg }/\opnm{A}$, we obtain from \cref{eq:thetak} that $\theta_k\leq \left(1+\sqrt{\sgf \sgg }/\opnm{A}\right)^{-k}$.
	\end{itemize}
	Note that these estimates hold simultaneously.  This leads to \cref{eq:L-gap-x-y} and finishes the proof.
\end{proof}

According to convergence rate estimate \cref{eq:L-gap-x-y}, we conclude that \cref{algo:ABPD-PS} admits the optimal iteration complexity \cite{ouyang_lower_2021,zhang_lower_2022} under the {\it proximal algorithm class} \cite[Definition 2.1]{zhang_lower_2022} $(x_k,y_k)\in\mathcal X^k\times \mathcal Y^k$, where
\[
\left\{
\begin{aligned}
	\mathcal X^{k+1}: = {}&{\bf span}\left\{x_i,\,\prox_{\delta_if}(\widetilde{x}_i-\delta_i A^* \widetilde{y}_i):\,\widetilde{x}_i\in\mathcal X^{i},\widetilde{y}_i\in\mathcal Y^i,\,0\leq i\leq k\right\},\\
	\mathcal Y^{k+1}: ={}& {\bf span}\left\{y_i,\,\prox_{\tau_ig}(\widetilde{y}_i+\tau_i A \widetilde{x}_i):\,\widetilde{x}_i\in\mathcal X^{i},\widetilde{y}_i\in\mathcal Y^i,\,  0\leq i\leq k\right\},	
\end{aligned}
\right.
\]
with $\mathcal X^{0}: = {\bf span}\left\{x_0\right\}$ and $\mathcal Y^{0}: = {\bf span}\left\{y_0\right\}$. More precisely, to achieve the accuracy $	\mathcal L(x_k,\widehat{y})-	\mathcal L(\widehat{x},y_k)\leq \epsilon$, the iteration complexity is
\[
\mathcal O	\left(\min\left\{\frac{\opnm{A}}{\epsilon},\,\frac{\opnm{A}}{\sqrt{\sgf \epsilon}},\,\frac{\opnm{A}}{\sqrt{\sgg \epsilon}},\,\frac{\opnm{A}}{\sqrt{\sgf \sgg}}\snm{\ln \epsilon}\right\}\right).
\]
When $\sgf >0$, the iteration complexity for finding $\nm{x_k-\widehat{x}}^2\leq \epsilon$ is 
\[
\mathcal O\left(\min\left\{\frac{\opnm{A}}{\sgf \sqrt{\epsilon}},\,\frac{\opnm{A}}{\sqrt{\sgf \sgg}}\snm{\ln \frac{\sgf }{\epsilon}}\right\}\right).
\]
When $\sgg >0$, the iteration complexity for finding $\nm{y_k-\widehat{y}}^2\leq \epsilon$ is
\[
\mathcal O\left(\min\left\{\frac{\opnm{A}}{\sgg \sqrt{\epsilon}},\,\frac{\opnm{A}}{\sqrt{\sgf \sgg}}\snm{\ln \frac{\sgg }{\epsilon}}\right\}\right).
\] 
\section{Accelerated Bregman Primal-Dual Proximal Gradient Splitting}
\label{sec:ABPD-PGS}
In this section, we mainly focus on the case that both $f$ and $g$ admit the composite structure that satisfies the following setting.  
\begin{assum}\label{assum:sym-ABPD-PGS} 
	$f=f_1+f_2$ where 		$f_2:\mathcal X\to\R\cup\{+\infty\}$ is closed proper convex and	$f_1:\mathcal X\to\R$ is $C_{\Lf }^1$-smooth and relatively $\mu_f$-convex w.r.t. $\phi$ with  $\sgf\geq0$. 	$g=g_1+g_2$ where 		$g_2:\mathcal Y\to\R\cup\{+\infty\}$ is closed proper convex and	$g_1:\mathcal Y\to\R$ is $C_{L_g}^1$-smooth and relatively $\mu_g$-convex w.r.t. $\psi$ with  $\sgg\geq0$. 	
\end{assum} 
\subsection{Algorithm presentation and Lyapunov analysis}
Based on \cref{eq:ABPD-PS}, we propose an variant SIE scheme for the continuous ABPDG flow model \cref{eq:ABPDG-sys} with composite objectives. Given $Z_k = (x_k,v_k,y_k,w_k)\in\mathcal Z$ and $\Theta_k = (\gamma_k,\beta_k)\in\R_+^2$, find $\alpha_k>0$ and update $Z_{k+1} = (x_{k+1},v_{k+1},y_{k+1},w_{k+1})\in\mathcal Z$ as follows 
\begin{subnumcases}{\label{eq:sym-ABPD-PGS}}
	{}	\frac{x_{k+1}-x_k}{\alpha_k}=v_{k+1}-x_{k+1},
	\label{eq:sym-ABPD-PGS-x}\\		
	{}	\frac{y_{k+1}-y_k}{\alpha_{k}}=\eta_k(w_{k+1}-y_{k+1}),
	\label{eq:sym-ABPD-PGS-y}\\
	{}	\gamma_k\frac{\nabla\phi(v_{k+1})-\nabla\phi(v_k)}{\alpha_k}\in
	\sgf  (\nabla\phi(\bar x_{k})-\nabla\phi(v_{k+1}))-\mathcal F(\bar x_k,w_k,v_{k+1}),
	\label{eq:sym-ABPD-PGS-v}\\
	{}	\beta_k\frac{\nabla\psi(w_{k+1})-\nabla\psi(w_k)}{\alpha_k} \in\sgg  ( \nabla\psi(\bar y_{k})-\nabla\psi(w_{k+1}))-\eta_k\mathcal G(\bar y_k,\bar v_{k+1},w_{k+1}),
	\label{eq:sym-ABPD-PGS-w}
\end{subnumcases}
where $\eta_k$ and $\bar v_{k+1}$ are defined by \cref{eq:etak-b-vk1} and  
\[
\begin{aligned}
	\mathcal F(\bar x_k,w_k,v_{k+1})={}& \nabla f_1(\bar x_k)+\partial f_2(v_{k+1})+A^* w_k,\quad \bar x_k = \frac{x_k+\alpha_kv_k}{1+\alpha_k},\\
	\mathcal G(\bar y_k,\bar v_{k+1},w_{k+1})={}& \nabla g_1(\bar y_k)+\partial g_2(w_{k+1})-A\bar v_{k+1},\quad \bar y_k = \frac{y_k+\eta_k\alpha_kw_k}{1+\eta_k\alpha_k}.
\end{aligned}
\]
Note that $(\bar{x}_k,\bar{y}_k)$ are convex combinations of $(x_k,y_k)$ and $(v_k,w_k)$. The parameter updating  relation \cref{eq:ABPD-PS-gk-bk} leaves unchanged.  

We take the step size $\alpha_k^2(L_f \beta_k+L_g\gamma_k\tau_0+\opnm{A}^2)=\gamma_{k}\beta_k$ with 
\begin{equation}\label{eq:tau0}
	\tau_0:=\Big(1+\sqrt{\gamma_0\beta_0}/\sqrt{\opnm{A}^2+L_f\beta_0+L_g\gamma_0}\Big)^2\geq 1,
\end{equation}
and rearrange \cref{eq:sym-ABPD-PGS} as a primal-dual proximal gradient splitting method and list it below in \cref{algo:ABPD-PGS}, which is called Accelerated Bregman Primal-Dual Proximal Gradient Splitting (ABPD-PGS). It is not hard to find that, with the initial setting: $x_0,v_0\in\dom f_2$ and  $y_0,w_0\in\dom g_2$, we have $\bar{x}_k,x_k,v_k\in\dom f_2$ and  $\bar{y}_k,y_k,w_k\in\dom g_2$. 

\begin{algorithm}[H]
	\caption{ABPD-PGS for solving \cref{eq:minmax-L} under \cref{assum:ABPD-PS-prox-fun,assum:sym-ABPD-PGS} }
	\label{algo:ABPD-PGS} 
	\begin{algorithmic}[1] 
		\REQUIRE  Problem parameters: $\sgf   ,\,\sgg   \geq 0,\,L_f ,\,L_g $ and $\opnm{A}$.			\\
		\quad ~~~Initial guesses: $x_0,v_0\in\dom f_2$ and  $y_0,w_0\in\dom g_2$.\\
		\quad ~~~Initial parameters: $\gamma_0,\,\beta_0>0$ and $\tau_0$ given by \cref{eq:tau0}.
		\FOR{$k=0,1,\ldots$} 
		\STATE Compute $\alpha_k = \sqrt{\gamma_k\beta_k/(L_f \beta_k+L_g\gamma_k\tau_0+\opnm{A}^2)}$. \label{line:Sym-ABPD-PGS-ak}
		\STATE Update $\gamma_{k+1} = (\sgf  \alpha_k+\gamma_k)/(1+\alpha_k)$ and $
		\beta_{k+1}=(\sgg  \alpha_k+\beta_k)/(1+\alpha_k)$.		
		\STATE Set $\alpha_{k+1}= \sqrt{\gamma_{k+1}\beta_{k+1}/(L_f \beta_{k+1}+L_g\gamma_{k+1}\tau_0+\opnm{A}^2)}$ and $	\eta_k=\frac{\alpha_{k+1}}{\alpha_k}(1+\alpha_k)$. 
		\STATE
		\STATE Set $ \bar x_k = \frac{x_k+\alpha_kv_k}{1+\alpha_k}$ and $\widetilde{v}_k = \nabla f_1(\bar x_k)+A^* w_k-\sgf   ( \nabla\phi(\bar x_{k})-\nabla\phi(v_{k}))$.		
		\STATE Update $v_{k+1}=\mathop{\argmin}\limits_{v\in \mathcal X}\!\!\left\{
		f_2(v)+\dual{\widetilde{v}_k,v}+\frac{\gamma_k+\sgf    \alpha_k}{\alpha_k}D_\phi(v,v_k)\right\}$. 		
		\label{line:sym-ABPD-PGS-vk1}
		\STATE Update $ x_{k+1}=( x_k+\alpha_kv_{k+1})/(1+\alpha_k)$ and  $\bar v_{k+1} = v_{k+1}+(v_{k+1} -v_k)/\eta_k$. 
		\STATE
		\STATE Set  $\bar y_k = \frac{y_k+\eta_k\alpha_kw_k}{1+\eta_k\alpha_k}$ and $\widetilde{w}_k=\nabla g_1(\bar y_k)-A\bar v_{k+1}-\sgg  /\eta_k( \nabla\psi(\bar y_{k})-\nabla\psi(w_{k}))$.
		\STATE Update  $w_{k+1}=\mathop{\argmin}\limits_{w\in \mathcal Y}\!\!\left\{
		g_2(w)+\dual{ \widetilde{w}_k,w}+\frac{\beta_k+\sgg   \alpha_k}{\eta_k\alpha_k}D_\psi(w,w_k)\right\}$. 	 			 \label{line:sym-ABPD-PGS-wk1}
		\STATE Update $ y_{k+1}=( y_k+\alpha_k\eta_kw_{k+1})/(1+\alpha_k\eta_k)$.
		\ENDFOR
	\end{algorithmic}
\end{algorithm}
\begin{rem}\label{rem:ABPD-PGS}
	We mention that if the smooth components $f_1$ and $f_2$ vanish, then \cref{eq:sym-ABPD-PGS,algo:ABPD-PGS} reduce to \cref{eq:ABPD-PS,algo:ABPD-PS}, respectively. Moreover, if we choose quadratic prox-functions in \cref{eq:phi-R-psi-S}, then similarly with \cref{rem:xk1-yk1-eulidean}, the subproblems in Lines \ref{line:sym-ABPD-PGS-vk1} and \ref{line:sym-ABPD-PGS-wk1} of \cref{algo:ABPD-PGS} agree with standard preconditioned proximal calculations:
	\begin{align}
		\label{eq:ABPD-PGS-vk1}
		v_{k+1}={}&\mathop{\argmin}\limits_{v\in \mathcal X}\!\!\left\{
		f_2(v)+\dual{\widetilde{v}_k ,v}+\frac{\gamma_k+\sgf   \alpha_k}{2\alpha_k}\nm{v-v_k}_{R}^2\right\},\\
		\label{eq:ABPD-PGS-wk1}
		w_{k+1}={}&\mathop{\argmin}\limits_{w\in \mathcal Y}\!\!\left\{
		g_2(w)+\dual{ \widetilde{w}_k ,w}+\frac{\beta_k+\sgg  \alpha_k}{2\alpha_k\eta_k}\nm{w-w_k}_{S}^2\right\}.		 
	\end{align}
\end{rem}
\begin{lem}\label{lem:ABPD-PGS-est-i2}
	Let	$\{\Theta_k=(\gamma_k,\beta_k)\}\subset\R_+^2$  and $\{Z_k=(x_k,v_k,y_k,w_k)\}\subset \mathcal Z$ be generated by \cref{algo:ABPD-PGS}. Then under \cref{assum:sym-ABPD-PGS} it holds that
	\begin{equation}\label{eq:ABPD-PGS-est-i2}
		\begin{aligned}
			%
			\mathbb I_{21}+	\mathbb I_{22}\leq & -\alpha_k\dual{A^*(w_k-\widehat{y}),v_{k+1}-\widehat{x}}
			-\gamma_k D_\phi(v_{k+1},v_{k})+\frac{L_f\alpha_k^2}{2(1+\alpha_k)}\nm{v_{k+1}-v_k}^2,
		\end{aligned}
	\end{equation}
	where $\mathbb I_{21}$ and $\mathbb I_{22}$ are defined in \cref{eq:diff-Hk} with any $\widehat{\bm z}=(\widehat{x},\widehat{y})\in\dom \mathcal L$.
\end{lem}
\begin{proof}
	In view of \cref{eq:hat-f-hat-g} and \eqref{eq:sym-ABPD-PGS-x}, it holds that   
	\[
	\begin{aligned}
		\mathbb I_{21}
		=	{}&	f(x_{k+1})-f(x_k)+	\alpha_k\left(f(x_{k+1})-f(\widehat{x})\right) \\
		{}&\quad + \dual{A^*\widehat{y},x_{k+1}- x_{k}}+\alpha_k \dual{A^*\widehat{y},x_{k+1}-\widehat{x}}\\
		=
		{}&	f(x_{k+1})-f(x_k)+	\alpha_k\left(f(x_{k+1})-f(\widehat{x})\right) +  \alpha_k \dual{A^*\widehat{y},v_{k+1}-\widehat{x}} .
	\end{aligned}
	\]
	In view of the composite structure $f = f_1+f_2$, let us establish the upper bounds for $f_1$ and $f_2$ separately. By \cref{assum:sym-ABPD-PGS,lem:L-sg-descent}, we have 
	\begin{align}
		{}&		\alpha_k\left(f_1(x_{k+1})-f_1(\widehat{x})\right)+f_1(x_{k+1})-f_1(x_k)\notag \\ 
		\leq{}  & 
		\alpha_k\left(	\dual{\nabla f_1(\bar x_{k}), x_{k+1} -\widehat{x}}-\sgf   D_\phi(\widehat{x},\bar x_{k})+\frac{L_f}{2}\nm{x_{k+1}-\bar x_k}^2\right)\notag \\
		{}&\quad+
		\dual{\nabla f_1(\bar x_{k}), x_{k+1} -x_k}+ \frac{L_f}{2}\nm{x_{k+1}-\bar x_k}^2 \notag \\
		\overset{\text{by \eqref{eq:sym-ABPD-PGS-x}}}{=}{} &\alpha_k	\dual{\nabla f_1(\bar x_k), v_{k+1} -\widehat{x}}-\sgf\alpha_kD_\phi(\widehat{x},\bar x_{k}) +\frac{L_f(1+\alpha_k)}{2}\nm{x_{k+1}-\bar x_k}^2.\label{eq:ABPD-PGS-est-f1}
	\end{align}
	For the nonsmooth part $f_2$, we use the convexity property to obtain
	\[
	\begin{aligned}
		{}&	\alpha_k\left(f_2(x_{k+1})-f_2(\widehat{x})\right)
		+f_2(x_{k+1})-f_2(x_{k})\\
		\overset{\text{by \eqref{eq:sym-ABPD-PGS-x} }}{=}	{}&	\alpha_k\left(f_2 \left(\frac{x_k+\alpha_kv_{k+1}}{1+\alpha_k}\right)-f_2(\widehat{x})\right)
		+f_2\left(\frac{x_k+\alpha_kv_{k+1}}{1+\alpha_k}\right)-f_2(x_{k})\\
		\leq 	{}&	\alpha_k\left(\frac{1}{1+\alpha_k}f_2(x_k)+\frac{\alpha_k}{1+\alpha_k}f_2(v_{k+1}) -f_2(\widehat{x})\right) \\
		{}&\quad +\frac{1}{1+\alpha_k}f_2(x_k)+\frac{\alpha_k}{1+\alpha_k}f_2(v_{k+1})  -f_2(x_{k})\\
		={}&\alpha_k\left(f_2(v_{k+1}) -f_2(\widehat{x})\right).
	\end{aligned}
	\] 
	Observing \eqref{eq:sym-ABPD-PGS-v}, we find $p_{k+1}\in \partial f_2(v_{k+1})$, where 
	\[
	p_{k+1}:=\gamma_k\frac{\nabla\phi(v_k)-\nabla\phi(v_{k+1})}{\alpha_k}+
	\sgf   (\nabla\phi(\bar x_{k})-\nabla\phi(v_{k+1}))-\nabla f_1(\bar x_k)-A^* w_k .
	\]
	Hence, it follows immediately that 
	\begin{equation}\label{eq:ABPD-PGS-est-f2}
		\begin{aligned}
			{}&	\alpha_k\left(f_2(x_{k+1})-f_2(\widehat{x})\right)
			+f_2(x_{k+1})-f_2(x_{k}) 
			\leq{} \alpha_k\dual{p_{k+1},v_{k+1}-\widehat{x}}\\
			={}&\gamma_k\dual{ \nabla\phi(v_k)-\nabla\phi(v_{k+1}),v_{k+1}-\widehat{x}}+\sgf\alpha_k\dual{ 
				\nabla\phi(\bar x_{k})-\nabla\phi(v_{k+1}),v_{k+1}-\widehat{x}}\\
			{}&\quad -\alpha_k\dual{  \nabla f_1(\bar x_k)+ A^* w_k ,v_{k+1}-\widehat{x}}.
		\end{aligned}
	\end{equation}
	Applying \cref{lem:3-id} gives
	\[
	\begin{aligned}
		\sgf   \alpha_k\dual{\nabla\phi(\bar x_{k})-\nabla\phi(v_{k+1}), v_{k+1} -\widehat{x}} 
		=	{}&  \sgf   \alpha_k \left( D_\phi(\widehat{x},\bar x_{k})-D_\phi(\widehat{x},v_{k+1})-D_\phi(v_{k+1},\bar x_{k})\right),\\
		\gamma_k\dual{\nabla\phi(v_{k})-\nabla\phi(v_{k+1}), v_{k+1} -\widehat{x}} 
		=	{}& \gamma_k \left( D_\phi(\widehat{x},v_{k})-D_\phi(\widehat{x},v_{k+1})-D_\phi(v_{k+1},v_{k})\right).
	\end{aligned}
	\]
	Combining this with \cref{eq:ABPD-PGS-est-f1,eq:ABPD-PGS-est-f2}, we get
	\[
	\begin{aligned}
		\mathbb I_{21} 
		\leq & -\gamma_k D_\phi(v_{k+1},v_{k})+\frac{L_f(1+\alpha_k)}{2}\nm{x_{k+1}-\bar x_k}^2-\alpha_k\dual{  A^*( w_k-\widehat{y}) ,v_{k+1}-\widehat{x}}\\
		{}&\quad  -\sgf   \alpha_k  D_\phi(v_{k+1},\bar x_{k}) +\underbrace{\gamma_k  D_\phi(\widehat{x},v_{k})-(\gamma_k +\mu_f\alpha_k)D_\phi(\widehat{x},v_{k+1})}_{=-\mathbb I_{22}}.
	\end{aligned}
	\]
	In view of \eqref{eq:sym-ABPD-PGS-x} and the setting
	$\bar x_k = (x_k+\alpha_kv_k)/(1+\alpha_k)$, we have $x_{k+1}-\bar x_k = \alpha_k(v_{k+1}-v_k)/(1+\alpha_k)$.  Since  $\gamma_k+\sgf\alpha_k=\gamma_{k+1}(1+\alpha_k)$, we can move the second line to the left hand side to get \cref{eq:ABPD-PGS-est-i2}. 	This  completes the proof of this lemma.
\end{proof}
\begin{lem}\label{lem:sym-ABPD-PGS-est-i3}
	Let	$\{\Theta_k=(\gamma_k,\beta_k)\}\subset\R_+^2$  and $\{Z_k=(x_k,v_k,y_k,w_k)\}\subset \mathcal Z$ be generated by \cref{algo:ABPD-PGS}. If  $\gamma_0\geq \mu_f$ and $\beta_0\geq \mu_g$, then under \cref{assum:sym-ABPD-PGS} we have 
	\begin{equation}\label{eq:sym-ABPD-PGS-est-i3} 
		\begin{aligned}
			%
			\mathbb I_{31}+	\mathbb I_{32}\leq{} & \alpha_k\eta_k\dual{A(\bar v_{k+1}-\widehat{x}),w_{k+1}-\widehat{y}}
			-\beta_k D_\psi(w_{k+1},w_{k})  +\frac{L_g\alpha_k^2\eta_k^2}{2(1+\alpha_k\eta_k)}\nm{w_{k+1}-w_k}^2,
		\end{aligned}
	\end{equation}
	where $\mathbb I_{31}$ and $\mathbb I_{32}$ are defined in \cref{eq:diff-Hk} with $\widehat{\bm z}=(\widehat{x},\widehat{y})$ being a saddle point to \cref{eq:minmax-L}.
\end{lem}
\begin{proof}
	
	The proof is almost in line with that of \cref{lem:ABPD-PGS-est-i2} but with subtle difference for handling the additional constant $\eta_k$ in \eqref{eq:sym-ABPD-PGS-y} and \eqref{eq:sym-ABPD-PGS-w}. To move on, let us first prove $\eta_k\geq 1$.  Since $\alpha_k^2(L_f \beta_k+L_g\gamma_k\tau_0+\opnm{A}^2)=\gamma_{k}\beta_k$ (cf. Line \ref{line:Sym-ABPD-PGS-ak} of  \cref{algo:ABPD-PGS}), using \cref{eq:gk1-bk1} gives
	\begin{equation}\label{eq:etak-pgs}
		\begin{aligned}
			\eta_k ={}& \frac{1+\alpha_k}{\alpha_k}\alpha_{k+1}=\frac{1+\alpha_k}{\alpha_k}\frac{\sqrt{\gamma_{k+1}\beta_{k+1}}}{\sqrt{L_f \beta_{k+1}+L_g\gamma_{k+1}\tau_0+\opnm{A}^2}}\\ \geq{}& \frac{\sqrt{(\gamma_{k}+\sgf  \alpha_k)(\beta_{k}+\sgg  \alpha_k)}}{\alpha_k\sqrt{L_f \beta_{k}+L_g\gamma_{k}\tau_0+\opnm{A}^2}} ={} \sqrt{\left(1+\sgf \alpha_k/\gamma_k\right)\left(1+\sgg  \alpha_k/\beta_k\right)}\geq 1.
		\end{aligned}
	\end{equation}
	Above, we used the monotone relations $\gamma_{k}\geq \gamma_{k+1}$ and $\beta_{k}\geq\beta_{k+1}$, which holds true by \cref{eq:gk-bk-lb-ub} and our assumption: $\gamma_0\geq \mu_f$ and $\beta_0\geq \mu_g$. 
	
	On the other hand, by \cref{eq:hat-f-hat-g} and \eqref{eq:sym-ABPD-PGS-y} and \cref{lem:L-gap-obj-res} (ii), it is clear that 
	\[
	\begin{aligned}
		\mathbb I_{31}=	{}&	\widehat{g}(y_{k+1};\widehat{z})-\widehat{g}(y_k;\widehat{z})+	\alpha_k\left(\widehat{g}(y_{k+1};\widehat{z})-g(\widehat{y})\right)\\
		\leq {}&	\widehat{g}(y_{k+1};\widehat{z})-\widehat{g}(y_k;\widehat{z})+	\alpha_k\eta_k\left(\widehat{g}(y_{k+1};\widehat{z})-g(\widehat{y})\right)\\
		=	{}&	g(y_{k+1})-g(y_k)+	\alpha_k\eta_k\left(g(y_{k+1})-g(\widehat{y})\right) - \dual{A \widehat{x},y_{k+1}- y_{k}}-\alpha_k\eta_k \dual{A \widehat{x},y_{k+1}-\widehat{y}}\\
		=
		{}&	g(y_{k+1})-g(y_k)+	\alpha_k\eta_k\left(g(y_{k+1})-g(\widehat{y})\right) -  \alpha_k\eta_k \dual{A \widehat{x},w_{k+1}-\widehat{y}}.
	\end{aligned}
	\]
	In the following, we consider $g_1$ and $g_2$ subsequently. By \cref{lem:L-sg-descent,assum:sym-ABPD-PGS} and \eqref{eq:sym-ABPD-PGS-y}, we have 
		\begin{align}
			{}&		\alpha_k\eta_k\left(g_1(y_{k+1})-g_1(\widehat{y})\right)+g_1(y_{k+1})-g_1(y_k)\notag \\ 
			\leq{}  & 
			\alpha_k\eta_k\left(	\dual{\nabla g_1(\bar y_{k}), y_{k+1} -\widehat{y}}-\sgg  D_\psi(\widehat{y},\bar y_{k})+\frac{L_g}{2}\nm{y_{k+1}-\bar y_k}^2\right)\notag \\
			{}&\quad+
			\dual{\nabla g_1(\bar y_{k}), y_{k+1} -y_k}+ \frac{L_g}{2}\nm{y_{k+1}-\bar y_k}^2 \notag \\
			={} &\alpha_k\eta_k	\dual{\nabla g_1(\bar y_k), w_{k+1} -\widehat{y}}-\sgg\alpha_k\eta_kD_\psi(\widehat{y},\bar y_{k}) +\frac{L_g(1+\alpha_k\eta_k)}{2}\nm{y_{k+1}-\bar y_k}^2.\label{eq:est-g1-sym-ABPD-PGS}
		\end{align}
	Also, it follows that
	\[
	\begin{aligned}
		{}&	\alpha_k\eta_k\left(g_2(y_{k+1})-g_2(\widehat{y})\right)
		+g_2(y_{k+1})-g_2(y_{k})\\
		\overset{\text{by \eqref{eq:sym-ABPD-PGS-y} }}{=}	{}&	\alpha_k\eta_k\left(g_2 \left(\frac{y_k+\alpha_k\eta_kw_{k+1}}{1+\alpha_k\eta_k}\right)-g_2(\widehat{y})\right)
		+g_2\left(\frac{y_k+\alpha_k\eta_kw_{k+1}}{1+\alpha_k\eta_k}\right)-g_2(y_{k})\\
		\leq 	{}&	\alpha_k\eta_k\left(\frac{1}{1+\alpha_k\eta_k}g_2(y_k)+\frac{\alpha_k\eta_k}{1+\alpha_k\eta_k}g_2(w_{k+1}) -g_2(\widehat{y})\right) \\
		{}&\quad +\frac{1}{1+\alpha_k\eta_k}g_2(y_k)+\frac{\alpha_k\eta_k}{1+\alpha_k\eta_k}g_2(w_{k+1})  -g_2(y_{k})\\
		={}&\alpha_k\eta_k\left(g_2(w_{k+1}) -g_2(\widehat{y})\right).
	\end{aligned}
	\] 
	Observing \eqref{eq:sym-ABPD-PGS-w}, we find $q_{k+1}\in \partial g_2(w_{k+1})$, where 
	\[
	q_{k+1}:=\frac{1}{\eta_k}\left[\beta_k\frac{\nabla\psi(w_k)-\nabla\psi(w_{k+1})}{\alpha_k}+
	\sgg   (\nabla\psi(\bar y_{k})-\nabla\psi(w_{k+1}))\right]-\nabla g_1(\bar y_k)+ A\bar v_{k+1}.
	\]
	Hence, we obtain
	\begin{equation}\label{eq:est-g2-sym-ABPD-PGS}
		\begin{aligned}
			{}&	\alpha_k\eta_k\left(g_2(y_{k+1})-g_2(\widehat{y})\right)
			+g_2(y_{k+1})-g_2(y_{k}) 
			\leq{} \alpha_k\eta_k\dual{q_{k+1},w_{k+1}-\widehat{y}}\\
			={}&\beta_k\dual{ \nabla\psi(w_k)-\nabla\psi(w_{k+1}),w_{k+1}-\widehat{y}}+\sgg\alpha_k\dual{ 
				\nabla\psi(\bar y_{k})-\nabla\psi(w_{k+1}),w_{k+1}-\widehat{y}}\\
			{}&\quad -\alpha_k\eta_k\dual{  \nabla g_1(\bar y_k)- A \bar v_{k+1} ,w_{k+1}-\widehat{y}}.
		\end{aligned}
	\end{equation}
	Furthermore, invoking the tree-term identity \cref{eq:3-id} yields that
	\[
	\begin{aligned}
		\sgg   \alpha_k\dual{\nabla\psi(\bar y_{k})-\nabla\psi(w_{k+1}), w_{k+1} -\widehat{y}} 
		=	{}&	  \sgg   \alpha_k \left( D_\psi(\widehat{y},\bar y_{k})-D_\psi(\widehat{y},w_{k+1})-D_\psi(w_{k+1},\bar y_{k})\right),\\	 \beta_k\dual{\nabla\psi(w_{k})-\nabla\psi(w_{k+1}), w_{k+1} -\widehat{y}} 
		=	{}&	 \beta_k \left( D_\psi(\widehat{y},w_{k})-D_\psi(\widehat{y},w_{k+1})-D_\psi(w_{k+1},w_{k})\right).
	\end{aligned}
	\]
	Combining these two estimates with \cref{eq:est-g1-sym-ABPD-PGS,eq:est-g2-sym-ABPD-PGS} gives 
	\[
	\begin{aligned}
		\mathbb I_{31}
		\leq {} &\underbrace{\beta_k  D_\psi(\widehat{y},w_{k})-(\beta_k +\mu_f\alpha_k)D_\psi(\widehat{y},w_{k+1})}_{=-	\mathbb I_{32}}+\alpha_k\eta_k\dual{  A (\bar v_{k+1}-\widehat{x}) ,w_{k+1}-\widehat{y}} \\
		{}&\quad -\beta_k D_\psi(w_{k+1},w_{k})+\frac{L_g(1+\alpha_k\eta_k)}{2}\nm{y_{k+1}-\bar y_k}^2  -\sgg   \alpha_k  D_\psi(w_{k+1},\bar y_{k}).
	\end{aligned}
	\]
	Notice that by \eqref{eq:sym-ABPD-PGS-y} and the setting
	$\bar y_k = (y_k+\alpha_k\eta_kw_k)/(1+\alpha_k\eta_k)$, we have $y_{k+1}-\bar y_k = \alpha_k\eta_k(w_{k+1}-w_k)/(1+\alpha_k\eta_k)$.
	Since  $\beta_k+\sgg\alpha_k=\beta_{k+1}(1+\alpha_k)$, we can move the first line to the left hand side to get \cref{eq:sym-ABPD-PGS-est-i3}. This  completes the proof of \cref{lem:sym-ABPD-PGS-est-i3}.
	%
\end{proof}
\begin{thm}
	\label{thm:conv-sym-ABPD-PGS} 
	Let	$\{\Theta_k=(\gamma_k,\beta_k)\}\subset\R_+^2$  and $\{Z_k=(x_k,v_k,y_k,w_k)\}\subset \mathcal Z$ be generated by \cref{algo:ABPD-PGS}. If  $\gamma_0\geq \mu_f$ and $\beta_0\geq \mu_g$, then under \cref{assum:ABPD-PS-prox-fun,assum:sym-ABPD-PGS} we have 
	\begin{equation}\label{eq:diff-sym-ABPD-PGS}
		\mathcal H_{k+1}(\widehat{\bm z})-\mathcal H_k(\widehat{\bm z})+\alpha_k\mathcal H_{k+1}(\widehat{\bm z}) \leq 0,
	\end{equation}
	where  $\mathcal H_k(\widehat{\bm z})$ is define by \cref{eq:Hk} with $\widehat{\bm z}=(\widehat{x},\widehat{y})$ being a saddle point to \cref{eq:minmax-L}.  
\end{thm}
\begin{proof}
		Combining \cref{eq:diff-Hk} with  \cref{lem:ABPD-PGS-est-i2,lem:sym-ABPD-PGS-est-i3} gives
		\[
		\begin{aligned} 
			\mathcal H_{k+1}(\widehat{\bm z})-\mathcal H_k(\widehat{\bm z})+\alpha_k\mathcal H_{k+1}(\widehat{\bm z})    
			\leq{}& \frac{L_f\alpha_k^2}{2(1+\alpha_k)}\nm{v_{k+1}-v_k}^2+\frac{L_g\alpha_k^2\eta_k^2}{2(1+\alpha_k\eta_k)}\nm{w_{k+1}-w_k}^2\\
			{}&\quad  +\Delta_k
			-	\gamma_{k}
			D_\phi(v_k,v_{k+1})- \beta_{k} 
			D_\psi(w_k, w_{k+1}),
		\end{aligned}
		\]
		where $\Delta_k$ is the same as \cref{eq:Deltak} and can be simplified as $	\Delta_k=\alpha_k\dual{A(v_{k+1}-v_k), w_{k+1}- w_k}$ by \cref{eq:etak-b-vk1}. Since $\gamma_0 \geq \mu_f$ and $\beta_0\geq\mu_g$, it follows from \cref{eq:gk1-bk1} that $\{\gamma_k\}$ and $\{\beta_k\}$ are non-increasing. According to our choice $\gamma_k\beta_k =  \alpha_k^2\big(\opnm{A}^2+L_f \beta_k+L_g \gamma_k\tau_0\big)$ (cf. Line \ref{line:Sym-ABPD-PGS-ak} of  \cref{algo:ABPD-PGS}), $\{\alpha_k\}$ is also non-increasing and thus by \cref{eq:tau0} we have  $\eta_k=\alpha_{k+1}(1+1/\alpha_k)\leq 1+\alpha_k\leq 1+\alpha_0$.
		In view of \cref{eq:tau0}, it is easy to obtain $1+\alpha_0 \leq\sqrt{ \tau_0}$. This together with the fact $\eta_k\geq1$ (cf.\cref{eq:etak-pgs}) implies immediately that
		\[
		\begin{aligned} 
			{}&	\mathcal H_{k+1}(\widehat{\bm z})-\mathcal H_k(\widehat{\bm z})+\alpha_k\mathcal H_{k+1}(\widehat{\bm z})    \\
			\leq{}& \Delta_k-	\frac{\gamma_k (1+\alpha_k)-L_f  \alpha_k^2}{2(1+\alpha_k)}\nm{v_{k+1}-v_k}^2-\frac{\beta_k(1+\alpha_k)-L_g \tau_0\alpha_k^2}{2(1+\alpha_k)}\nm{w_{k+1}-w_k}^2.
		\end{aligned}
		\]
		Observe the relations $\gamma_k (1+\alpha_k)>L_f  \alpha_k^2,\, \beta_k(1+\alpha_k)>L_g \tau_0\alpha_k^2$ and 
		\[
		\left(\gamma_k (1+\alpha_k)-L_f  \alpha_k^2\right)\left( \beta_k(1+\alpha_k)-L_g \tau_0\alpha_k^2\right)\geq \opnm{A}^2\alpha_k^2(1+\alpha_k)^2.
		\]
		Therefore, by using $	\Delta_k\leq \alpha_k\opnm{A}\nm{v_{k+1}-v_k}\nm{w_{k+1}-w_k}$ and the basic inequality $ab\leq (a^2+b^2)/2$, we obtain \cref{eq:diff-sym-ABPD-PGS} and conclude the proof. 
	\end{proof}
	\subsection{Convergence rate and iteration complexity}
	\begin{thm}
		\label{thm:rate-sym-ABPD-PGS}
		Let $\{(x_k,y_k)\}\subset\dom \mathcal L$ be generated by  \cref{algo:ABPD-PGS}.  Then under \cref{assum:ABPD-PS-prox-fun,assum:sym-ABPD-PGS}, we have 	\begin{equation}\label{eq:L-gap-x-y-sym-ABPD-PGS}
			\mathcal L(x_k,\widehat{y})-	\mathcal L(\widehat{x},y_k)+	\frac{\sgf  }{2}\nm{x_k-\widehat{x}}^2+\frac{\sgg }{2}\nm{y_k-\widehat{y}}^2\leq 2\theta_k \mathcal H_0(\widehat{\bm z}),
		\end{equation}
		where $\{\theta_k\}$ is defined by \cref{eq:thetak} and  $\widehat{\bm z}=(\widehat{x},\widehat{y})$ is a saddle point to \cref{eq:minmax-L}.
		If  $\gamma_0\geq \sgf ,\,\beta_0\geq\sgg  $ and $\gamma_0\beta_0\leq L_f \beta_{0}+L_g \gamma_0+\opnm{A}^2$, then we have the following decay estimates of $\theta_k$.
		\begin{itemize}
			\item If  $\sgf =\sgg  = 0$, then it holds that
			\begin{equation}\label{eq:ABPD-PGS-rate1}
				\theta_k\leq a_k:=\frac{64(L_f \beta_0+4L_g\gamma_0) }{\gamma_0\beta_0k^2}+ \frac{4\opnm{A} }{\sqrt{\gamma_0\beta_0}k}.
			\end{equation}
			\item If $\sgf+\sgg>0$, then we have
			\begin{equation}\label{eq:ABPD-PGS-rate2}\small
				\theta_k\leq c_k:=\min\left\{a_k,\,b_k,\,\exp\left(-\frac{k}{8}\sqrt{\frac{\sgg  }{L_g}}\right)+ \frac{64(L_f \sgg  +\opnm{A}^2)}{\sgg  \gamma_0k^2}\right\},
			\end{equation}
			where 
			\[
			b_k:=\exp\left(-\frac{k}{4}\sqrt{\frac{\sgf  }{L_f }}\right)+ \frac{64(4L_g \sgf  +\opnm{A}^2)}{\sgf  \beta_0k^2}.
			\]
			\item If $\sgf   \sgg >0$, then $		\theta_k\leq \min\left\{c_k,\,\big(1+\sqrt{\sgf  \sgg  }/\sqrt{L_f \sgg  +4L_g \sgf  +\opnm{A}^2}\big)^{-k}\right\}$.
		\end{itemize}
	\end{thm}
	\begin{proof}
		Since $\opnm{A}^2\alpha_k^2\leq \gamma_k\beta_k$ (cf. Line \ref{line:Sym-ABPD-PGS-ak} of  \cref{algo:ABPD-PGS}), by \cref{assum:sym-ABPD-PGS} we obtain
		\[
		\alpha_k \snm{\dual{A(v_k-\widehat{x}), w_k-\widehat{y}}}\leq  \gamma_k D_\phi(\widehat{x},v_k)	+ \beta_k  D_\psi(\widehat{y}, w_k).
		\]
		Applying \cref{lem:decomp-diff-Hk,thm:conv-sym-ABPD-PGS}  proves \cref{eq:L-gap-x-y-sym-ABPD-PGS}. Recall the relation \cref{eq:diff-tk} for the sequence $\{\theta_k\}$. It remains to establish the decay estimate of $\theta_k$. In view of $\gamma_0\beta_0\leq L_f \beta_{0}+L_g \gamma_0+\opnm{A}^2$, it follows from \cref{eq:tau0} that $1\leq \tau_0\leq 4$. By \cref{eq:gk-bk-lb-ub}, we have $\gamma_k\leq\gamma_0$ and $\beta_k\leq \beta_0$ and thus  
		\[
		\alpha_k =\frac{ \sqrt{\gamma_{k}\beta_k}}{	\sqrt{{\Lf}\beta_k+L_g\gamma_k\tau_0+\opnm{A}^2}}  \leq 
		\frac{\sqrt{\gamma_0\beta_0}}{\sqrt{\Lf\beta_{0}+L_g\gamma_0+\opnm{A}^2}}\leq 1.
		\]
		Analogously to the proof of \cref{eq:L-gap-x-y}, according to the lower bounds of $\gamma_k$ and $\beta_k$, we consider four cases.
		\begin{itemize}
			\item[(i)] Firstly, for $\alpha_k\geq \sqrt{\gamma_0\beta_0}\theta_k/\sqrt{(L_f \beta_0+4L_g\gamma_0)\theta_{k}+\opnm{A}^2}$, it follows from \cref{eq:diff-tk,lem:app-est} that 
			\[
			\theta_k\leq\frac{64(L_f \beta_0+4L_g\gamma_0) }{\gamma_0\beta_0k^2}+ \frac{4\opnm{A} }{\sqrt{\gamma_0\beta_0}k}.
			\] 
			\item[(ii)] Secondly, for $\alpha_k\geq \sqrt{\sgf  \beta_0\theta_k}/\sqrt{L_f \beta_0\theta_{k}+4L_g \sgf  +\opnm{A}^2}$, we obtain
			\[
			\theta_k\leq\exp\left(-\frac{k}{4}\sqrt{\frac{\sgf  }{L_f }}\right)+ \frac{64(4L_g \sgf  +\opnm{A}^2)}{\sgf  \beta_0k^2}.
			\] 
			\item[(iii)] Thirdly, for $\alpha_k\geq \sqrt{\sgg  \gamma_0\theta_k}/\sqrt{4L_g \gamma_0\theta_k+L_f \sgg   +\opnm{A}^2}$, applying \cref{lem:app-est} again yields 
			\[
			\theta_k\leq\exp\left(-\frac{k}{8}\sqrt{\frac{\sgg  }{L_g}}\right)+ \frac{64(L_f \sgg  +\opnm{A}^2)}{\sgg  \gamma_0k^2}.
			\] 
			\item[(iv)] Finally, for $\alpha_k\geq \sqrt{\sgf  \sgg  }/\sqrt{L_f \sgg  +4L_g \sgf  +\opnm{A}^2}=\rho$, we obtain from \cref{eq:thetak} that $\theta_k\leq (1+\rho)^{-k}$.
		\end{itemize}
		Collection the above four estimates leads to the desired decay estimates. This finishes the proof.
	\end{proof}
	
	Thanks to \cref{thm:rate-sym-ABPD-PGS}, our \cref{algo:ABPD-PGS} also achieves the optimal iteration complexity \cite{ouyang_lower_2021,zhang_lower_2022} under the {\it pure first-order algorithm class} \cite[Definition 2.2]{zhang_lower_2022}. That is, to get  $	\mathcal L(x_k,\widehat{y})-	\mathcal L(\widehat{x},y_k)\leq \epsilon$, the iteration complexity is
	\[
	\begin{aligned}
		{}&	\mathcal O	\Bigg(\min\left\{\frac{\opnm{A}}{\epsilon},\,\frac{\opnm{A}}{\sqrt{\sgf  \epsilon}},\,\frac{\opnm{A}}{\sqrt{\sgg  \epsilon}},\,\frac{\opnm{A}}{\sqrt{ \sgf\sgg }}\snm{\ln \epsilon}\right\}\\
		{}&\qquad +\min\left\{\sqrt{\frac{L_f }{\epsilon}},\,\sqrt{\frac{L_f}{\mu_f}}\snm{\ln \epsilon}\right\}+\min\left\{\sqrt{\frac{L_g }{\epsilon}},\,\sqrt{\frac{L_g}{\mu_g}}\snm{\ln \epsilon}\right\}\Bigg).
	\end{aligned}
	\]
	When $\sgf  >0$, the iteration complexity for finding $\nm{x_k-\widehat{x}}^2\leq \epsilon$ is 
	\[
	\mathcal O	\left(\min\left\{\frac{\opnm{A}}{\sgf  \sqrt{\epsilon}},\,\frac{\opnm{A}}{\sqrt{\sgf  \sgg}}\snm{\ln \frac{\sgf  }{\epsilon}}\right\}+\sqrt{\frac{L_f}{\mu_f}}\snm{\ln \frac{\sgf  }{\epsilon}}+\min\left\{\sqrt{\frac{L_g }{\sgf  \epsilon}},\,\sqrt{\frac{L_g}{\mu_g}}\snm{\ln \frac{\sgf  }{\epsilon}}\right\}\right).
	\]
	When $\sgg  >0$, the iteration complexity for finding $\nm{y_k-\widehat{y}}^2\leq \epsilon$ is
	\[
	\mathcal O	\left(\min\left\{\frac{\opnm{A}}{\sgg  \sqrt{\epsilon}},\,\frac{\opnm{A}}{\sqrt{\sgf\sgg}}\snm{\ln \frac{\sgg  }{\epsilon}}\right\}+\sqrt{\frac{L_g}{\mu_g}}\snm{\ln \frac{\sgg  }{\epsilon}}+\min\left\{\sqrt{\frac{L_f }{ \sgg  \epsilon}},\,\sqrt{\frac{L_f}{\mu_f}}\snm{\ln \frac{\sgg  }{\epsilon}}\right\}\right).
	\]
	\section{Restart}
	\label{sec:restart} 
	Recently, the restart strategy has been applied to linear programming \cite{Applegate2023} and quadratic programming \cite{Chen2025d}. However, it is rare to see restart schemes for general saddle-point problems. In \cite{tran-dinh_smooth_2018}, Tran-Dinh et al. proposed a fixed restart for the proposed methods but without theoretical analysis. In this section, we introduce a fixed restart idea to \cref{algo:ABPD-PS,algo:ABPD-PGS} and prove the linear convergence results.
	
	For given $K\in\mathbb N_+$ and $\gamma,\beta>0,x,\,v\in\dom f,y,\,w\in\dom g$, let 
	\[
	(\gamma^+,\beta^+,Z^+) = \texttt{ABPD-PS}(\gamma,\beta,Z),\quad Z = (x,v,y,w),
	\]
	be the output of \cref{algo:ABPD-PS} after {\it one} iteration
	with the initial setting: $\gamma_0=\gamma,\beta_0=\beta,x_0=x,y_0=y,v_0=v,w_0=w$ and $\mu_f=\mu_g=0$. A restart variant of \cref{algo:ABPD-PS} is given in \cref{algo:ABPD-PS-restart}.
	\begin{algorithm}[H]
		\caption{ABPD-PS with Restart}
		\label{algo:ABPD-PS-restart} 
		\begin{algorithmic}[1] 
			\REQUIRE  Restart period $K\in\mathbb N_{\geq 2}$.			\\
			\quad ~~~Initial guesses: $x_0=v_0\in\dom f$ and  $y_0=w_0\in\dom g$.\\
			\quad ~~~Initial parameters: $\gamma_0,\,\beta_0>0$ and  problem parameter: $\opnm{A}$.
			\FOR{$j=1,\ldots$} 
			\STATE $n_{j,0} = (j-1)K,\,Z_{n_{j,0}} = (x_{n_{j,0}},v_{n_{j,0}},y_{n_{j,0}},w_{n_{j,0}})$
			\FOR{$i=0,1,\cdots,K-1$}
			\STATE $n_{j,i} = (j-1)K+i,\,n_{j,i+1} = (j-1)K+i+1$
			\STATE $(\gamma_{n_{j,i+1}},\beta_{n_{j,i+1}},Z_{n_{j,i+1}}) = \texttt{ABPD-PS}(\gamma_{n_{j,i}},\beta_{n_{j,i}},Z_{n_{j,i}})$
			\ENDFOR
			\STATE $\gamma_{n_{j,K}} = \gamma_0, \,x_{n_{j,K}}= x_{n_{j,K-1}},\,v_{n_{j,K}}= x_{n_{j,K-1}} $\label{algo:re-x}
			\STATE $ \beta_{n_{j,K}}  =\beta_0, \, y_{n_{j,K}}= y_{n_{j,K-1}},\,w_{n_{j,K}}= y_{n_{j,K-1}}$ \label{algo:re-y}
			\STATE $ Z_{n_{j,K}} = (x_{n_{j,K}},v_{n_{j,K}},y_{n_{j,K}},w_{n_{j,K}})$		
			\ENDFOR
		\end{algorithmic}
	\end{algorithm}
	\begin{thm}\label{thm:conv-ABPD-PS-R}
		Suppose \cref{assum:ABPD-PS-f-g} holds with $\mu_f\mu_g>0$. For \cref{algo:ABPD-PS-restart}, if we choose 
		\begin{equation}\label{eq:K}
			K
			=1+ \left\lceil\frac{2e\opnm{A}}{\sqrt{\gamma_0\beta_0}}\left(\frac{\mu_f\mu_g+2\mu_g\gamma_0+2\mu_f\beta_0}{\mu_f\mu_g}\right)\right\rceil,
		\end{equation}
		then for any $N\in\mathbb N_{\geq1 }$, we have 
		\begin{equation}\label{eq:linear-rate-ps-restart}
			\mathcal L(x_{N},\widehat{y})-	\mathcal L(\widehat{x},y_{N})\leq K 	e^{-N/K }\left[ \mathcal L(x_0, \widehat{y})-\mathcal L(\widehat{x}, y_0)\right],
		\end{equation}
		where $\widehat{\bm z}=(\widehat{x},\widehat{y})$ a saddle point to \cref{eq:minmax-L}.
	\end{thm}
	\begin{proof}
		According to the initial setting $x_0=v_0,y_0=w_0$ and the restart setting (cf.Lines \ref{algo:re-x} and \ref{algo:re-y} of \cref{algo:ABPD-PS-restart}), we claim that 
		\[
		v_{n_{j,0}} = x_{n_{j,0}},\quad 	w_{n_{j,0}} = y_{n_{j,0}},\quad \gamma_{n_{j,0}} = \gamma_0,\quad \beta_{n_{j,0}}=\beta_0,\quad \text{ for all }j\geq 1.
		\]  
		Thanks to \cref{thm:rate-ABPD-PS}, we have 
		\[
		\mathcal L(x_{n_{j,i}},\widehat{y})-	\mathcal L(\widehat{x},y_{n_{j,i}})\leq 	
		\frac{2\opnm{A}}{\opnm{A}+\sqrt{\gamma_{n_{j,0}}\beta_{n_{j,0}}}i}
		\left(\mathcal L(x_{n_{j,0}}, \widehat{y})-\mathcal L(\widehat{x}, y_{n_{j,0}})+	R_{j,0}(\widehat{\bm z})\right),
		\]
		for all $0\leq i\leq K-1$, where 
		\[
		\begin{aligned}
			R_{j,0}(\widehat{\bm z}) := \gamma_{n_{j,0}} D_{\phi}(\widehat{x},v_{n_{j,0}})
			+\beta_{n_{j,0}} D_\psi(\widehat{y},w_0)-\alpha_{n_{j,0}}\dual{A(v_{n_{j,0}}-\widehat{x}),w_{n_{j,0}}-\widehat{y}}.
		\end{aligned}
		\]
		In view of \cref{eq:cross}, it follows that
		\[
		R_{j,0}(\widehat{\bm z})\leq 2\gamma_{n_{j,0}} D_{\phi}(\widehat{x},v_{n_{j,0}})
		+2\beta_{n_{j,0}} D_\psi(\widehat{y},w_{n_{j,0}})=2\gamma_{n_{j,0}} D_{\phi}(\widehat{x},x_{n_{j,0}})
		+2\beta_{n_{j,0}} D_\psi(\widehat{y},y_{n_{j,0}}).
		\]
		Thanks to \cref{lem:L-gap-obj-res} (ii), it is clear that 
		\begin{equation}\label{eq:ub-H0}
			\begin{aligned}
				R_{j,0}(\widehat{\bm z}) \leq	{}& \frac{2\gamma_{n_{j,0}}}{\mu_f}\left[ \mathcal L(x_{n_{j,0}}, \widehat{y})-\mathcal L(\widehat{x}, y_{n_{j,0}})\right]
				+\frac{2\beta_{n_{j,0}}}{\mu_g}\left[ \mathcal L(x_{n_{j,0}}, \widehat{y})-\mathcal L(\widehat{x}, y_{n_{j,0}})\right]\\
				={}& \frac{2\mu_g\gamma_{n_{j,0}}+2\mu_f\beta_{n_{j,0}}}{\mu_f\mu_g} \left[ \mathcal L(x_{n_{j,0}}, \widehat{y})-\mathcal L(\widehat{x}, y_{n_{j,0}})\right].
			\end{aligned}
		\end{equation}
		Consequently, we obtain 
		\begin{equation}\label{eq:key-nji}
			\mathcal L(x_{n_{j,i}},\widehat{y})-	\mathcal L(\widehat{x},y_{n_{j,i}})\leq 	\frac{2\opnm{A}(\mu_f\mu_g+2\mu_g\gamma_0+2\mu_f\beta_0)}{\sqrt{\gamma_0\beta_0}\mu_f\mu_gi}   \left[ \mathcal L(x_{n_{j,0}}, \widehat{y})-\mathcal L(\widehat{x}, y_{n_{j,0}})\right],
		\end{equation}
		for all $j\geq 1$ and $1\leq i\leq K -1$. This implies a key estimate
		\begin{equation}\label{eq:key-nj}
			\begin{aligned}
				{}&	\mathcal L(x_{n_{j+1,0}},\widehat{y})-	\mathcal L(\widehat{x},y_{n_{j+1,0}})=
				\mathcal L(x_{n_{j,K }},\widehat{y})-	\mathcal L(\widehat{x},y_{n_{j,K }})\\
				={}&
				\mathcal L(x_{n_{j,K -1}},\widehat{y})-	\mathcal L(\widehat{x},y_{n_{j,K -1}})\leq 	\frac{C}{ K -1} \left[ \mathcal L(x_{n_{j,0}}, \widehat{y})-\mathcal L(\widehat{x}, y_{n_{j,0}})\right],
			\end{aligned}
		\end{equation}
		where 
		\begin{equation}\label{eq:C}
			C:=\frac{2\opnm{A}(\mu_f\mu_g+2\mu_g\gamma_0+2\mu_f\beta_0)}{\sqrt{\gamma_0\beta_0}\mu_f\mu_g} .
		\end{equation}
		
		Now for any $N\in\mathbb N_{\geq1}$, we have $N = \bar{j}K +\bar{i}$ with $\bar{j}\in\mathbb N$ and $1\leq \bar{i}\leq K -1$. In other words, to get $\{(x_N,y_N)\}$, we have already done restart operation exactly $\bar{j}$ times. Thus by \cref{eq:key-nj}, it follows that
		\[
		\begin{aligned}
			{}&	\mathcal L(x_{n_{\bar{j}+1,0}},\widehat{y})-	\mathcal L(\widehat{x},y_{n_{\bar{j}+1,0}})\leq \frac{C}{K -1}\mathcal L(x_{n_{\bar{j},0}},\widehat{y})-	\mathcal L(\widehat{x},y_{n_{\bar{j},0}})\\
			\leq{}& \cdots\leq \left(\frac{C}{K -1}\right)^{\bar{j}}\left[ \mathcal L(x_{n_{1,0}}, \widehat{y})-\mathcal L(\widehat{x}, y_{n_{1,0}})\right]
			=\left(\frac{C}{K -1}\right)^{\bar{j}}\left[ \mathcal L(x_{0}, \widehat{y})-\mathcal L(\widehat{x}, y_{0})\right].
		\end{aligned}
		\]
		In addition, in view of \cref{eq:key-nji}, we also have
		\[
		\begin{aligned}
			{}&	\mathcal L(x_{N},\widehat{y})-	\mathcal L(\widehat{x},y_{N})=	\mathcal L(x_{n_{\bar{j}+1,\bar{i}}},\widehat{y})-	\mathcal L(\widehat{x},y_{n_{\bar{j}+1,\bar{i}}})\leq \frac{C}{\bar{i}}\left[\mathcal L(x_{n_{\bar{j}+1,0}},\widehat{y})-	\mathcal L(\widehat{x},y_{n_{\bar{j}+1,0}})\right]\\
			\leq{}& C\left[\mathcal L(x_{n_{\bar{j}+1,0}},\widehat{y})-	\mathcal L(\widehat{x},y_{n_{\bar{j}+1,0}})\right]\leq C\left(\frac{C}{K -1}\right)^{\bar{j}}\left[ \mathcal L(x_{0}, \widehat{y})-\mathcal L(\widehat{x}, y_{0})\right].
		\end{aligned}
		\]
		Observing \cref{eq:K,eq:C}, it is clear that $K -1\geq eC$ and $\bar{j}\geq N/K -1$, which yields
		\[
		C\left(\frac{C}{K -1}\right)^{\bar{j}}\leq K e^{-N/K }.
		\]
		This leads to \cref{eq:linear-rate-ps-restart} and completes the proof of this theorem.
	\end{proof}
	
	Similarly, we can consider a restart variant of \cref{algo:ABPD-PGS}. To do this, for given $K\in\mathbb N_+$ and $\gamma,\beta>0,x,\,v\in\dom f,y,\,w\in\dom g$, let 
	\[
	(\gamma^+,\beta^+,Z^+) = \texttt{ABPD-PGS}(\gamma,\beta,Z),\quad Z = (x,v,y,w),
	\]
	be the output of \cref{algo:ABPD-PGS} after {\it one} iteration
	with the initial setting: $\gamma_0=\gamma,\beta_0=\beta,x_0=x,y_0=y,v_0=v,w_0=w$ and $\mu_f=\mu_g=0$. 
	\begin{algorithm}[H]
		\caption{ABPD-PGS with Restart}
		\label{algo:ABPD-PGS-restart} 
		\begin{algorithmic}[1] 
			\REQUIRE  Restart period $K\in\mathbb N_{\geq 2}$ and  problem parameters: $L_f,L_g,\,\opnm{A}$.			\\
			\quad ~~~Initial guesses: $x_0=v_0\in\dom f_2$ and  $y_0=w_0\in\dom g_2$.\\
			\quad ~~~Initial parameters: $\gamma_0,\,\beta_0>0,\,\tau_0=\Big(1+\sqrt{\gamma_0\beta_0}/\sqrt{\opnm{A}^2+L_f\beta_0+L_g\gamma_0}\Big)^2$.
			\FOR{$j=1,\ldots$} 
			\STATE $n_{j,0} = (j-1)K,\,Z_{n_{j,0}} = (x_{n_{j,0}},v_{n_{j,0}},y_{n_{j,0}},w_{n_{j,0}})$
			\FOR{$i=0,1,\cdots,K-1$}
			\STATE $n_{j,i} = (j-1)K+i,\,n_{j,i+1} = (j-1)K+i+1$
			\STATE $(\gamma_{n_{j,i+1}},\beta_{n_{j,i+1}},Z_{n_{j,i+1}}) = \texttt{ABPD-PGS}(\gamma_{n_{j,i}},\beta_{n_{j,i}},Z_{n_{j,i}})$
			\ENDFOR
			\STATE $\gamma_{n_{j,K}} = \gamma_0, \,x_{n_{j,K}}= x_{n_{j,K-1}},\,v_{n_{j,K}}= x_{n_{j,K-1}} $ 
			\STATE $ \beta_{n_{j,K}}  =\beta_0, \, y_{n_{j,K}}= y_{n_{j,K-1}},\,w_{n_{j,K}}= y_{n_{j,K-1}}$  
			\STATE $ Z_{n_{j,K}} = (x_{n_{j,K}},v_{n_{j,K}},y_{n_{j,K}},w_{n_{j,K}})$		
			\ENDFOR
		\end{algorithmic}
	\end{algorithm}
	
	Applying the proof of \cref{thm:conv-ABPD-PS-R}, we conclude the following.
	\begin{thm}\label{thm:conv-ABPD-PGS-R}
		Suppose \cref{assum:sym-ABPD-PGS}  holds with $\mu_f\mu_g>0$. For \cref{algo:ABPD-PGS-restart}, if we choose 
		\begin{equation}\label{eq:K-pgs}\small
			K
			=1+ \left\lceil\frac{16e\opnm{A}}{\sqrt{\gamma_0\beta_0}}\frac{C_0}{\mu_f\mu_g}+\frac{16\sqrt{e}\sqrt{L_f \beta_0+4L_g\gamma_0} }{\sqrt{\gamma_0\beta_0}}\sqrt{\frac{C_0}{\mu_f\mu_g}}\right \rceil,
		\end{equation}
		with $C _0:= \mu_f\mu_g+2\mu_g\gamma_0+2\mu_f\beta_0$,
		then for any $N\in\mathbb N_{\geq1 }$, we have 
		\begin{equation}\label{eq:linear-rate-pgs-restart}
			\mathcal L(x_{N},\widehat{y})-	\mathcal L(\widehat{x},y_{N})\leq \frac{K^2}{2} 	e^{-N/K }\left[ \mathcal L(x_0, \widehat{y})-\mathcal L(\widehat{x}, y_0)\right],
		\end{equation}
		where $\widehat{\bm z}=(\widehat{x},\widehat{y})$ a saddle point to \cref{eq:minmax-L}.
	\end{thm}
	\begin{proof}
		Under \cref{assum:sym-ABPD-PGS}, we claim that $f$ is relatively $\mu_f$-convex w.r.t. $\phi$ with  $\sgf\geq0$ and $g$ is relatively $\mu_g$-convex w.r.t. $\psi$ with  $\sgg\geq0$. 	Therefore, following the proof of \cref{eq:key-nji}, by \cref{thm:rate-sym-ABPD-PGS,eq:ub-H0}, we have
		\begin{equation}\label{eq:key-pgs}
			\mathcal L(x_{n_{j,i}},\widehat{y})-	\mathcal L(\widehat{x},y_{n_{j,i}})\leq 	
			C_3\left(\frac{C_1 }{ i^2}+ \frac{C_2 }{i}\right) \left[ \mathcal L(x_{n_{j,0}}, \widehat{y})-\mathcal L(\widehat{x}, y_{n_{j,0}})\right],
		\end{equation}
		for all $1\leq i\leq K-1$, where
		\begin{equation}\label{eq:C-pgs}
			C_1:=\frac{128(L_f \beta_0+4L_g\gamma_0) }{\gamma_0\beta_0}
			,\quad C_2:= \frac{8\opnm{A} }{\sqrt{\gamma_0\beta_0}},\quad C_3:=
			\frac{C_0}{\mu_f\mu_g} .
		\end{equation}
		This yields that 
		\[
		\begin{aligned}
			{}&	\mathcal L(x_{n_{j+1,0}},\widehat{y})-	\mathcal L(\widehat{x},y_{n_{j+1,0}})=
			\mathcal L(x_{n_{j,K }},\widehat{y})-	\mathcal L(\widehat{x},y_{n_{j,K }})\\
			={}&
			\mathcal L(x_{n_{j,K -1}},\widehat{y})-	\mathcal L(\widehat{x},y_{n_{j,K -1}})\leq 		C_3\left(\frac{C_1 }{ (K-1)^2}+ \frac{C_2 }{K-1}\right)\left[ \mathcal L(x_{n_{j,0}}, \widehat{y})-\mathcal L(\widehat{x}, y_{n_{j,0}})\right].
		\end{aligned}
		\]
		Again, let us consider any $N = \bar{j}K +\bar{i}$ with $\bar{j}\in\mathbb N$ and $1\leq \bar{i}\leq K -1$. It follows that
		\[
		\begin{aligned}
			\mathcal L(x_{n_{\bar{j}+1,0}},\widehat{y})-	\mathcal L(\widehat{x},y_{n_{\bar{j}+1,0}})\leq 	{}&C_3\left(\frac{C_1 }{ (K-1)^2}+ \frac{C_2 }{K-1}\right)\mathcal L(x_{n_{\bar{j},0}},\widehat{y})-	\mathcal L(\widehat{x},y_{n_{\bar{j},0}})\\
			\leq{}& \cdots\leq \left(\frac{C_1C_3 }{ (K-1)^2}+ \frac{C_2C_3 }{K-1}\right)^{\bar{j}} \left[ \mathcal L(x_{0}, \widehat{y})-\mathcal L(\widehat{x}, y_{0})\right].
		\end{aligned}
		\]
		Thanks to \cref{eq:key-pgs}, we obtain
		\[
		\begin{aligned}
			{}&	\mathcal L(x_{N},\widehat{y})-	\mathcal L(\widehat{x},y_{N})=	\mathcal L(x_{n_{\bar{j}+1,\bar{i}}},\widehat{y})-	\mathcal L(\widehat{x},y_{n_{\bar{j}+1,\bar{i}}})\\
			\leq{}& 	C_3\left(\frac{C_1 }{\bar{i}^2}+ \frac{C_2 }{\bar{i}}\right) \left[\mathcal L(x_{n_{\bar{j}+1,0}},\widehat{y})-	\mathcal L(\widehat{x},y_{n_{\bar{j}+1,0}})\right]\\
			\leq{}& (C_1C_3+C_2C_3)\left[\mathcal L(x_{n_{\bar{j}+1,0}},\widehat{y})-	\mathcal L(\widehat{x},y_{n_{\bar{j}+1,0}})\right]\\
			\leq {}&(C_1C_3+C_2C_3) \left(\frac{C_1C_3 }{ (K-1)^2}+ \frac{C_2C_3 }{K-1}\right)^{\bar{j}} \left[ \mathcal L(x_{0}, \widehat{y})-\mathcal L(\widehat{x}, y_{0})\right].
		\end{aligned}
		\]
		Observing \cref{eq:K-pgs,eq:C-pgs}, it is clear that $K -1\geq 2eC_2C_3,\,(K-1)^2\geq 2eC_1C_3$ and $\bar{j}\geq N/K -1$, which implies that 
		\[
		(C_1C_3+C_2C_3) \left(\frac{C_1C_3 }{ (K-1)^2}+ \frac{C_2C_3 }{K-1}\right)^{\bar{j}}\leq \frac{K^2}{2} e^{-N/K }.
		\]
		This gives \cref{eq:linear-rate-pgs-restart} and finishes the proof of this theorem.
	\end{proof}
	\begin{rem}
		According to \cref{thm:conv-ABPD-PS-R,thm:conv-ABPD-PGS-R}, with the priori restart period  choices \cref{eq:K,eq:K-pgs}, \cref{algo:ABPD-PS-restart,algo:ABPD-PGS-restart} admit linear converge rates. Moreover, if $\gamma_0=\mu_f$ and $\beta_0=\mu_g$, then the corresponding iteration complexity bounds are optimal
		\[
		\mathcal O\left(\frac{\opnm{A}}{\sqrt{\mu_f\mu_g}}\snm{\ln\epsilon}\right),\quad \text{and }
		\mathcal O\left(\left(\frac{\opnm{A}}{\sqrt{\mu_f\mu_g}}+\sqrt{\frac{L_f}{\mu_f}}+\sqrt{\frac{L_g}{\mu_g}}\right)\snm{\ln\epsilon}\right).
		\] 
	\end{rem}
\section{Numerical Experiments}
\label{sec:numer} 
In this part, we provide some numerical examples. 
We compare our \cref{algo:ABPD-PGS} (ABPD-PGS) and \cref{algo:ABPD-PGS-restart} (ABPD-PGS(Restart)) with the following baseline methods:
\begin{itemize}
	\item Accelerated primal-dual (APD) method \cite[Algorithm 2]{chen_optimal_2014},
	\item Chambolle--Pock's (CP) method \cite[Algorithms 1, 4 and 5]{chambolle_ergodic_2016},
	\item New Accelerated Smoothed GAp ReDuction (ASGARD+) \cite[Algorithm 1]{tran-dinh_unified_2021},
	\item Inertial corrected primal-dual proximal splitting (IC-PDPS) \cite[Algorithm 4.1]{valkonen_inertial_2020}.
\end{itemize}
As mentioned in \cref{rem:ABPD-PGS}, \cref{algo:ABPD-PS} can be viewed a special case of \cref{algo:ABPD-PGS} when the objective have no smooth components. Therefore, we omit the test for \cref{algo:ABPD-PS}. For all examples, we set the fixed restart period $K = 50$ for \cref{algo:ABPD-PGS-restart}.  To measure the performance, we consider the KKT residual
\[
{\rm KKT}(x_k,y_k):= \nm{x-\prox_{f}(x_k-{A^*} y_k)}+\nm{y_k-\prox_{g}(y_k+Ax_k)},
\]
and for the composite case: $f=f_1+f_2,\,g=g_1+g_2$, we use a modified one
\[
{\rm KKT}(x_k,y_k):=\nm{x_k-\prox_{f_2}(x_k-\nabla f_1(x_k)-{A^*} y_k)}+\nm{y_k-\prox_{g_2}(y_k-\nabla g_1(y_k)+Ax_k)}.
\]

\subsection{Simplex constrained least square problem}

Let us first consider the simplex-constrained least squares	problem \cite[Section 7.2]{chambolle_ergodic_2016}
\begin{equation}\label{eq:simplex-ls}
	\min_{x\in\Delta_n} \frac{1}{2}\|Ax-b\|^2,
\end{equation}
where $A\in\mathbb{R}^{m\times n}$ and $b\in\mathbb{R}^m$. The standard saddle-point formulation of \cref{eq:simplex-ls} is given by
$$\min\limits_{x\in\R^n}\max\limits_{y\in\R^m}\,\mathcal L(x,y):= \delta_{\Delta_n}(x)+\langle Ax,y\rangle-\dual{b,y}-\frac{1}{2}\|y\|^2.$$
Since 
\[
\sup_{x\in\Delta_n} \frac{1}{2}\|Ax-b\|^2\leq \max_{1\leq j\leq m} \frac{1}{2}\|A_{:,j}-b\|^2:=R,
\]
the saddle point problem is also equivalent to 
\[
\min\limits_{x\in\R^n}\max\limits_{y\in\R^m}\,\mathcal L(x,y):= \delta_{\Delta_n}(x) +\langle Ax,y\rangle-\delta_{B_R}(y)-\dual{b,y}-\frac{1}{2}\|y\|^2,
\]
where $B_R:=\{y\in\R^m:\,\nm{y}\leq \sqrt{2R}\}$. We consider two settings with different prox-functions.
\begin{itemize}
	\item For the Euclidean case, set $\mathcal X =\mathcal X^*= (\R^n,\nm{\cdot}_2)$ with $\phi(x) = 1/2\nm{x}^2_2$  and $\mathcal Y=\mathcal Y^*=(\R^m,\nm{\cdot}_2)$ with $\psi(y)=1/2\nm{y}_2^2$. Then $\opnm{A}=\nm{A}_2$.
	\item For the non-Euclidean case, set $\mathcal X = (\R^n,\nm{\cdot}_1)$ with $\phi(x) = -\sum_{i=1}^{n}x_i\ln x_i$ and $\mathcal Y=(\R^m,\nm{\cdot}_2)$ with $\psi(y) = 1/2\nm{y}_2^2$. Then $\mathcal X^*=(\R^n,\nm{\cdot}_\infty)$ and $\mathcal Y^*=\mathcal Y$. Moreover, $	\opnm{A} = \sup_{x\in\R^n}\nm{Ax}_{2}/\nm{x}_1 = \max_{1\leq j\leq m}\nm{A_{:,j}}$.
\end{itemize}


Numerical results are summarized in \cref{tab:simplex_lsp_normal,tab:simplex_lsp_uniform}, with normally and uniformly generated $A$ and $b$, respectively. The stopping criterion is either $k\geq 10^5$ or ${\rm KKT}(x_k,y_k)\leq 10^{-6}$. If the maximum iteration $10^5$ is reached, then the KKT residual at the final step is reported. Our ABPD-PGS is not competitive with the other methods, but its restart variant is highly efficient and outperforms all of them. For CP and ABPD-PGS, the entropy distance indeed improves performance; for APD and the restart variant of ABPD-PGS, however, the reverse case is observed.
\begin{table}[H]
	\caption{Numerical results for the simplex-constrained least squares problem \cref{eq:simplex-ls} with standard normally distributed data.}
	\label{tab:simplex_lsp_normal}
	\small  
	\renewcommand{\arraystretch}{1.2}
	\setlength{\tabcolsep}{4pt}  
	\setlength{\aboverulesep}{0pt}
	\setlength{\belowrulesep}{0pt}
	\begin{tabular}{l|cc|cc|cc}
		\toprule
		\multicolumn{1}{c|}{\multirow{2}{*}{Algorithm}} & \multicolumn{2}{c|}{$m=n=300$} & \multicolumn{2}{c|}{$m=n=500$} & \multicolumn{2}{c}{$m=n=1000$} \\
		\cmidrule{2-7}
		& Euclidean & Entropy & Euclidean & Entropy & Euclidean & Entropy \\
		\midrule
		CP & 3.73e-3 & 5.10e-4 & 8.66e-3 & 1.11e-3 & 2.65e-2 & 3.50e-3 \\
		\hline
		APD & 18326 & 3.60e-6 & 17318 & 9.74e-6 & 17959 & 2.02e-5 \\
		\hline
		IC-PDPS & 2.76e-5 & -- & 2.82e-5 & -- & 2.89e-5 & -- \\
		\hline
		ASGRAD+ & 12299 & -- & 12555 & -- & 13565 & -- \\
		\hline
		ABPD-PGS & 3.97e-3 & 2.50e-3 & 4.93e-3 & 3.05e-3 & 6.63e-3 & 4.06e-3 \\
		\hline
		ABPD-PGS(Restart) & 2535 & \textbf{2029} & \textbf{3381} & 11824 & \textbf{5266} & 6060 \\
		\bottomrule
	\end{tabular}
\end{table} 

\begin{table}[H]
	\caption{Numerical results for the  simplex-constrained least squares	problem \cref{eq:simplex-ls} with uniformly distributed data.}
	\label{tab:simplex_lsp_uniform}
	\small  
	\renewcommand{\arraystretch}{1.2}
	\setlength{\tabcolsep}{4pt}  
	\setlength{\aboverulesep}{0pt}
	\setlength{\belowrulesep}{0pt}
	\begin{tabular}{l|cc|cc|cc}
		\toprule
		\multicolumn{1}{c|}{\multirow{2}{*}{Algorithm}} & \multicolumn{2}{c|}{$m=n=300$} & \multicolumn{2}{c|}{$m=n=500$} & \multicolumn{2}{c}{$m=n=1000$} \\
		\cmidrule{2-7}
		& Euclidean & Entropy & Euclidean & Entropy & Euclidean & Entropy \\
		\midrule
		CP & 5.06e-4 & 5.97e-5 & 1.34e-3 & 1.39e-4 & 4.57e-3 & 4.77e-4 \\
		\hline
		APD & 12875 & 6.43e-6 & 11975 & 87258 & 13766 & 7.42e-6 \\
		\hline
		IC-PDPS & 2.73e-5 & -- & 2.84e-5 & -- & 2.88e-5 & -- \\
		\hline
		ASGRAD+ & 9460 & -- & 8563 & -- & 10268 & -- \\
		\hline
		ABPD-PGS & 2.50e-3 & 1.52e-3 & 3.22e-3 & 1.91e-3 & 4.33e-3 & 2.57e-3 \\
		\hline
		ABPD-PGS(Restart) & \textbf{1424} & 4878 & 1836 & \textbf{1774} & \textbf{2678} & 3373 \\
		\bottomrule
	\end{tabular}
\end{table}

\subsection{Elastic net problem} 
Then let us focus on the elastic net problem  \cite[Section 7.3]{chambolle_ergodic_2016}
\begin{equation}\label{eq:enp}
	\min_{x\in B_1}\frac12\|Ax-b\|^2+\sigma\|x\|_1+\frac{\mu}2\|x\|^2,
\end{equation}
where $B_1:=\{x\in\R^n:\,\nm{x}\leq 1\}$ is the unit ball, $A\in\mathbb{R}^{m\times n}$ and $b\in\mathbb{R}^m$ are given data, and $\sigma,\mu>0$ are parameters. It can be rewritten as the following standard saddle-point problem:
$$\min_{x\in \R^n}\max_{y\in\R^m}\,\delta_{B_1}(x)+\sigma\|x\|_1+\frac{\mu}{2}\|x\|^2+\langle Ax,y\rangle-\frac12\|y\|^2-\dual{b,y}.$$
Since 
\[
\max_{x\in B_1}\frac12\|Ax-b\|^2\leq \frac{1}{2}\left(\nm{A}+\nm{b}\right)^2:=R,
\]
the saddle point problem is also equivalent to 
$$\min\limits_{x\in\R^n}\max\limits_{y\in\R^m}\,\mathcal L(x,y):= \delta_{B_1}(x)+\sigma\|x\|_1+\frac{\mu}{2}\|x\|^2+\langle Ax,y\rangle-\delta_{B_R}(y)-\dual{b,y}-\frac{1}{2}\|y\|^2,$$
where $B_R:=\{y\in\R^m:\,\nm{y}\leq \sqrt{2R}\}$. Hence, for this example, we consider only the Euclidean setting.
%

\begin{figure}[H]
	\centering 
	\includegraphics[width=12cm]{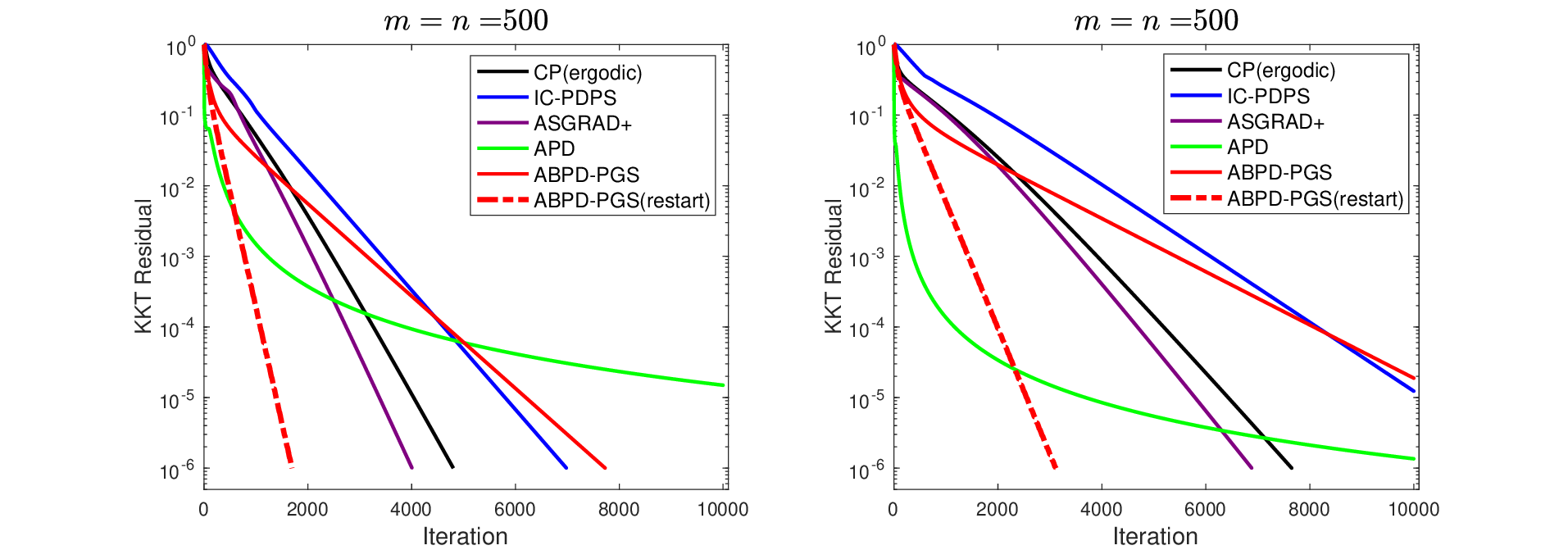}
	\caption{KKT residuals of all methods for the  elastic net problem \cref{eq:enp} with $\sigma=1,\mu=10^{-2}$. The data $A$ and $b$ are generated from the normal distribution (Left) and the uniform distribution (Right).}
	\label{fig:elasticnet}
\end{figure}

Numerical results are summarized in \cref{tab:elastic_net_normal,tab:elastic_net_uniform}, with $A$ and $b$ generated from normal and uniform distributions, respectively. The stopping criterion is either $k\geq 10^4$ or ${\rm KKT}(x_k,y_k)\leq 10^{-6}$. Again, if the maximum iteration $10^4$ is reached, then we report the KKT residual at the final step. It is observed that the restart variant of ABPD-PGS outperforms all other methods, with iteration counts less than half of those of ASGRAD+. In \cref{fig:elasticnet}, we also present the convergence behavior of the KKT residual. Since the problem under consideration is strongly convex-concave, all methods exhibit linear convergence except APD, which is designed for general convex-concave problems and thus achieves only a sublinear rate.
\begin{table}[H]
	\caption{Numerical results for the  elastic net problem \cref{eq:enp} with $\sigma=1,\mu=10^{-2}$ and standard normally distributed data.}
	\label{tab:elastic_net_normal}
	\renewcommand{\arraystretch}{1.3}
	\setlength{\aboverulesep}{0pt}
	\setlength{\belowrulesep}{0pt}
	\begin{tabular}{l|c|c|c}
		\toprule
		Algorithm & $m=n=300$ & $m=n=500$ & $m=n=1000$ \\
		\midrule
		CP & 6192 & 7901 & 2.13e-6 \\
		\hline
		APD & 3.97e-6 & 1.40e-6 & 9322 \\
		\hline
		IC-PDPS & 9500 & 1.71e-5 & 3.52e-4 \\
		\hline
		ASGRAD+ & 5392 & 7107 & 9787 \\
		\hline
		ABPD-PGS & 1.76e-6 & 2.78e-5 & 2.65e-4 \\
		\hline
		ABPD-PGS(Restart) & \textbf{2339} & \textbf{3236} & \textbf{4881} \\
		\bottomrule
	\end{tabular}
\end{table}

\begin{table}[H]
	\caption{Numerical results for the  elastic net problem \cref{eq:enp} with $\sigma=1,\mu=10^{-2}$ and uniformly distributed data.}
	\label{tab:elastic_net_uniform}
	\renewcommand{\arraystretch}{1.3}
	\setlength{\aboverulesep}{0pt}
	\setlength{\belowrulesep}{0pt}
	\begin{tabular}{l|c|c|c}
		\toprule
		Algorithm & $m=n=300$ & $m=n=500$ & $m=n=1000$ \\
		\midrule
		CP & 3974 & 4779 & 6487 \\
		\hline
		APD & 1.08e-4 & 1.68e-5 & 4.01e-6 \\
		\hline
		IC-PDPS & 5392 & 6939 & 1.09e-6 \\
		\hline
		ASGRAD+ & 3137 & 3988 & 5721 \\
		\hline
		ABPD-PGS & 5995 & 7677 & 3.18e-6 \\
		\hline
		ABPD-PGS(Restart) & \textbf{1330} & \textbf{1683} & \textbf{2491} \\
		\bottomrule
	\end{tabular}
\end{table}

\subsection{Total variation problem}
In this example, let us consider the total variation regularized linear inverse problem for image reconstruction:
\begin{equation}\label{eq:rof}
	\min_{x\in X}\,\frac\mu2\|x-b\|^2+\|Ax\|_{2,1},
\end{equation}
where $X:=\{x\in\mathbb{R}^N:x_l\leq x\leq x_u\}$ is a box constraint, $b\in\R^N$ with $N=mn$  represents the observed data from an $m$-by-$n$ picture and $A:\R^{N}\to\R^{2N}$ denotes the discrete gradient operator.
	Note that $\nm{A}\leq 2\sqrt{2}$.
	The problem \cref{eq:rof} can be reformulated as the following saddle point problem 
	$$
	\min_{x\in \R^{N}}\max_{y\in \R^{2N}}\,\mathcal L(x,y):=\delta_X(x)+\frac\mu2\|x-b\|^2+\langle Ax,y\rangle - \delta_Y(y).
	$$		 	
	%
	
\begin{figure}[H]
		\centering 
		\includegraphics[width=12cm]{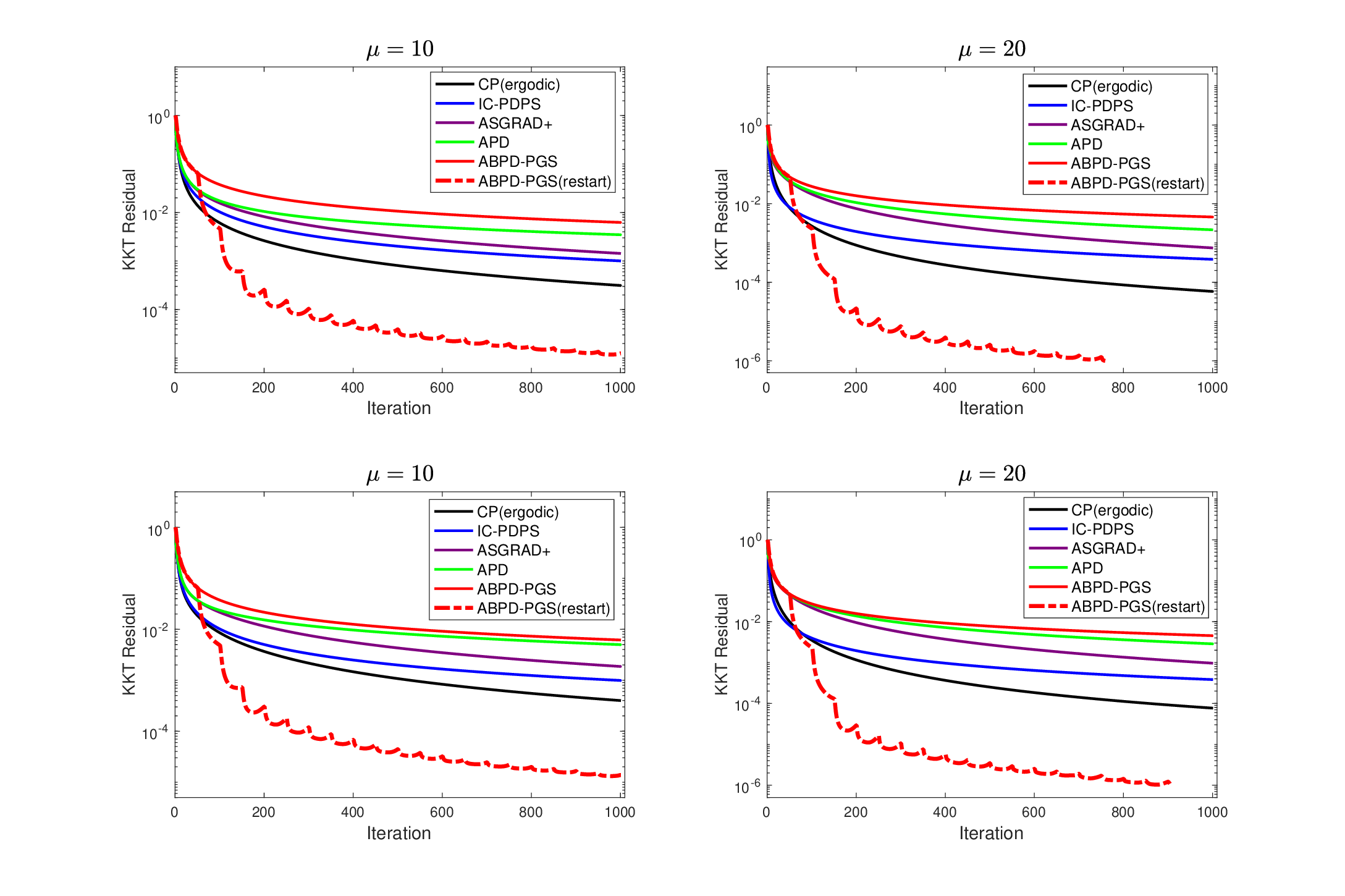}
		\caption{First row: \texttt{boat}. Second row: \texttt{peppers}}
		\label{fig:rof}
	\end{figure}

	We take two sample images of size $512 \times 512$: \texttt{boat} and \texttt{peppers}, which are widely used in the literature. The original images are corrupted by Gaussian noise with variance $0.1$. The stopping criterion is either $k\geq 10^3$ or ${\rm KKT}(x_k,y_k)\leq 10^{-6}$. Numerical results are given in \cref{tab:test-boat,tab:test-peppers}, reporting the number of iterations, the KKT residual at the final step, and the peak signal-to-noise ratio (PSNR). As the parameter $\mu$ increases, all methods perform better. The residuals of ABPD-PGS and APD are very close to each other and are larger than those of the other methods. Our restart variant of ABPD-PGS significantly outperforms all other methods, as shown in \cref{fig:rof}.

	\begin{table}[H]
		\caption{Numerical results for the image reconstruction problem \cref{eq:rof} with the original picture from \texttt{boat}}
		\label{tab:test-boat}
		\renewcommand{\arraystretch}{1.3}
		\setlength{\aboverulesep}{0pt}
		\setlength{\belowrulesep}{0pt}
		\begin{tabular}{l|cc|cc|cc}
			\toprule
			\multicolumn{1}{c|}{\multirow{2}{*}{Algorithm}} & \multicolumn{2}{c|}{$\mu=1$} & \multicolumn{2}{c|}{$\mu=10$} & \multicolumn{2}{c}{$\mu=20$} \\
			\cmidrule{2-7}
			& KKT & PSNR & KKT & PSNR & KKT & PSNR \\
			\midrule
			CP & 4.02e-3 & 20.7483 & 3.19e-4 & 28.2797 & 5.77e-5 & 27.3258 \\
			\hline
			APD & 6.43e-3 & 20.7793 & 3.49e-3 & 28.2837 & 2.18e-3 & 27.3248 \\
			\hline
			IC-PDPS & 1.46e-2 & 20.7039 & 1.00e-3 & 28.2780 & 3.83e-4 & 27.3220 \\
			\hline
			ASGRAD+ & 4.78e-3 & 20.7499 & 1.44e-3 & 28.2799 & 7.51e-4 & 27.3258 \\
			\hline
			ABPD-PGS & 2.10e-2 & 20.7233 & 6.24e-3 & 28.2281 & 4.57e-3 & 27.2969 \\
			\hline
			ABPD-PGS(Restart) & \textbf{7.29e-4} & 20.7498 & \textbf{1.28e-5} & 28.2798 & \textbf{774} & 27.3258 \\
			\bottomrule
		\end{tabular}
	\end{table}
	
	\begin{table}[H]
		\caption{Numerical results for the image reconstruction problem \cref{eq:rof} with the original picture from \texttt{peppers}}
		\label{tab:test-peppers}
		\renewcommand{\arraystretch}{1.3}
		\setlength{\aboverulesep}{0pt}
		\setlength{\belowrulesep}{0pt}
		\begin{tabular}{l|cc|cc|cc}
			\toprule
			\multicolumn{1}{c|}{\multirow{2}{*}{Algorithm}} & \multicolumn{2}{c|}{$\mu=1$} & \multicolumn{2}{c|}{$\mu=10$} & \multicolumn{2}{c}{$\mu=20$} \\
			\cmidrule{2-7}
			& KKT & PSNR & KKT & PSNR & KKT & PSNR \\
			\midrule
			CP & 6.08e-3 & 21.3838 & 3.95e-4 & 30.2346 & 7.60e-5 & 28.0272 \\
			\hline
			APD & 1.09e-2 & 21.4007 & 4.95e-3 & 30.2512 & 2.85e-3 & 28.0269 \\
			\hline
			IC-PDPS & 1.40e-2 & 21.3156 & 9.91e-4 & 30.2292 & 3.83e-4 & 28.0224 \\
			\hline
			ASGRAD+ & 8.04e-3 & 21.3857 & 1.84e-3 & 30.2351 & 9.59e-4 & 28.0273 \\
			\hline
			ABPD-PGS & 2.07e-2 & 21.3685 & 6.17e-3 & 30.1591 & 4.53e-3 & 27.9984 \\
			\hline
			ABPD-PGS(Restart) & \textbf{8.93e-4} & 21.3785 & \textbf{1.39e-5} & 30.2349 & \textbf{909} & 28.0272 \\
			\bottomrule
		\end{tabular}
	\end{table} 
	
\section{Conclusion}
\label{sec:conclu}
In this work, we start with the spectral analysis of a quadratic game and show that a {\it naive} application of the unconstrained model \cref{eq:nag-sys} to saddle-point problems does not work. We then propose a new template model \cref{eq:APDG-spec} and prove its stability and improved spectral radius that matches the theoretical optimal convergence rate. After that, we combine general Bregman distances with proper time rescaling parameters and present two novel continuous models called APDG flow \cref{eq:APDG} and ABPDG flow \cref{eq:ABPDG-sys}, for which the exponential convergence of the Lagrange gap along the continuous trajectory is established via a compact Lyapunov analysis. 

Based on the ABPDG flow, we consider two semi-implicit Euler discretizations and obtain two first-order primal-dual splitting methods for solving the convex-concave saddle-point problem \cref{eq:minmax-L}, as summarized in \cref{algo:ABPD-PS,algo:ABPD-PGS}. By using careful discrete Lyapunov analysis, we establish the optimal convergence rate of the Lagrange gap with sharp dependence on the Lipschitz and convexity constants. This also matches the theoretical lower iteration complexity bound of first-order primal-dual methods. Also, efficient restart variants are proposed and analyzed. Numerical results validate the performances of the restart methods. 

\appendix
\section{Proof of \cref{thm:spec-gs}}
\label{app:spec-gs}
A key property is that $\mathcal G_{\rm new}$ can be transferred into a block tridiagonal matrix
\[
\mathcal P\mathcal G_{\rm new} \mathcal P^\top= 
\begin{bmatrix}
	-I&I&O&O\\
	O&-I&A/\mu&O\\	
	O&-A^\top/\mu&-I&O\\
	O& O&I&-I	
\end{bmatrix},\quad \text{with }\,
\mathcal P = 
\begin{bmatrix}
	O&I&O&O\\
	O&O&O&I\\	
	O&O &I&O\\
	I& O&O&O	
\end{bmatrix}.
\]
Note that $\mathcal P$ is a permutation matrix and $\mathcal G_{\rm new} = -I +\mathcal L + \mathcal U$, where $\mathcal L$ and $\mathcal U$ are respectively the strictly lower triangular part and upper triangular part. According to \cite[Chapter 4.1]{varga_matrix_2000}, we find that $\mathcal L+\mathcal U$ is consistently ordered, weakly cyclic of index 2, which yields the spectrum relation
\[
\sigma(\mathcal L+\mathcal U) = \sigma\left(\eta\mathcal L+\mathcal U/\eta\right),\quad\forall\,\eta\in\mathbb C\backslash\{0\}.
\]
Moreover, by \cite[Theorem 4.4]{varga_matrix_2000}, we have 
\begin{equation}\label{eq:vargathm4.4}
	\det\left(\gamma I-\eta\mathcal L-\beta\mathcal U\right)
	=\det\left(\gamma I-\sqrt{\eta\beta}\left(\mathcal L+\mathcal U\right)\right),\quad\forall\,\gamma,\eta,\beta\in\mathbb C.
\end{equation}

We then compute the eigenvalues of $\widehat{\mathcal G} = \left[(1+\alpha)I-\alpha \mathcal L\right]^{-1}\left(I+\alpha \mathcal U\right)$. For any $\tau\in\sigma(\widehat{\mathcal G})$, it holds that 
\[
0 = \det\left(\tau I-\widehat{\mathcal G}\right) = \det\left[((1+\alpha)\tau-1)I-\alpha\tau\mathcal L-\alpha\mathcal U\right].
\]
In view of \cref{eq:vargathm4.4}, we obtain 
\[
\begin{aligned}
	0 ={}& \det\left[\left((1+\alpha)\tau-1\right)I-\alpha\sqrt{\tau}\left(\mathcal L+\mathcal U\right)\right] 
	={} 
	\det\left[\left((1+\alpha)\tau-\alpha\sqrt{\tau}-1\right)I-\alpha\sqrt{\tau}\mathcal G_{\rm new}\right].
\end{aligned}
\]
From this we claim that $\tau\neq 0$. Then it follows
\begin{equation}\label{eq:eig-H}
	(1+\alpha)\tau - \alpha(1+\lambda)\sqrt{\tau}-1= 0,\quad \lambda\in\sigma(\mathcal G_{\rm new}).
\end{equation}
By \cref{eq:eig-Gnew}, we have either $\lambda = -1$ or $\lambda=-1\pm i\sqrt{\delta}/\mu$, with $\delta\in\sigma(A A^\top  )$. If $\lambda=-1$, then $\tau = (1+\alpha)^{-1}$. Otherwise, \cref{eq:eig-H} becomes 
\[
(1+\alpha)\tau \pm i\frac{\alpha\sqrt{\delta}}{\mu}\sqrt{\tau}-1= 0\quad \Longleftrightarrow\quad  
(1+\alpha)\tau -1=\pm i\frac{\alpha\sqrt{\delta}}{\mu}\sqrt{\tau}.
\] This further implies that $\tau\in\sigma(\widehat{\mathcal G})$ satisfies the quadratic equation (with real coefficients)
\begin{equation}\label{eq:eig-H-}
	(1+\alpha)^2\tau^2 + (\delta\alpha^2/\mu^2-2(1+\alpha))\tau+1=0.
\end{equation}
Since $\alpha\leq 2\mu/\nm{A}\leq 2\mu/\sqrt{\delta}$ for any $\delta\in\sigma(A A^\top )$, we find that
\[
\Delta =(\delta\alpha^2/\mu^2-2(1+\alpha))^2-4(1+\alpha)^2= \frac{\delta^2\alpha^4}{\mu^4}-\frac{4\delta}{\mu^2}(1+\alpha)\alpha^2 = \frac{\delta\alpha^2}{\mu^4}\left(\delta\alpha^2-4\mu^2(1+\alpha)\right)< 0.
\]
Hence, the quadratic equation \cref{eq:eig-H-}  has two conjugate complex roots and thus by Vieta's theorem, we get $\snm{\tau}^2 = (1+\alpha)^{-2}$. Therefore, this implies $\rho(\widehat{\mathcal G})=(1+\alpha)^{-1}$ and completes the proof of \cref{thm:spec-gs}.
\section{A Technical Estimate}
We list a useful estimate of the positive sequence $\{\theta_k\}$ that satisfy
\begin{equation}\label{eq:tk-case2}
	\theta_{k+1}-\theta_k \leq  -\frac{ \theta^{\nu}_{k}\theta_{k+1}}{\sqrt{Q\theta_{k}+R^2}},\quad \theta_0=1,
\end{equation}
where $ Q>0,\,R\geq 0 $ and $\nu\geq 1/2$. The detailed proof can be found in \cite[Lemma C.2]{Luo2025}.
\begin{lem}[\cite{Luo2025}]
	\label{lem:app-est}
	Let $\{\theta_k\}$ be a positive real sequence satisfying \cref{eq:tk-case2}. If  $\theta_{k+1}/\theta_{k}\geq \delta>0$ for $k\in\mathbb N$, then for all $k\geq 0$, it holds that
	\begin{equation*} 
		\theta_k
		\leq \left\{
		\begin{aligned}
			{}&\exp\left(-\frac{\delta k}{2\sqrt{Q}} \right)+
			\left(	\frac{4R}{4R+\delta k}\right)^{2},&&\text{if }\nu=1/2,\\	
			{}&	\left(\frac{4\sqrt{Q}}{4\sqrt{Q}+\delta(2\nu-1) k}\right)^{\frac{2}{2\nu-1}}+
			\left(\frac{2R}{2R+\delta\nu k}\right)^{\frac{1}{\nu}}			,&&\text{if }\nu>1/2.
		\end{aligned}
		\right.
	\end{equation*}
\end{lem}

\section{Ergodic v.s. nonergodic: An illustrative example}
\label{sec:counter}
\renewcommand{\arraystretch}{1.1}
In this part, we discuss the iterative behavior of the nonergodic and ergodic sequences of CP. Let $A\in\R^{n\times n}$ be given and consider the matrix game 
\begin{equation}\label{eq:counter}
	\min_{x\in \R^n}\max_{y\in \R^n}\,y^\top Ax=\frac{1}{2}(x^\top,-y^\top) \mathcal M\begin{pmatrix}
		x\\ y
	\end{pmatrix},\quad \text{where }\,\mathcal M = 	 \begin{bmatrix}
		O&A^\top \\ -A&O
	\end{bmatrix}. 
\end{equation} 
The set of all saddle points is the null space $\mathcal M^{-1}(0) = \{(x^*,y^*):\,x^*\in{\rm Null}(A),\,y^*\in{\rm Null}(A^\top)\}$.
Applying CP \cref{eq:cp}
with $\alpha=1$ to this problem leads to the iteration scheme
\begin{equation}\label{eq:cp-counter}
	\left\{
	\begin{aligned}
		x_{k+1} = {}&x_k-\tau A^\top y_k, \\ 
		y_{k+1} ={}&y_k+\theta A(2x_{k+1}-x_k) ,
	\end{aligned}
	\right.
\end{equation}
where $\tau\theta\nm{A}^2<1$. According to \cite[Theorem 1]{chambolle_first-order_2011}, for {\it any} initial guesses $x_0,\,y_0\in\R^n$, we have the $O(1/k)$ rate for the {\it ergodic sequence} $(X_k,Y_k):=\frac{1}{k}\sum_{i=1}^{k}(x_i,y_i)$. In the following, by using the spectral analysis, we confirm that the nonergodic sequence $\{(x_k,y_k)\}$ is convergent indeed for any invertible $A$. However, we will provide a simple example (cf.\cref{eq:counter-A}) to show the dramatically slow convergence of  $\{(x_k,y_k)\}$.

Introduce $z_k=(x_k,y_k)$ and rewrite \cref{eq:cp-counter} as follows
\begin{equation}\label{eq:zk}
	z_{k+1} = \mathcal B_{\rm cp}z_k,\qquad \text{ with } \,\mathcal B_{\rm cp}:=\begin{bmatrix}
		I&-\tau A^\top\\\theta A&I-2\tau\theta AA^\top
	\end{bmatrix}.
\end{equation}
We claim that $\sigma(\mathcal B_{\rm cp}) = \{\lambda(\xi) = 1-\tau\theta\xi\pm i\sqrt{\tau\theta\xi(1-\tau\theta\xi)}:\,\xi\in\sigma(AA^\top)\}$. Indeed, for any $\lambda\in\sigma(\mathcal B_{\rm cp})$, if $\lambda=1$, then $0=\det(I-\mathcal B_{\rm cp})=\snm{\det A}^2\tau^n\theta^n$, yielding $0\in\sigma(AA^\top)$ and $\lambda=\lambda(0)$. If $\lambda\neq 1$, then we have
\[
\begin{aligned}
	0={}&\det(\lambda I-\mathcal B_{\rm cp})=\det\begin{bmatrix}
		(\lambda-1)I&\tau A^\top\\-\theta A&(\lambda-1)I+2\tau\theta AA^\top
	\end{bmatrix} \\
	={} &\det\left((\lambda-1)^2I+\tau
	\theta(2\lambda-1)AA^\top\right),
\end{aligned}
\]
which implies that $(\lambda-1)^2-\tau\theta\xi(1-2\lambda)=0$, with $\xi\in\sigma(AA^\top)$. From this and the fact $\tau\theta\xi\leq \tau\theta\nm{A}^2<1$, we obtain $\lambda = \lambda(\xi)$. According to the above discussion, we conclude that $\rho(\mathcal B_{\rm cp}) = \sqrt{1-\tau\theta\lambda_{\min}(AA^\top)}$, where $\lambda_{\min}(AA^\top)$ is the smallest eigenvalue of $AA^\top$. 
If $A$ is invertible, then $\rho(\mathcal B_{\rm cp})<1$ and \cref{eq:zk} is always convergent for any initial setting $z_0$.
However, if the smallest eigenvalue $\lambda_{\min}(AA^\top)=O(\epsilon)$ is close to 0, then $\rho(\mathcal B_{\rm cp})  = \sqrt{1-O(\epsilon)}$ and this might lead to very slow convergence rate.

\begin{figure}[H]
	\centering 
	\includegraphics[width=12cm]{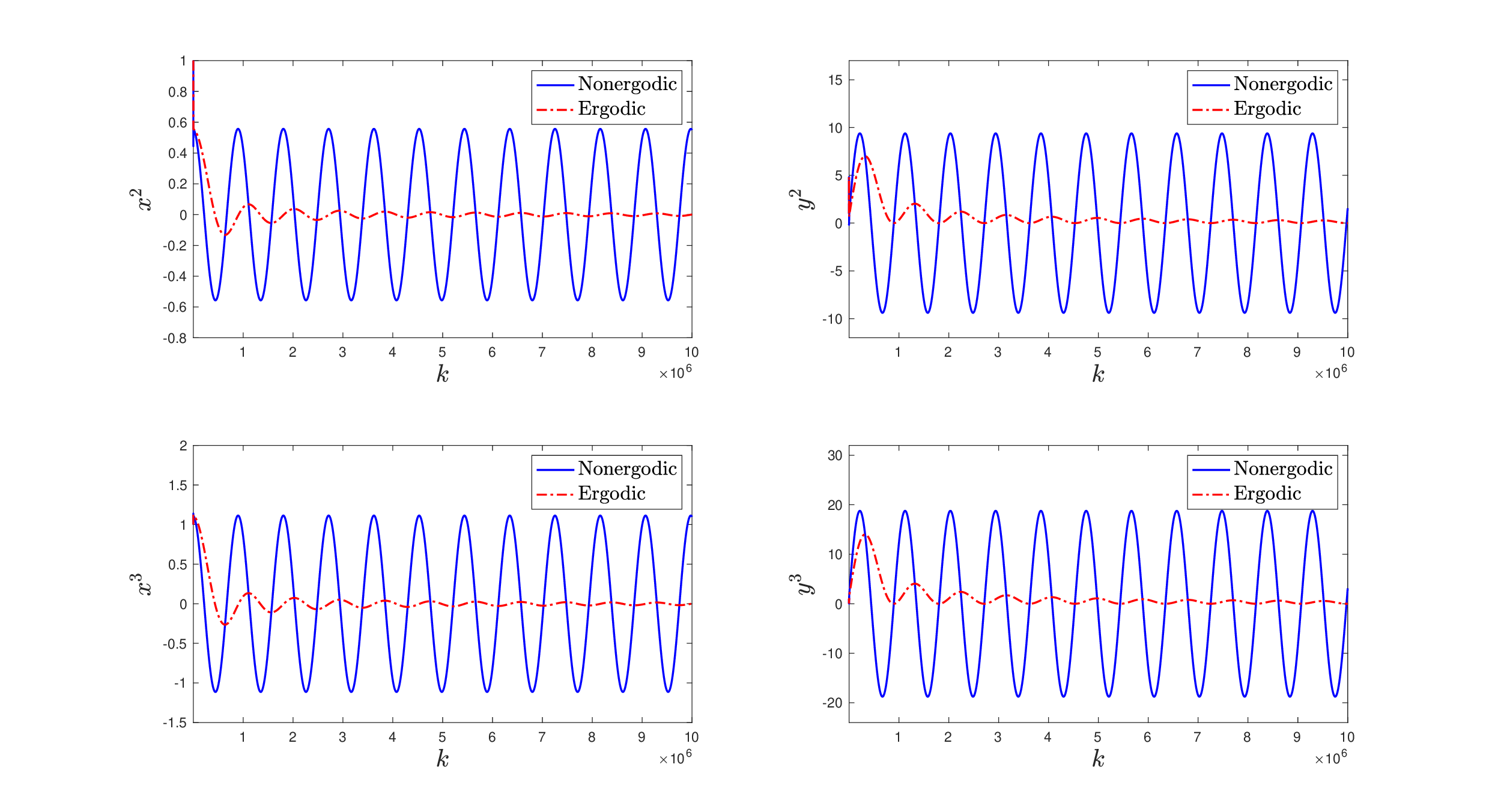}
	\caption{Iterative cycle behavior of CP for solving \cref{eq:counter} with \cref{eq:counter-A} and $\epsilon=10^{-4}$. The initial points are $x_0=y_0=[1,1,1]^\top$ and the parameters are $\alpha=1,\,\theta=1,\,\tau = 0.5/\nm{A}^2$. The first components of the nonergodic sequence $\{(x_k,y_k)\}$ and the ergodic sequence $\{(X_k,Y_k)\}$ have already been reported in \cref{fig:cp-limit-cycle}. In the current figure, the left column and the right column are respectively the rest components of $\{(x_k,X_k)\}$ and $\{(y_k,Y_k)\}$.}
	\label{fig:cp-limit-cycle-x-y}
\end{figure}

\begin{figure}[H]
	\centering  
	\includegraphics[width=14cm]{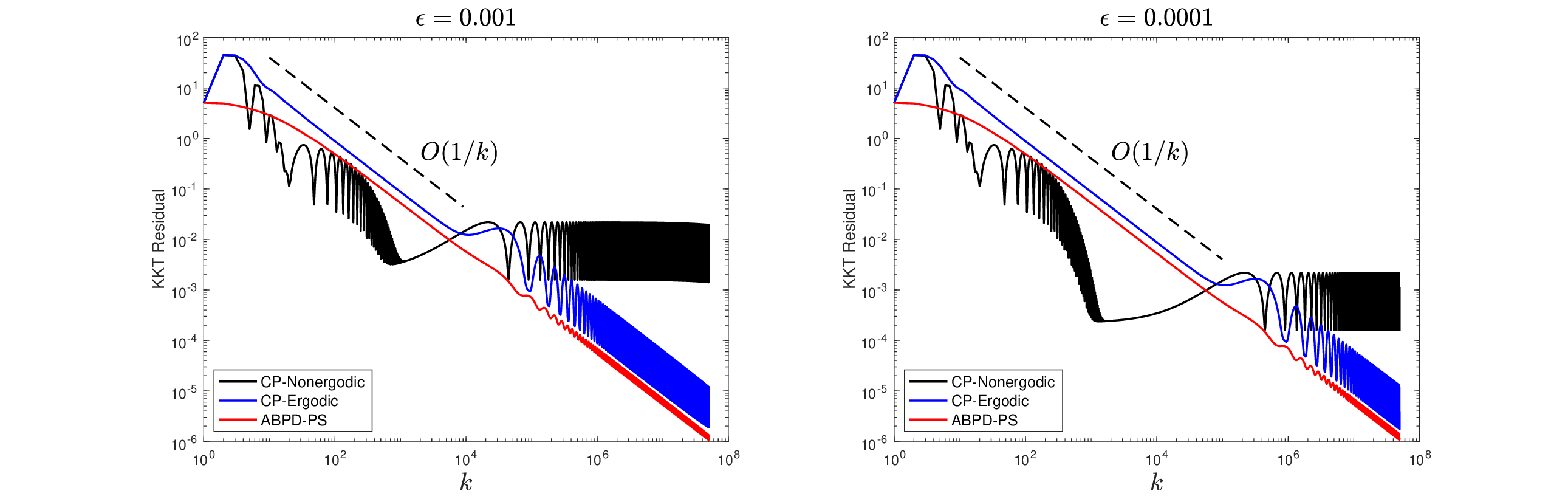} 
	\caption{KKT residuals of CP and ABPD-PS for solving \cref{eq:counter} with \cref{eq:counter-A} and different $\epsilon$. The initial points are $x_0=y_0=[1,1,1]^\top$ and the parameters for CP are $\alpha=1,\,\theta=1,\,\tau = 0.5/\nm{A}^2$. 
	}
	\label{fig:cp-limit-cycle-res}
\end{figure}

To further support our claim, we provide a simple example 
\begin{equation}\label{eq:counter-A}
	A = \epsilon I + A_0,\quad A_0 =\begin{bmatrix}
		2 & -4 & 0  \\
		-2 & 8 & -2 \\
		0 & -2 & 1 \\
	\end{bmatrix}.
\end{equation}
Note that $A$ is invertible for all $\epsilon>0$ and $ \lambda_{\min}(AA^\top)= O(\epsilon^2),\,\nm{A}\leq 13+\epsilon$,
which implies
\[
\rho(\mathcal B_{\rm cp})= \sqrt{1-\tau\theta O(\epsilon^2)}\geq  \sqrt{1-O(\epsilon^2)/\nm{A}^2}=1-O(\epsilon^2).
\]
In \cref{fig:cp-limit-cycle-x-y}, we report the iterative behavior of CP for solving \cref{eq:counter} with \cref{eq:counter-A} and $\epsilon=10^{-4}$. The initial points are $x_0=y_0=[1,1,1]^\top$ and the parameters are $\alpha=1,\,\theta=1,\,\tau = 0.5/\nm{A}^2$. Moreover, in \cref{fig:cp-limit-cycle-res}, we plot the historical KKT residuals of the ergodic sequence $\{(X_k,Y_k)\}$ and the nonergodic sequence $\{(x_k,y_k)\}$.

It is observed from \cref{fig:cp-limit-cycle-x-y} that the nonergodic sequence $\{(x_k,y_k)\}$ is highly oscillated around the unique saddle point $0$ and is extremely close to a limit cycle. According to \cref{fig:cp-limit-cycle-res}, the residuals corresponding to ergodic sequences decrease like the theoretical rate $O(1/k)$. On the contrary, the residuals of the nonergodic sequence decay very quickly within the first few steps but then remain flat in an oscillatory manner for a very long period. As a comparison, we also apply \cref{algo:ABPD-PS} to \cref{eq:counter} and in \cref{fig:cp-limit-cycle-res}, we report the KKT residual which agrees with the theoretical prediction $O(1/k)$.
%

\bibliographystyle{abbrv}


\end{document}